\documentclass{elsart}

\usepackage{amssymb}
\usepackage{color}
\usepackage{graphicx}
\usepackage{subfigure}
\usepackage{amsmath}
\usepackage{amsfonts}
\usepackage[title]{appendix}
\usepackage[colorlinks=true,citecolor=red,linkcolor=red]{hyperref}
\usepackage{multirow}

\def\dfr#1#2{\displaystyle{\frac{#1}{#2}}}
\renewcommand{\vec}[1]{\mbox{\boldmath \small $#1$}}
\newcommand{\smallvec}[1]{\mbox{\boldmath \scriptsize $#1$}}

\newtheorem{example}{Example}[section]

\newtheorem{remark}{Remark}[section]
\numberwithin{equation}{section}
\numberwithin{figure}{section}
\numberwithin{table}{section}

\begin{document}
%\journal{Applied Mathematics and Mechanics (English Edition)}

\begin{frontmatter}

\title{A Newton multigrid method for  steady-state shallow water equations
with topography and dry areas}
  \author{Kailiang Wu},
\ead{wukl@pku.edu.cn}
\author[label2]{Huazhong Tang}
\thanks[label2]{Corresponding author. Tel:~+86-10-62757018;
Fax:~+86-10-62751801.}
\ead{hztang@math.pku.edu.cn}
\address{HEDPS, CAPT \& LMAM, School of Mathematical Sciences, Peking University,
Beijing 100871, P.R. China}
 \date{\today{}}

\maketitle

\begin{abstract}
The paper develops a Newton multigrid  ({MG}) method for one- and two-dimensional steady-state shallow water equations (SWEs) with topography and dry areas. It solves the nonlinear system arising from the well-balanced finite volume discretization of the steady-state SWEs by using Newton's method as the outer iteration and a geometric MG method  with the  block symmetric Gauss-Seidel smoother as the inner iteration.
%Instead of adopting the time-stepping relaxation technique based on
%the local CFL number,
The proposed {Newton MG}  method makes use of the local residual to regularize the Jacobian matrix of the Newton iteration,  and can  handle the steady-state problem with wet/dry transitions.
Several numerical experiments are conducted to demonstrate the
efficiency, robustness, and  well-balanced property  of the proposed  method.
The relation between the convergence behavior of the {Newton MG}  method and the distribution
of the eigenvalues of the iteration matrix is detailedly discussed. % by considering different numerical fluxes.
\end{abstract}


\end{frontmatter}

%%%%%%%%%%%%%%%%%%%%%%%%%%%%%%%%%%%%%%%%%%%%%%%%%%%%

\section{Introduction}
\label{sec:intro}

%The mathematical model of fluid flow  plays an important role in a large variety of natural or man-made situations.
%Originated by Saint-Venant \cite{BarredeSaint1987} in the study of
%floods and tides,
The shallow water equations ({SWEs}) are commonly used to describe the motion of ``shallow'' free-surface flows
 subject to gravitational force and have played
a critical role in the modeling and simulation of the flows
in rivers or channels, the ocean tides, and the tsunami, etc.
Under the assumption of incompressible
fluid and hydrostatic pressure distribution, with the vertical acceleration of water particles
neglected,
the SWEs may be derived by depth-integrating
the %Reynolds-averaged 
Navier-Stokes equations as follows  \cite{Vreugdenhil1994}
\begin{equation}\label{eq:GNeqs}
\frac{{\partial \vec U}}{{\partial t}} + \nabla  \cdot \vec F(\vec U) = \vec S(\vec x,\vec U),
\end{equation}
where $t$ denotes time, and  the effect of bed slope on the flow has been modeled by the inclusion of source terms at the right hand
side of  \eqref{eq:GNeqs}  which modifies the momentum   equations.
In the 2D Cartesian coordinates $\vec x =(x,y)$,  the conservative vector $\vec U$,
the flux vector $\vec F=(\vec F_1, \vec F_2)$, and the source term $\vec S$ in \eqref{eq:GNeqs} are given by
%\[
%\vec U = \left( {\begin{array}{*{20}c}
%	h  \\
%	{hu}
%	\end{array}} \right),
%\quad
%\vec F = \left( {\begin{array}{*{20}c}
%	{hu}  \\
%	{hu^2  + \frac{1}{2}gh^2 }
%	\end{array}} \right),
%\quad
%\vec S = \left( {\begin{array}{*{20}c}
%	0  \\
%	{ - gh\frac{{\partial B}}{{\partial x}}}
%	\end{array}} \right),
%\]
\[
\vec U = \left( {\begin{array}{*{20}c}
	h  \\
	{h\vec u^{\rm T} } \\
	\end{array}} \right),
\quad
\vec F = \left( {\begin{array}{*{20}c}
	{h\vec u}  \\
	{h\vec u^{\rm T} \vec u  + \frac{1}{2}gh^2\vec I }
	\end{array}} \right),
%\quad
%\vec F_2 = \left( {\begin{array}{*{20}c}
%	{hv}  \\
%	{huv} \\
%	{hv^2  + \frac{1}{2}gh^2 }
%	\end{array}} \right),
\quad
\vec S = \left( {\begin{array}{*{20}c}
	0  \\
	{ - gh\nabla  z} \\
	\end{array}} \right),
\]
in which $\vec u=(u,v)$ denotes the velocity vector, $g$ is the acceleration due to gravity, $\vec I$ is the identity matrix,
and $h$ is the water depth, i.e. the height of  water above the
 riverbed topography $z(\vec x)$.
Dropping the time derivatives in   \eqref{eq:GNeqs} gives the  the steady-state SWEs
\begin{equation}\label{eq:GNeqsSTEEADY}
\nabla  \cdot \vec F (\vec U) = \vec S(\vec x,\vec U).
\end{equation}
If the water is at rest, i.e. $\vec u(\vec x)=\vec 0$,
then the momentum parts in the above equations reduce to
 \begin{equation}\label{eq:GNeqsSTEEADY2}
\frac12(gh^2)_x=-gh z_x, \quad \frac12(gh^2)_y=-gh z_y,
 \end{equation}
which imply the so-called well-balanced property,
i.e. ``steady state of the water at rest''
$$
\vec u(\vec x)=\vec 0,  \quad h(\vec x)+z(\vec x)=\mbox{const}.
$$

% The  shallow water equations (SWEs) have been widely applied in ocean and hydraulic engineering
%to describe hydrodynamic phenomena such as the hydraulic jump, bore wave propagation, surface
%irrigation, dam break, tidal flows in estuary and coastal water regions, and open channel flows,
% and surface runoff.
%
%The shallow water equations are derived from equations of the mass and linear momentum conservation, which hold even when the assumptions of shallow water break down, such as across a hydraulic jump.
%The equations are derived from depth-integrating the Navier-Stokes equations, in the case where the horizontal length scale is much greater than the vertical length scale.

%Generally, the  analytical solutions to the SWEs cannot be obtained.
% are often established on the basis of
%simple geometries and mostly suffer from simplified underlying assumptions. %e.g. \cite{Delestre2013},
%The paper is concerned with numerical methods for  the quasilinear SWEs.
% with topography and dry areas.
%Currently, the primary approach for the study of the quasilinear SWEs is numerical simulation.
 Up to now, there exist various numerical methods  for the SWEs,
 such as the finite difference scheme based on flux-difference splitting \cite{Glaister1988},
 the generalized Riemann problem scheme \cite{LiChen2006},
high-order WENO schemes \cite{Noellea2007,Xing2005},  the gas-kinetic schemes \cite{TangTangXu2004,Xu2002},
the moving mesh method \cite{Tang2004}, and
the Runge-Kutta discontinuous Galerkin methods \cite{Kesserwani2011,Xing2010}, and so on.
Most of them  are explicit for the time-dependent SWEs,
and can successfully simulate the evolution of the time-dependent solutions
with good accuracy in time.
%In practice,   it is also interesting and important to
%investigate the steady-state behavior of the flow.
If using the explicit time advancing method to
investigate the steady-state behavior of the flow,
then the long time simulation is needed and becomes time-consuming due to the small time step size
satisfying a Courant-Friedrichs-Lewy  condition to guarantee stability.
For example, in  the sediment transport and morphodynamic change model \cite{Bilanceri2012,Hudson2001},  %coupling the SWEs  with a sediment transport equation,
 %the characteristic time scales of the flow and of the sediment transport may be very
%different  and
the time stiffness arising from the characteristic time scales in flow
and sediment transport seriously challenges the  long time simulation
if the interaction of  water flow with  bed topography is very weak.
%A scheme with large time step is excepted to fast obtain the ``steady-state'' solution of the
%SWE in each step of fixed bed topography.
%It is well known that for explicit time advancing the time step is limited by numerical stability,
%and thus it usually takes expensive time cost for a explicit scheme to capture %the steady-state flow or
%the long time steady-state behavior of a flow.
%As a result of the requirements for practical applications,
For such case, the implicit or semi-implicit schemes are attractive,
such as the multigrid semi-implicit finite difference method \cite{Spitaleri1997},
the
%{\color{red} sign matrix based [it was not mentioned specifically in the paper]}
linearized implicit scheme with a modified Roe flux \cite{Bilanceri2012},
the space-time discontinuous residual distribution scheme \cite{Sarmany2013},
%the implicit finite volume method developed in \cite{Wu2011} on unstructured quadtree rectangular mesh,
the implicit higher-order compact scheme \cite{Bagheri2013}, and
the semi-implicit discontinuous Galerkin methods \cite{Dumbser2013,Tavellia2014} etc.
In fact, the time step size for an implicit scheme is also often constrained by convergence.
Even for the same time step size as used in the explicit case, unsteady
solutions of the implicit schemes may be less accurate.
An implicit scheme usually requires solving a nonlinear equation by some iteration method,
and  thus it is also very time-consuming.
For the steady-state behavior of the {SWEs} flow,
another way is to directly solve the steady-state SWEs,
 see e.g. \cite{Glaister1994}, and thus
 developing robust and efficient solver for corresponding nonlinear algebraic system is  key.
Unfortunately there are few such study for the steady-state SWEs,
but  the steady-state Euler equations and Navier-Stokes equations have been
 well solved numerically, see e.g. \cite{Chen2000,HuLiTang2010,Jameson1995,LiWang2008,Mavriplis1989}.
For example, a multigrid block lower-upper
symmetric Gauss-Seidel (LU-SGS) algorithm was proposed for the 2D steady-state Euler equations on unstructured grid \cite{LiWang2008}.  Unlike the existing methods which add the pseudo-time terms to the steady-state equations, the norm of the local residual in each cell was used  to  regularize the nonlinear algebraic system arising from the  spatial discretization of the steady-state Euler equations. The Newton iteration was then adopted to solve
the nonlinear algebraic system and the multigrid method was used as the inner iteration with the BLU-SGS
 smoother.
%The algorithm has been extended to high-order in \cite{HuLiTang2010} by using the quadratic reconstruction with the WENO hierarchical limiting strategy \cite{Xu2009}.
%In addition to the efficiency of the solvers for those nonlinear systems, the robustness of a scheme for
%the steady-state SWEs should also be taken into consideration.

The aim of the paper is to extend  the Newton multigrid method  \cite{LiWang2008} to the steady-state SWEs
and investigate its convergence,
in which the steady-state SWEs are
discretized by a well-balanced hydrostatic reconstruction, % \cite{Audusse2004},
 the wet/dry transition is numerically handled,
 %for the topography source term
 %
%especially the one done by Li et al. in \cite{LiWang2008}.
%
%In order to be consistent with the lake at rest,
%a suitable discretization technique well-known as   will be
the resulting nonlinear algebraic  system
 is iteratively solved by using the Newton  multigrid  iteration,
%By considering different numerical fluxes for numerical experiments, it will be
%observed that using different numerical fluxes shows
and  convergence behavior of the
 method for different numerical fluxes is  detailedly investigated in numerical experiments
%which is an interesting phenomenon and
%may be explained
through the distribution of the eigenvalues of the iteration matrix.
%Moreover,  well-known challenge for the numerical methods of SWEs
%is the simulation of flows involve wet/dry transitions.
%Although the paper focuses on the numerical method for the steady-state SWEs instead of time-dependent SWEs,
%the proposed Newton multigrid method can also be applied to the implicit or semi-implicit
%schemes for time-dependent SWEs to improve the efficiency, see Remark \ref{rem:implicit}.
The  paper is organized as follows.
Section \ref{sec:scheme} presents the Newton multigrid method, including the well-balanced spatial discretization of the steady-state SWEs
  in Subsection \ref{sec:Discre}, the regularization of the resulting nonlinear system
and its Newton iteration linearization in Subsection \ref{sec:newton},
the geometric multigrid method  in Subsection \ref{sec:mg-solver},
and the solution procedure of the Newton multigrid method in Subsection \ref{sec:procedure}.
Section \ref{sec:numerical-results} conducts several numerical experiments to demonstrate the
efficiency and robustness of the proposed Newton multigrid method and presents
a detailed discussion on the relation between the
convergence behavior of the proposed method and the distribution of the eigenvalues of   the block symmetric Gauss-Seidel iteration matrix  for different numerical fluxes.
Section \ref{sec:conclud} concludes the paper with several remarks.

%\section{Governing equations}
%\label{sec:GovernEqns}

%This section introduces

\section{Numerical method}
\label{sec:scheme}

This section is devoted to present the Newton multigrid method for %one- and two-dimensional
steady-state SWEs \eqref{eq:GNeqsSTEEADY}.
%For the sake of simplicity,  the method will be only  described
%on the 2D structured mesh.
% with the spatial step sizes $\Delta x$ and $\Delta y$.
Let $\cal T$ be a partition of the spatial domain $\Omega_p$,
and ${\cal K}_i \in {\cal T}$
be the $i$th cell, whose  centroid is $\vec x_i$.
%, which is an interval (resp.  a quadrangle) in the 1D (resp. 2D) case in our code.
Use $\partial {\cal K}_i$ to denote the edge set of ${\cal K}_i$,
 $e_{ij}$ to be the edge   of   ${\cal K}_i$ sharing with
the neighboring cell ${\cal K}_j$, i.e. $e_{ij} = \partial {\cal K}_i\cap\partial {\cal K}_j$,
$\vec n_{ij}=( n_{ij}^{x}, n_{ij}^{y} )^{\rm T}$ to be the unit normal vector of $e_{ij}$, pointing from ${\cal K}_i$
to ${\cal K}_j$, and $\left| {e_{ij} } \right|$ to denote the length of edge $e_{ij}$.

\subsection{Well-balanced finite volume discretization}
\label{sec:Discre}
This section presents the well-balanced finite volume discretization of the steady-state SWEs \eqref{eq:GNeqsSTEEADY}.
%The cell-centered scheme is adopted here to discretize the steady-state shallow water equations \eqref{eq:GNeqsSTEEADY}.
Integrating   \eqref{eq:GNeqsSTEEADY} over the
cell ${\cal K}_i$ and using the divergence theorem give
\begin{equation}\label{eq:IntEQS}
\sum\limits_{e_{ij}  \in \partial {\cal K}_i } {\int_{e_{ij} }
{{   {{\vec F}_{{\smallvec n}_{ij}}}    (\vec U)}~ {\rm d} S} }  =
 {\iint_{ {\cal K}_i } { \vec S(\vec x,\vec U)~ {\rm d} \vec x} },
\end{equation}
where  ${ {{\vec F}_{{\smallvec n}_{ij}}} (\vec U)}:=  {\vec F (\vec U)} \cdot \vec n_{ij}$.
If let $\vec u=\vec 0$, then \eqref{eq:IntEQS}
reduces to
\begin{equation}\label{eq:IntEQS2}
\sum\limits_{e_{ij}  \in \partial {\cal K}_i }
{\int_{e_{ij} } {\vec S_{ij}(h)}  ~ {\rm d} S}:=\sum\limits_{e_{ij}  \in \partial {\cal K}_i } {\int_{e_{ij} }
\frac12 gh^2  \widetilde{\vec n}_{ij} ~ {\rm d} S}   =
 {\iint_{ {\cal K}_i } { \vec S(\vec x,\vec U)~ {\rm d} \vec x} },
\end{equation}
where $\widetilde{\vec n}_{ij}=(0, n_{ij}^{x}, n_{ij}^{y} )^{\rm T}$.
It is a  starting point of the hydrostatic reconstruction method \cite{Audusse2004} to derive a well-balanced method.

%Combining \eqref{eq:IntEQS} with \eqref{eq:IntEQS2} gives
%\begin{equation}\label{eq:IntEQS3}
%\sum\limits_{e_{ij}  \in \partial {\cal K}_i } {\int_{e_{ij} }
 %{   \left( {\vec F_{\small n_{ij}} (\vec U)}-{\vec S_{ij}(h)}\right) ~ {\rm %d} S} }   =0.
%\end{equation}
Reconstruct a piecewise polynomial $\vec U_h(\vec x)$
to  approximate the  solution $\vec U(\vec x)$,  e.g.
%where
%each component of $\vec U_h$ is {piecewise constant function} (discontinuous at cell interface) with
\begin{equation}\label{eq:APsolu}
\vec U_h (\vec x) = \overline{\vec U}_i ,\quad \forall \vec x \in {\cal K}_i,
\end{equation}
%or
%\begin{equation}\label{eq:APsolu2}
%\vec U_h (\vec x) = \overline{\vec U}_i +   \overline{\nabla\vec U}_i
%(\vec x-\vec x_i),\quad \forall \vec x \in {\cal K}_i,
%\end{equation}
where $\overline{\vec U}_i$ % and  $\overline{\nabla\vec U}_i$ are
is the cell average approximation
of $\vec U(\vec x)$  over   ${\cal K}_i$.
Replacing $\vec U(\vec x)$ in \eqref{eq:IntEQS} with $\vec U_h(\vec x)$,
using the midpoint quadrature to evaluate the integral, and
approximating  the values of  ${\vec F_{   n_{ij}} (\vec U)}$ and
$\vec S_{ij}(h)$ at the midpoint $\vec x_{e_{ij}}$ of the edge $e_{ij}$ by numerical fluxes
give the finite volume discretization of \eqref{eq:GNeqsSTEEADY} as follows
\begin{equation}\label{eq:DSeqs}
\sum\limits_{e_{ij}  \in \partial {\cal K}_i } {\left| {e_{ij} } \right| \widehat{\vec {\cal F}} \left( {\overline{\vec U}_i ,\overline{\vec U}_j ,z_i ,z_j } \right)} = \vec 0,
\end{equation}
with
%$\widehat{\vec {\cal F}}\left( {\overline{\vec U}_i ,\overline{\vec U}_j %,z_i ,z_j } \right)$ is the ``numerical flux'' taking the discretization of the source term
%into account, and has the following form
\begin{equation}\label{eq:totalFlux}
%\widehat{\vec {\cal F}}\left( {\overline{\vec U}_i ,\overline{\vec U}_j ,z_i %,z_j } \right)
% := \widehat {\vec F}_{n_{ij}}\left( \vec U_i^{[j]} ,\vec U_j^{[i]}  \right)   -  {\vec S}_{ij} \left( {h_i^{[j]} } \right) +  {\vec S}_{ij} \left( {\overline {h}_i } \right),
 \widehat{\vec {\cal F}}\left( {\overline{\vec U}_i ,\overline{\vec U}_j ,z_i ,z_j } \right)
  := {\widehat{\vec F}_{{\smallvec n}_{ij}}}   \left( \vec U_{e_{ij},-},\vec U_{e_{ij},+} \right)
  -   {\vec S}_{ij} {  \left( \hat{h}_{e_{ij},-} \right) },
\end{equation}
where $\vec U_{e_{ij},\pm} { =: \left( {h}_{e_{ij},-}, {h}_{e_{ij},-} {\vec u}_{e_{ij},-}\right)^{\rm T}}$ denote the left and right limit values of $\vec U_h(\vec x)$
at the point $\vec x_{e_{ij}}$ in the direction ${\vec n}_{ij}$,
$\widehat {\vec F}_{n_{ij}}\left( \vec U_{e_{ij},-},\vec U_{e_{ij},+}  \right) $
is any given numerical flux
satisfying the consistency
$$ {\widehat{\vec F}_{{\smallvec n}_{ij}}}   \left( \vec U ,\vec U  \right)
={\vec F}_{n_{ij}}\left( \vec U   \right)  ,
$$ %the estimate of outward normal flux across edge $e_{ij}$ of cell ${\cal K}_i$,
 $$ {\vec S}_{ij}
 {   \left( \hat{h}_{e_{ij},-} \right) }
 =\left(0, n^x_{ij} \frac{g}{2} {   {h}_{e_{ij},-}^2 } , n^y_{ij}\frac{g}{2}  {   {h}_{e_{ij},-}^2 } \right)^{\rm T},
 $$
 and
 \begin{align}
 \begin{aligned}
{   {h}_{e_{ij},-} }=&\max\big\{0,\bar h_i+z_i-\max\{z_i,z_{j}\}\big\},\quad
 {  \vec u_{e_{ij},-}=\bar {\vec u}_i, }
 \\
{  {h}_{e_{ij},+} }=&\max\big\{0,\bar h_{j}+z_{j}-\max\{z_i,z_{j}\}\big\},
 \quad
  {
 \vec u_{e_{ij},+}=\bar {\vec u}_{j}. }
 \end{aligned}
 \end{align}
 Such locally reconstructed heights can roughly capture dry regions where
 $h=0$, see \cite{Audusse2004} for a more detailed discussion.
%$\vec U_i^{[j]}$ is defined by
%$$
%\vec U_i^{[j]}  =
%h_i^{[j]} \left( 1, \bar{u}_{i},  \bar{v}_{i} \right)^{\rm T}, \quad
%h_i^{[j]} : = \max \left\{ {0,\bar{h}_i  + z_i  - \max \left\{ {z_i ,z_j } \right\}} \right\}.
%$$

\begin{remark}
For 1D steady-state SWEs,
first-order accurate well-balanced scheme may be described as follows
\begin{align}
%\Delta x\frac{d}{dt} \overline{\vec U}_i+
\hat{\vec F}_{i+\frac12}-\hat{\vec F}_{i-\frac12}=\vec S_i,
\end{align}
where numerical flux $\hat{\vec F}_{i+\frac12}$ and source term $\vec S_i$ are
\begin{align}\hat{\vec F}_{i+\frac12}=\hat{\vec F}(\vec U_{i+\frac12,-},\vec U_{i+\frac12,+}),\quad
\vec S_i=\left(0, \frac{g}{2} {(
{ {h}^2_{i+\frac12,-}-{h}^2_{i-\frac12,+} } )  }
\right)^{\rm T}.
\end{align}
and the left and right limits of the approximate solution at $x_{i+\frac12}$ are
\begin{align}
\begin{aligned}
{{h}_{i+\frac12,-} }=&\max\big\{0,  {  \bar h_i}+z_i-\max\{z_i,z_{i+1}\}\big\},\ \ \qquad u_{i+\frac12,-}= { \bar u_i},\\
{  {h}_{i+\frac12,+} }=&\max\big\{0,  {  \bar h_{i+1} }+z_{i+1}-\max\{z_i,z_{i+1}\}\big\},
\quad u_{i+\frac12,+}={ \bar u_{i+1}}.
\end{aligned}
\end{align}

\end{remark}

%and
%$$
%{\vec S}_{ij}(\vec U) := \left( {0,gh^2 n_{ij}^{x} /2,gh^2 n_{ij}^{y} /2} \right)^{\rm T}.
%$$
%Define the rotation matrix and its inverse
%\[
%\vec T_{\vec n}  =
%\left[ {\begin{array}{*{20}c}
%   1 & 0 & 0  \\
%   0 & {n^{x} } & {n^{y} }  \\
%   0 & { - n^{y} } & {n^{x} }
%\end{array}} \right], \quad
%\vec T_{\vec n}^{-1}  = \left[ {\begin{array}{*{20}c}
%   1 & 0 & 0  \\
%   0 & {n^{x} } & { - n^{y} }  \\
%   0 & {n^{y} } & {n^{x} }
%\end{array}} \right]
%\]
%for a given normal vector $\vec n = \left( n^{x},n^{y}\right)^{\rm T}$.
%Then using the rotational invariance property \cite{LeVeque2002} of the 2D SWEs gives
%\[
%\widehat {\vec F}\left( {\vec U_L ,\vec U_R } \right) \cdot \vec n
%= \vec T_{\vec n}^{ - 1} \widehat{\vec F_1 }\left( {\vec T_{\vec n} \vec U_L ,\vec T_{\vec n} \vec U_R } \right)
%= :\vec T_{\vec n}^{ - 1} \widehat{\vec F_1 }\left( {\widetilde {\vec U}_{\vec n,L} ,\widetilde {\vec U}_{\vec n,R} } \right),
%\]
%where $\widetilde {\vec U}_{\vec n}  = \vec T_{\vec n} \vec  U = \left( {h,hu_n ,hu_\tau  } \right)^{\rm T}$, and
%$u_n  := \vec u \cdot \vec n,~u_\tau   :=  - u n^{y}  + v n^{x}$ respectively denote the
%velocities perpendicular and tangential to the cell face with the outer unit normal vector $\vec n$.
%In this way, the evaluation of the numerical flux $\widehat {\vec F} \cdot \vec n_{ij}$ in \eqref{eq:totalFlux}
%can be given by an (exact or approximate) Riemann solver of
%the corresponding local one-dimensional Riemann problem in the direction normal to cell
%interface $e_{ij}$.

\begin{remark}\label{Rem2.2}
Three numerical fluxes will be considered in this work. The first is
 the HLLC (resp. HLL) flux for 2D (resp. 1D) SWEs, see \cite{Liang2009,Toro1994}.
 The 2D HLL flux is defined by
\begin{align}\nonumber
&\widehat {\vec F}^{\mbox{\tiny HLL}}_{{  \smallvec n}_{ij}}\left( \vec U_{e_{ij},-},\vec U_{e_{ij},+}  \right)
\\[2mm]
&
= \begin{cases}
{ {\vec F}_{{\smallvec n}_{ij}}} \left( \vec U_{e_{ij},-} \right), & \mbox{} s_L \ge 0,\\
\dfr{ s_R  { {\vec F}_{{\smallvec n}_{ij}}}  \left( \vec U_{e_{ij},-} \right) 
- s_L  { {\vec F}_{{\smallvec n}_{ij}}}  \left( \vec U_{e_{ij},+} \right) + s_L s_R
\left(  \vec U_{e_{ij},+} - \vec U_{e_{ij},-} \right) } { s_R - s_L },  & \mbox{} s_L < 0 < s_R,  \\
 { {\vec F}_{{\smallvec n}_{ij}}} \left( \vec U_{e_{ij},+} \right), & \mbox{} s_R \le 0,
\end{cases}
\label{eq:HLLflux}
\end{align}
where
$$
s_L = \begin{cases}
u_{n_{ij},+}  - 2 \sqrt{g h_{e_{ij},+}}, & \mbox{} h_{e_{ij},-} =0,\\
\min\left\{ u_{n_{ij},-}  -  \sqrt{g h_{e_{ij},-}} , u_{n_{ij},*}  -  \sqrt{g h_{e_{ij},*}} \right\}, & \mbox{}{h_{e_{ij},-} \neq 0},
\end{cases}
$$
and $$
s_R = \begin{cases}
u_{n_{ij},-}  + 2 \sqrt{g h_{e_{ij},-}}, & \mbox{} h_{e_{ij},+} =0,\\
\max\left\{ u_{n_{ij},+}  +  \sqrt{g h_{e_{ij},+}} , u_{n_{ij},*}  +  \sqrt{g h_{e_{ij},*}} \right\}, & \mbox{}{h_{e_{ij},+} \neq 0},
\end{cases}
$$
here $u_{n_{ij},\pm} := \vec u_{e_{ij},\pm} \cdot \vec n_{ij}$, and
\begin{align*}
&
u_{n_{ij},*} = \frac{ u_{n_{ij},-} + u_{n_{ij},+} }{2}  + \sqrt{g h_{e_{ij},-} } - \sqrt{g h_{e_{ij},+} } , \\
&
h_{e_{ij},*} = \frac{1}{4g} \left[ \sqrt{g h_{e_{ij},-} } + \sqrt{g h_{e_{ij},+} } + \frac{ u_{n_{ij},-} - u_{n_{ij},+} }{2}   \right]^2.
\end{align*}
Based on the rotational invariance property of 2D SWEs \cite{LeVeque2002}
$$
\vec T_{n_{ij}}
\widehat {\vec F}^{\mbox{\tiny HLL}}_{ { \smallvec n}_{ij}}\left( \vec U_{e_{ij},-},\vec U_{e_{ij},+}  \right)
= \widehat {\vec F}^{\mbox{\tiny HLL}}_{1}\left( \vec T_{  {\smallvec n}_{ij}} \vec U_{e_{ij},-}, \vec T_{ {\smallvec n}_{ij}}  \vec U_{e_{ij},+}  \right)
=:\left( \widehat { f}_{1}, \widehat { f}_{2}, \widehat {f}_{3} \right)^{\rm T},
$$
with
$$
\vec T_{ {\smallvec n}_{ij}}  =\begin{pmatrix}
   1 & 0 & 0  \\
   0 & {n^{x}_{ij} } & {n^{y}_{ij} }  \\
   0 & { - n^{y}_{ij} } & {n^{x}_{ij} }
\end{pmatrix},
$$
 the HLLC flux may be given by
$$
\widehat {\vec F}^{\mbox{\tiny HLLC}}_{  {\smallvec n}_{ij}}\left( \vec U_{e_{ij},-},\vec U_{e_{ij},+}  \right) = \vec T^{-1}_{ {\smallvec n}_{ij}}
\widehat {\vec F}^{\mbox{\tiny HLLC}}_{1}\left(  \vec T_{  {\smallvec n}_{ij}}  \vec U_{e_{ij},-},  \vec T_{ {\smallvec n}_{ij}}  \vec U_{e_{ij},+}  \right),
$$
where
$$
\widehat {\vec F}^{\mbox{\tiny HLLC}}_{1}\left(  \vec T_{ {\smallvec n}_{ij}}  \vec U_{e_{ij},-},  \vec T_{ {\smallvec n}_{ij}}  \vec U_{e_{ij},+}  \right)
= \begin{cases}
\left( \widehat { f}_{1}, \widehat { f}_{2}, \widehat {f}_{3} \right)^{\rm T}, & \mbox{} s_L \ge 0  \mbox{ or } s_R \le 0,\\
\left( \widehat { f}_{1}, \widehat { f}_{2},  u_{\tau_{ij},-}  \widehat {f}_{1} \right)^{\rm T}, & \mbox{} s_L < 0 \le s_M, \\
\left( \widehat { f}_{1}, \widehat { f}_{2},  u_{\tau_{ij},+}   \widehat {f}_{1} \right)^{\rm T}, & \mbox{} s_M < 0 < s_R,
\end{cases}
$$
with
$$
u_{\tau_{ij},\pm} := -  u_{e_{ij},\pm} n^{y}_{ij} +   v_{e_{ij},\pm} n^{x}_{ij} ,
$$
and
$$
s_M : = \dfr{ s_L h_{e_{ij},+} \left( u_{n_{ij},+} - s_R \right) -s_R h_{e_{ij},-} \left( u_{n_{ij},-} - s_L \right) }
{ h_{e_{ij},+} \left( u_{n_{ij},+} - s_R \right) - h_{e_{ij},-} \left( u_{n_{ij},-} - s_L \right) }.
$$

 The second   is the local Lax-Friedrichs (LLF) flux
\begin{equation}\label{eq:LLFflux}
\widehat {\vec F}^{\mbox{\tiny LLF}}_{ {\smallvec n}_{ij}}\left( \vec U_{e_{ij},-},\vec U_{e_{ij},+}  \right)
=
\frac{1}{2}\big(    {{\vec F}_{{\smallvec n}_{ij}}}  ( \vec U_{e_{ij},-}) +
{{\vec F}_{{\smallvec n}_{ij}}} ( \vec U_{e_{ij},+})
 - s_{\rm max} ( {\vec U_{e_{ij},+}  - \vec U_{e_{ij},-} })
  \big),
%
%\widehat{\vec F_1 }^{\rm LLF} \left( {\vec U_L ,\vec U_R } \right) = \frac{1}{2}\big( {\vec F_1 (\vec U_L ) + \vec F_1 (\vec U_R ) - s_{\rm max} \left( {\vec U_R  - \vec U_L } \right)} \big),
\end{equation}
where $s_{\rm max}$ denotes an estimation of the fastest wave speed
in the local 1D Riemann problem solution, and is usually taken as an upper bound for the absolute value
of eigenvalues of the Jacobian $\partial \vec F_{n_{ij}}/\partial \vec U$ as follows
$$
s_{\rm max} = \max \left\{ |u_{n_{ij},-}| + \sqrt{g h_{e_{ij},-}} ,  |u_{n_{ij},+}| + \sqrt{g h_{e_{ij},+}}   \right\}.
$$
If there exist wet/dry transitions, then
because the speed of a wet/dry front
is of form $S_{*+}=u_{n_{ij},+}-2 \sqrt{gh_{e_{ij},+}}$
for a left dry state and $S_{*-}=u_{n_{ij},-} + 2 \sqrt{g {  h_{e_{ij},-} } }$
for a right dry state \cite{Toro2001},
 $s_{\rm max}$ should be appropriately larger in the presence of wet/dry fronts to avoid
numerical instabilities and taken as follows
$$
s_{\rm max} = \max \left\{ |u_{n_{ij},-}| + \sqrt{gh_{e_{ij},-}} ,
|u_{n_{ij},+}| + \sqrt{gh_{e_{ij},+} }   \right\} +
\varepsilon_r \sqrt{g \max \left\{h_{e_{ij},-} , h_{e_{ij},+} \right\}},
$$%{\rm e}^{-\min \{h_L,h_R\}} \sqrt{gh_*},
where $\varepsilon_r = 0.03$ in our computations.
%, and $h_*$ is an estimate for the exact solution for $h$ in the star region.
%We use the two-rarefaction Riemann solver \cite{Toro2001}, which gives
%\[
%h_*  = \frac{1}{g}\left[ {\frac{1}{2}\left( {\sqrt {gh_L }  + \sqrt {gh_R } } \right) + \frac{1}{4}\left( {u_{L}  - u_{R} } \right)} \right]^2.
%\]
% and these speeds can be larger than those associated with the eigenvalues $u \pm \sqrt{gh}$.

The third   is the Roe flux with Harten's entropy fix
\begin{equation}\label{eq:ROEflux}
 \widehat {\vec F}^{\mbox{\tiny Roe}}_{ { \smallvec n}_{ij}}
\left( \vec U_{e_{ij},-},\vec U_{e_{ij},+}  \right)
=
\frac{1}{2}\big(  { {\vec F}_{{\smallvec n}_{ij}}} ( \vec U_{e_{ij},-}) +
{{\vec F}_{{\smallvec n}_{ij}}}( \vec U_{e_{ij},+})
 - \sum\limits_{k=1}^{m}
 Q \left( {\widehat \lambda_k } \right)\widehat\chi _k \widehat {\vec r}_k
  \big),
%{\widehat{\vec F_1 }}^{\rm ROE} \left( {\vec U_L ,\vec U_R } \right) = %\frac{1}{2}\left( {\vec F_1 (\vec U_L ) + \vec F_1 (\vec U_R ) - %\sum\limits_{k = 1}^m {Q^H \left( {\widehat \lambda_k } \right)\widehat\chi %_k \widehat {\vec r}_k } } \right),
\end{equation}
where  $m$ is equal to $2$ (resp. 3) for the 1D (resp. 2D) case,
 $\widehat \lambda_k = \widehat \lambda_k(\vec U_{e_{ij},-}, \vec U_{e_{ij},+})$ and
$\widehat {\vec r}_k = \widehat {\vec r}_k(\vec U_{e_{ij},-}, \vec U_{e_{ij},+})$ are the eigenvalues and  corresponding right eigenvectors
of the Roe matrix, respectively,
$\widehat {\chi}_k$ solves the system %are the components of the jump $\vec U_{e_{ij},+} - \vec U_{e_{ij},-}$ on
%these eigenvectors, namely,
\[
\vec U_{e_{ij},+} - \vec U_{e_{ij},-} = \sum\limits_{k = 1}^m {\widehat\chi _k \widehat {\vec r}_k },
\]
and
\begin{equation*}
Q (x)=
\begin{cases}
\frac{{x^2 }}{{4\varepsilon_f }} + \varepsilon_f, &  |x|<2\varepsilon_f,\\
|x|,&    |x|\geq 2\varepsilon_f,\end{cases}
\end{equation*}
where  $\varepsilon_f$  is a small positive constant, e.g. $0.4$ in later computations.

\end{remark}

%hydro-static reconstruction method
%
%and the boundary and volume integrals in \eqref{eq:IntEQS} are discretized %by midpoint rule
%as follows
%\begin{equation*}
% \int_{e_{ij} } {\vec F\left( {\vec U_h } \right) \cdot \vec n_{ij} {\rm d} s}  \approx \left| {e_{ij} } \right| \vec F\left( {\vec U_h \left(x_{ij}^M \right)} \right) \cdot \vec n_{ij},  \quad
% \int {\int_{   {\cal K}_i } { \vec S\left( {\vec x,\vec U_h } \right) {\rm d} \vec x} }  \approx \left| { {\cal K}_i } \right| \vec S_i ,
%\end{equation*}
%where $x_{ij}^M$ denote the midpoint of the edge $e_{ij}$, and $\vec S_i$ is approximation of $\vec S(\vec x,\vec U)$ at the center of the cell ${\cal K}_i$ which
%should be considered carefully to obtain a well-balanced scheme.
%and $\vec x_i^C$ respectively denote the midpoint of the edge $e_{ij}$ and
%
%the center of the cell ${\cal K}_i$.

\begin{remark} \label{rem:boundary}
%The solution of the steady-state SWEs \eqref{eq:GNeqsSTEEADY} is highly dependent on the
%boundary conditions. Besides, the boundary conditions are also essential to the well-posedness,
%convergence of a scheme for the discrete problem \eqref{eq:DSeqs}.
Let us discuss the boundary conditions for the discrete problem \eqref{eq:DSeqs}.
 If the boundary is open, according
to the local Froude number $F_r:=|\vec u|/\sqrt{gh}$,
the numerical boundary
conditions are specified  as follows:
\begin{itemize}
\item Subcritical inflow boundary $\{F_r<1,u_n <0\}$:
${\cal R}_{\mbox{\rm g}}^- = {\cal R}_I^-$, and $h_{\mbox{\rm g}} u_{n,{\mbox{\rm g}}}$ and $u_{\tau,{\mbox{\rm g}}}$ are prescribed.
\item Subcritical outflow boundary $\{F_r<1,u_n >0\}$:
${\cal R}_{\mbox{\rm g}}^- = {\cal R}_I^-$,  $u_{\tau,{\mbox{\rm g}}} = u_{\tau,I}$, and $h_{\mbox{\rm g}}$ is prescribed.
\item Supercritical inflow boundary $\{F_r>1,u_n <0\}$: $h_{\mbox{\rm g}}$, $u_{n,{\mbox{\rm g}}}$, and $u_{\tau,{\mbox{\rm g}}}$ are  prescribed.
\item Supercritical outflow boundary $\{F_r>1,u_n >0\}$: $h_{\mbox{\rm g}}=h_I$, $u_{n,{\mbox{\rm g}}} = u_{n,I}$, and $u_{\tau,{\mbox{\rm g}}} = u_{\tau,I}$.
\end{itemize}
Here ${\cal R}^\mp := u_n \pm 2 \sqrt{gh}$ denotes  1D Riemann invariant associated with the eigenvalues $u_n \mp \sqrt{gh}$, % along the  normal vector of the boundary.
 $u_n$ and $u_{\tau}$ are the normal and tangential velocity components  to the cell interface located on the domain boundary, respectively.
The quantities with subscripts ${\mbox{\rm g}}$ and $I$
denote the values from the ghost and interior cells adjacent to the domain boundary, respectively.

The slip boundary conditions
$\{h_{\mbox{\rm g}}=h_I$,$u_{n,{\mbox{\rm g}}}=0$, and $u_{\tau,{\mbox{\rm g}}}=u_{\tau,I}\}$
or the reflective boundary conditions
$\{h_{\mbox{\rm g}}=h_I$, $u_{n,{\mbox{\rm g}}}= - u_{n,I}$, and $u_{\tau,{\mbox{\rm g}}}=u_{\tau,I}\}$
 may be specified on the wall.

\end{remark}

\subsection{Newton's iterative method}
\label{sec:newton}

%Assume an initial solution over the cell ${\cal K}_i$ is approximated by
%\[
%\overline {\vec U}_i^{(0)}  = \frac{1}{{\left| {\cal K}_i  \right|}}\int {\int_{ {\cal K}_i } {\vec U^{(0)} (\vec x) {\rm d} \vec x} }
%\]

As soon as the boundary conditions are specified,
the approximate solutions of the SWEs \eqref{eq:GNeqsSTEEADY} may be obtained by
iteratively solving the nonlinear algebraic system \eqref{eq:DSeqs} with respect to
the unknown variables $\overline{\vec U}_i$.
Here Newton's iteration method is employed to solve \eqref{eq:DSeqs}
with the formula
%Newton's iteration formula for  the  system \eqref{eq:DSeqs} becomes
%, with which \eqref{eq:DSeqs} can be linearized as
\begin{equation}\label{eq:Newton}
 \sum\limits_{e_{ij}  \in \partial {\cal K}_i } {\left| {e_{ij} } \right|
	\left( {\frac{{\partial          \widehat {\vec {\cal F}}          }}{{\partial \overline{\vec U}_i }}} \right)^{(n)}_{ij} } \delta \vec U_i^{(n)}
+ \sum\limits_{e_{ij}  \in \partial {\cal K}_i } {\left| {e_{ij} } \right|
	\left( {\frac{{\partial \widehat {\vec {\cal F}} }}{{\partial\overline{\vec U}_j }}} \right)^{(n)}_{ij} } \delta \vec U_j^{(n)}  =-\vec R_i^{(n)} , \ \ n=0,1,\cdots,
\end{equation}
where the unknown $\delta \vec U_j^{(n)}$ is the  increment,
$\overline{\vec U}^{(0)}_i$ is the given initial guess,
$\vec R_i^{(n)}$ is the local residual at the $n$th Newtonian iterative step defined by
\begin{equation}\label{eq:Res}
\vec R_i^{(n)} := \sum\limits_{e_{ij}  \in \partial {\cal K}_i } {\left| {e_{ij} } \right| \widehat {\vec {\cal F}}\left( {\overline{\vec U}_i^{(n)} ,\overline{\vec U}_j^{(n)} ,z_i ,z_j } \right)},
\end{equation}
%= \overline{\vec U}_j^{(n+1)}
%-\overline{\vec U}_j^{(n)}$,
and  the partial derivatives
\begin{align}\label{eq:jacobi1}
&
\left( {\frac{{\partial \widehat {\vec {\cal F}} }}{{\partial \overline{\vec U}_i }}} \right)^{(n)}_{ij}  := \frac{{\partial \widehat {\vec {\cal F}}}}{{\partial \overline{\vec U}_i }}
\left( { \overline{\vec U}_i^{(n)} ,\overline{\vec U}_j^{(n)} ,z_i ,z_j } \right),\\ \label{eq:jacobi2}
&
\left( {\frac{{\partial \widehat {\vec {\cal F}} }}{{\partial\overline{\vec U}_j }}} \right)^{(n)}_{ij}  := \frac{{\partial\widehat { \vec {\cal F}}}}{{\partial \overline{\vec U}_j }}
\left( {\overline{\vec U}_i^{(n)} ,\overline{\vec U}_j^{(n)} ,z_i ,z_j } \right),
\end{align}
  contributing to  the Jacobian matrix in the Newtonian method
are approximately
calculated by using the numerical differentiation  as follows
\begin{equation}\label{eq:NumerDiff}
\left( \left( {\frac{{\partial \widehat {\vec {\cal F}} }}{{\partial \overline{\vec U}_i }}} \right)^{(n)}_{ij} \right)_{s,k}
\approx \frac{  \widehat {\cal F}_s\left( {\overline{\vec U}_i^{(n)} + \epsilon \vec e_k,\overline{\vec U}_j^{(n)} ,z_i ,z_j } \right)    -
	\widehat {\cal F}_s \left( {\overline{\vec U}_i^{(n)} ,\overline{\vec U}_j^{(n)} ,z_i ,z_j } \right) }
{ \epsilon },
\end{equation}
which denotes the derivative of the $s$th component of vector $\widehat {\vec {\cal F}}$,
denoted by  $\widehat {\cal F}_s$,  with respect to the $k$th
component of vector $\vec U_i$,
where
 $\vec e_k$ is the $k$th column vector of the identity matrix of the same size as
$\left( {{{\partial \widehat {\vec {\cal F}} }} / {{\partial\overline{\vec U}_i }}} \right)^{(n)}_{ij}$.
%and $\epsilon$ \eqref{eq:NumerDiff} is chosen as $10^{-8}$ in our implementation.
	If an edge of the interior cell ${\cal K}_i$ is located on the boundary of $\Omega_p$,
%	denoted by  $e_{ib}$,
then \eqref{eq:NumerDiff} should be replaced with
	\begin{equation*}
	\left( \left( {\frac{{\partial\widehat { \vec {\cal F}} }}{{\partial \overline{\vec U}_i }}} \right)^{(n)}_{i{\mbox{\rm g}}} \right)_{s,k}
	\approx \frac{ \widehat  {\cal F}_s\left( {\overline{\vec U}_i^{(n)} + \epsilon \vec e_k,\overline{\vec U}_{\mbox{\rm g}} (\overline{\vec U}_i^{(n)} + \epsilon \vec e_k) ,z_i ,z_j } \right)    -
		\widehat {\cal F}_s \left( {\overline{\vec U}_i^{(n)} ,\overline{\vec U}_{\mbox{\rm g}}(\vec U_i^{(n)}) ,z_i ,z_j } \right) }
	{ \epsilon },
	\end{equation*}
	where   $\overline{\vec U}_{\mbox{\rm g}} $ denotes the approximate  cell-average value of $\vec U$ in the ghost cell
	and is considered as a  function of $\overline{\vec U}_i^{(n)}$ according to the boundary
	conditions given in Remark \ref{rem:boundary}.
Although the calculation of ${ {{\partial \widehat {\vec {\cal F}} }}/{{\partial \overline{\vec U}_i }}}$
can also be obtained by using the
chain rule, but it may be seriously tedious. Moreover,
their analytical expressions are difficultly derived near the domain boundary.

The linear system \eqref{eq:Newton} is generally singular %so it should be regularized.
and  should be regularized. One way is to add an artificial time derivative term into \eqref{eq:Newton}.
%such as \begin{equation}\label{eq:EM-DSeqs}
%\frac{{\overline{\vec U}^{(n + 1)}_i  - \overline{\vec U}^{(n)}_i }}{{\Delta t_n }}
%+
%\sum\limits_{e_{ij}  \in \partial {\cal K}_i } {\left| {e_{ij} } \right|  \widehat{\vec {\cal F}}\left( {\overline{\vec U}_i^{(n)} ,\overline{\vec U}_j^{(n)} ,z_i ,z_j } \right)} = \vec 0,
%\end{equation}
%the iteration actually not able to carried out directly.
An alternative approach is to use the $l^1$--norm of the local residual %$ \left\| \vec R_i^{(n)} \right\|_{\ell^1 }$
to regularize \eqref{eq:Newton} as   follows
\begin{equation}\label{eq:NewtonREG}
\alpha  \left\| \vec R_i^{(n)}   \right\|_{\ell^1 }   \delta \vec U_i^{(n)}
+ \sum\limits_{e_{ij}  \in \partial {\cal K}_i } {\left| {e_{ij} } \right|
	\left( {\frac{{\partial \widehat {\vec {\cal F}} }}{{\partial \vec U_i }}} \right)^{(n)}_{ij} } \delta \vec U_i^{(n)}
+ \sum\limits_{e_{ij}  \in \partial {\cal K}_i } {\left| {e_{ij} } \right|
	\left( {\frac{{\partial \widehat {\vec {\cal F}} }}{{\partial \vec U_j }}} \right)^{(n)}_{ij} } \delta \vec U_j^{(n)}  = -\vec R_i^{(n)},
\end{equation}
where  $\alpha$ is the positive regularization parameter.
Such regularization technique has been used in solving steady-state Euler equations in \cite{LiWang2008}.
%instead of the regularization of the artificial
%time derivative term,
%
Solving the linear system \eqref{eq:NewtonREG} for the unknown  $\delta \vec U_j^{(n)}$ by the MG method
will be discussed  in Section \ref{sec:mg-solver}.
If the solution of  the linear system \eqref{eq:NewtonREG} is gotten,
then the approximate solution of  \eqref{eq:DSeqs}  can be updated by
\begin{equation}\label{eq:NewtonREG-12}
\overline{\vec U}_i^{(n + 1)}  = \overline{\vec U}_i^{(n)}  + \tau _i \delta \vec U_i^{(n)} ,
%n=0,1,\cdots,
\end{equation}
where $\tau_i$ is a relaxation parameter on cell ${\cal K}_i$. % and always taken as $1$ in this paper.

\subsection{The geometric multigrid solver}
\label{sec:mg-solver}
This section  extends the geometric multigrid method  \cite{LiWang2008}
to the linear system \eqref{eq:NewtonREG}.
%There are two main
%components of the multigrid method: projection operator and smoother.
%
The geometric multigrid methods described so far need
a hierarchy of geometric grids or meshes $\{{\cal T}_l, l=0,1,\cdots, N_L\}$ for the spatial domain $\Omega_p$,
from the coarsest one ($l = N_L$) to the finest one  ($l = 0$).
On all levels but the coarsest one, the smoother will be applied and on the coarsest
level, the system is usually solved exactly. Assume that the ratio of grid points on 		``neighboring''
grids is  constant throughout the grid hierarchy,  and
each coarser cell ${\cal K}_{i,l+1} \in {\cal T}_{l+1}$ is a union of several neighboring finer cells
(two cells for 1D case and four cells in 2D case) in mesh ${\cal T}_l$, i.e.,
 \[
 {\cal K}_{i,l + 1}  = \cup_{j \in I_{i,l+1}} {{\cal K}_{j,l} },
 \]
where $I_{i,l+1}$ is corresponding index set of  those finer cells ${\cal K}_{j,l}$.
 %
%Efficient solver for the linear system \eqref{eq:NewtonREG} arising from
%Newton linearization is critical to the methods.
%In the fairly used Newton-Krylov method \cite{Brown1990,Chisholm2009,Geuzaine2001}, the Krylov subspace iteration methods, such as the well-known GMRES method,
%are usually adopted to iteratively solve this linear system.

%The original partition $\cal T$ on the spatial domain $\Omega_p$ is denoted by ${\cal T}_0$, which is the finest mesh in the multigrid.
%A sequence of coarser meshes used in the multigrid method are denoted as ${\cal T}_l,~l=1,2,\cdots,N_L-1$.
Under the above assumptions,
the steady-state SWEs \eqref{eq:GNeqsSTEEADY} are discretized and solved by the Newton
iteration on the finest mesh ${\cal T}_0$, see Sections \ref{sec:Discre} and \ref{sec:newton}.
% For the sake of convenience,
 %Before introducing the geometric multigrid method for
The linear system \eqref{eq:NewtonREG} on ${\cal T}_0$ is reformulated in the following matrix-vector form
\begin{equation}\label{eq:MatrixForm}
%{\cal L}_l (\delta\vec U_l):=
\sum\limits_j {\vec A_{ij,l} \delta \vec U_{j,l}  }  =  - \vec R_{i,l} ,  \ l=0,
\end{equation}
where the subscript $l$ marks the mesh level, the superscript $(n)$ in \eqref{eq:NewtonREG} has been omitted for convenience,
and
\begin{align*}
&
\vec A_{ii,l}  = \alpha  \left\| \vec R_{i,l}    \right\|_{\ell^1 }
 + \sum\limits_{e_{ij,l}  \in \partial {\cal K}_{i,l} } {\left| {e_{ij,0} } \right|
\left( {\frac{{\partial \widehat{\vec {\cal F}} }}{{\partial \vec U_i }}} \right) _{ij,l} } , \\
&
\vec A_{ij,l} =
\begin{cases}
\left| {e_{ij,l} } \right|
\left( {\frac{{\partial\widehat{ \vec {\cal F}} }}{{\partial \vec U_j }}} \right)_{ij,l} , &  e_{ij}  \in \partial {\cal K}_{i,l}  ,\\
\vec 0,&    \rm{otherwise},
\end{cases}\qquad  i\neq j.
\end{align*}
On the coarser mesh,
the coarse mesh matrices  $\vec A_{ij,l+1}$
are defined by using the Galerkin projection,
and  the source term $ \vec R_{i,l+1}$ is
derived by using the restriction operator $\vec I^{l+1}_l$ that restricts
 those on the fine mesh ${\cal T}_l$ to
the coarse mesh  ${\cal T}_{l+1}$, $l\geq 0$.
In this paper, specifically, they are
%Then the projected linear system on coarser mesh level $l+1$ from fine mesh level $l$ is given as
\begin{equation}\label{eq:ProjectAR}
\vec A_{ij,l+1} = \sum\limits_{{\imath}  \in I_{i,l + 1} } {\sum\limits_{{\jmath}  \in I_{j,l + 1} } {\vec A_{{\imath} {\jmath} ,l} } } ,
\quad
\vec R_{i,l + 1} = \sum\limits_{j \in I_{i,l + 1} } {\left( {\vec R_{j,l}  + \sum\limits_{{\imath}}  {\vec A_{j{\imath} ,l} \delta \vec U_{{\imath} ,l} } } \right)},
\end{equation}
where $\delta \vec U_{j,l}$ is the (approximate) solution of \eqref{eq:MatrixForm} if $l=0$, otherwise solves the following  linear system
\begin{equation}\label{eq:Project}
\sum\limits_j {\vec A_{ij,l} \delta \vec U_{j,l}  }  =  - \vec R_{i,l},\quad l\geq 1.
\end{equation}

As soon as  the correction $\delta \vec U_{j,l+1}$
on the coarse mesh ${\cal T}_{l+1}$ is obtained,
the  correction  $\delta \vec U_{j,l}$
on the fine mesh  ${\cal T}_l$ will be improved as follows
%$\delta \vec U_{i,l+1}$ is added
%back to the solution on the mesh level $l$ as
\begin{equation}\label{eq:correction}
\delta \vec U_{j,l}\longleftarrow
 \delta \vec U_{j,l}  + \vec I_{l+1}^l \delta \vec U_{i,l + 1} ,\quad \forall j \in I_{i,l + 1},
\end{equation}
where $\vec I_{l+1}^l$ denotes the prolongation or coarse-to-fine operator that prolongates or interpolates the correction to the fine mesh from the coarse mesh.

 A multigrid cycle can be defined as a recursive procedure  that is applied at each mesh level as it moves through the grid hierarchy.
 For example, multigrid methods with $\gamma$-cycle has the
 following compact recursive definition.
% \begin{table}[h]

 \noindent
{\tt Algorithm 0}: $\delta \vec U_{j,l}\longleftarrow  \mbox{MG}^{\gamma}_{l}(\delta \vec U_{j,l}, \vec R_{j,l})$\\
 \begin{tabular}{rl}
 (1) & \parbox{13cm}{{\rm Pre smoothing}: Apply the smoother $\nu_1$ times to Eq. \eqref{eq:MatrixForm} or \eqref{eq:Project} with the initial guess $\delta \vec U_{j,l}$.}
 \\
 (2) & If ${\cal T}_l$ is the coarsest grid, i.e. $l=N_L$.
 \\
     & -- solve the problem \eqref{eq:Project} with $l=N_L$.\\
     & else\\
    & -- Restrict  to the next coarser grid ${\cal T}_{l+1}$ by \eqref{eq:ProjectAR}.\\
    &-- Set initial increment on the next coarser grid:
    $\delta \vec U_{j,l+1}=0$.\\
    &--   If ${\cal T}_l$ is the finest grid, set
  $\gamma= 1$.\\
    &--   Call the $\gamma$-cycle scheme
  $\gamma$ times for the next coarser grid ${\cal T}_{l+1}$:\\
 &\qquad  $\delta \vec U_{j,l+1}\longleftarrow  \mbox{MG}^{\gamma}_{l+1}(\delta \vec U_{j,l+1}, \vec R_{j,l+1}).$
 \\
 (3) & Correct with the prolongated update \eqref{eq:correction}.
 \\
 (4) & \parbox{13cm}{{\rm Post smoothing}: Apply the smoother $\nu_2$ times {to Eq.} \eqref{eq:MatrixForm} or \eqref{eq:Project} with the initial guess $\delta \vec U_{j,l}$.}
 \end{tabular}
% \end{table}

\begin{figure}[htbp]
  \centering
    \includegraphics[width=0.48\textwidth]{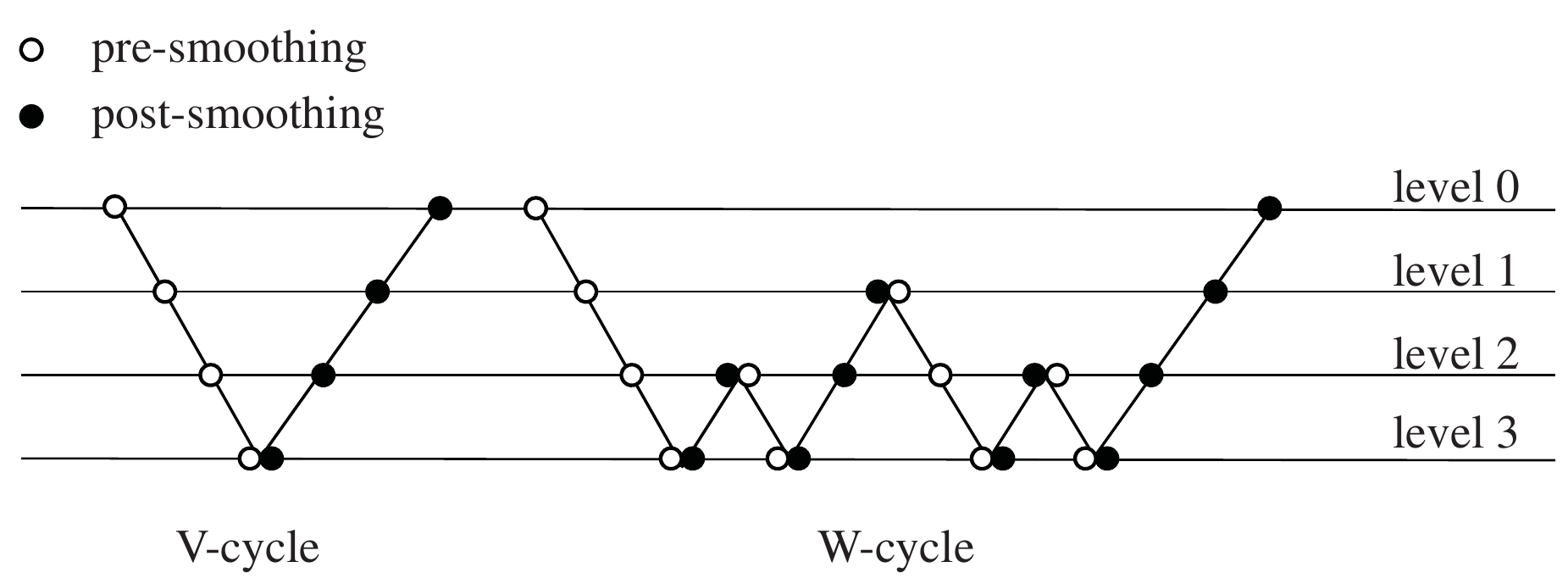}
  \caption{\small  Schematic description of a   V (left) and W cycle (right) Multigrid.
}\label{fig:mgcycle}
\end{figure}

\noindent   This paper only focus on two types of multigrid cycles, the V cycle ($\gamma = 1$) and W cycle ($\gamma = 2$).
%    In practice, only
%    (V-cycle) and
%    $\gamma  = 2$ (W-cycle) are used.
 Fig. \ref{fig:mgcycle} shows their schematic description with $N_L=3$, where symbols ``$\circ$'' and ``$\bullet$'' denote the pre- and post-smoothing, respectively, while the oblique lines
between two symbols ``$\circ$'' (resp. ``$\bullet$'')
correspond to the %fine-to-coarse or
restriction $\vec I^{l+1}_l$ (resp. prolongation   $\vec I_{l+1}^l$) steps.
{The smoother is taken as
the block symmetric Gauss-Seidel (SGS) iteration.}

\subsection{Solution procedure}
\label{sec:procedure}

This section summarizes the solution procedure of our Newton multigrid method for the steady-state SWEs \eqref{eq:GNeqsSTEEADY}.
%
%The Newton multigrid scheme iteratively solves the nonlinear system \eqref{eq:DSeqs} resulting from the well-balanced
%spatial discretization of the steady-state SWEs by using Newton's iterative method presented in subsection \ref{sec:newton} as the outer iteration
% and the geometric multigrid method given in \ref{sec:mg-solver} as the inner iteration.
%The block SGS iteration is adopted as the smoother of the multigrid method, and
%the nonlinear block SGS iteration referred as the BLU-SGS method in \cite{Chen2000} is also performed over the finest mesh as the smoother of the Newton multigrid method.
It is well known that the convergence of Newton¡¯s iteration is seriously dependent on the choice of the initial guess. To overcome this difficulty, the initial guess is obtained by using
{the improved block lower-upper SGS (BLU-SGS) method  \cite{Chen2000}
  to solve
 \eqref{eq:Project}  from the coarsest mesh ${\cal T}_{N_L}$ to the  mesh ${\cal T}_1$ successively.
The idea of BLU-SGS method is to retain the block diagonal matrices but employ LU-SGS-like
 backward and forward Gauss-Seidel iterations to include the implicit contributions from the off diagonal
blocks. This method may be regarded as a ``nonlinear extension'' of the block SGS method for solving the nonlinear algebraic system \eqref{eq:DSeqs}.}

The detailed solution procedure is illustrated by the following flowchart.
%the initial guess

\noindent{\tt Algorithm 1}: Newton multigrid method  (abbr. {\tt NMGM})
\begin{description}
  \item[Step 1: Initialization]
Give the ``initial data'' $\vec U\left( {\vec x } \right)$ and  successively refined partitions $\{{\cal T}_l: l = N_L,\cdots,1,0\}$ of the spatial domain $\Omega_p$.
    %, ${\cal T}_l\subset{\cal T}_{l-1}$. %
  \begin{enumerate}

    \item
    Compute $\overline{\vec U}_{i,N_L}$, the cell average value of $\vec U\left( {\vec x } \right)$ over the coarsest cell ${\cal K}_{i,N_L}$.
    \item For $l = N_L,\cdots, 1$, do the following: %with Block LU-SGS iteration
        \begin{itemize}
          \item Perform the BLU-SGS iteration on the mesh ${\cal T}_l$ %for ${\cal K}_{i,l} \in {\cal T}_l$, symmetrically loop for $i$ increasingly and then increasingly,
         by
            \begin{align*}
            \overline{\vec U}_{i,l} \longleftarrow \overline{\vec U}_{i,l} -
            \left( \alpha \left\| {\vec R_{i,l} } \right\|_{\ell^1 } \vec I  + \frac{{\partial \vec R_{i,l} }}{{\partial \overline{\vec U}_{i,l} }}   \right)^{ - 1} \vec R_{i,l},
            \end{align*}
             until $\sum\limits_i {\left\| {\vec R_{i,l} } \right\|_{\ell^1 } }  < \varepsilon _p 2^{ - l}$, where $\vec I$ denotes the identity matrix.
             % of the same size as $ {{\partial \vec R_{i,l} }}/{{\partial \vec U_{i,l} }} $.
          \item Prolongate the solution to the fine mesh from the coarse mesh by $\overline{\vec U}_{j,l-1}= \overline{\vec U}_{i,l}$ for all $j \in I_{i,l}$.
        \end{itemize}
    \item Set $n=0$ and the initial guess for Newton's iteration by $ {\vec U}_{j,0}^{(0)}= \overline{\vec U}_{i,0}$, for all $j \in I_{i,1}$.
  \end{enumerate}
  \item[Step 2: Newton multigrid iteration]
  For $n=1,2,\cdots, N_{\rm step}$, do the followings.
  \begin{enumerate}
    \item Pre smoothing: perform the BLU-SGS iteration on the finest mesh ${\cal T}_0$ by %: for ${\cal K}_{i,0} \in {\cal T}_0$, symmetrically loop for $i$ increasingly and then increasingly, do
            \begin{align*}
            \vec U_{i,0}^{(n)} \longleftarrow \vec U_{i,0}^{(n)} -
            \left( \alpha \left\| {\vec R_{i,0}^{(n)} } \right\|_{\ell^1 } \vec I  + \frac{{\partial \vec R_{i,0}^{(n)} }}{{\partial \overline{\vec U}_{i,0}^{(n)} }}   \right)^{ - 1}
            \vec R_{i,0}^{(n)}.
            \end{align*}
    \item  Solve \eqref{eq:MatrixForm} %or or \eqref{eq:Project}
     by calling {\tt Algorithm 0}  $N_{mg}$ times.
     %Taking the multigrid strategy of V-cycle for
   % example, it is performed as follows for $N_{mg}$ times:
    %recursively from the finest mesh ${\cal T}_0$ to the coa
 %   \\
 %   for $l = 0, 1, \cdots, N_L -1 $, do the following:
 %   \begin{itemize}
  %     \item Pre-smoothing: perform a few iterations of the block SGS method on the mesh ${\cal T}_l$, leading to solution $\delta \vec U_{i,l}^{(n)}$.
   %    \item Restriction:  Downsampling the residual error to a coarser mesh ${\cal T}_{l+1}$ according to \eqref{eq:Project}--\eqref{eq:ProjectAR} if $l<N_L-1$.
   % \end{itemize}
   % for $l = N_L-1, \cdots, 1,0  $, do the following:
   % \begin{itemize}
    %   \item Prolongation: interpolating a correction computed on a coarser grid into a finer mesh.
     %  $\delta \vec U_{i,l}^{(n)}$ is added back to the solution on the mesh level $l-1$ according to \eqref{eq:correction} if $l>0$.
      % \item Post-smoothing: perform additional block SGS iterations on level $l$ to smooth the corrected solution.
   % \end{itemize}
    \item Update the solution  $\vec U_{i,0}^{(n + 1)}  = \vec U_{i,0}^{(n)}  + \tau _i \delta \vec U_{i,0}^{(n)}$.
  \end{enumerate}
   \item[{Step 3}: ] Check  $\sum\limits_i {\left\| {\vec R_{i,0}^{(n)} } \right\|_{l^1 } }  < \varepsilon $. If yes, output the results and stop; otherwise
 set    $n \leftarrow n+1$ and go to \textbf{Step 2}.
\end{description}

Before ending this section, several remarks are given below.

\begin{remark}
The parameter $\epsilon_p$ in {\bf Step 1} of {\tt Algorithm 1}
 should be chosen appropriately. If  $\epsilon_p$ is very small, the
cost of {\bf Step 1} becomes huge so that  the steady-state {SWEs} solver is inefficient. According to the numerical experiences, {\tt Algorithm 1} works satisfactory if
$\epsilon_p$ is chosen about one percent of the residual given by the initial data.
\end{remark}

\begin{remark}
If the approximate solution $\vec U_{i,0}^{(n )}$ given in the Newton multigrid iteration satisfies $ h_{i,0}^{(n )} < h_{\epsilon}$, the cell ${\cal K}_{i,0}$ is
temporarily regarded as a dry cell, and  the value of $\vec U_{i,0}^{(n )}$ is reset as zero. In our computations, $h_{\epsilon} = 10^{-6}$.
\end{remark}

\begin{remark}
\label{rem:implicit}
The Newton multigrid scheme in {\tt Algorithm 1}
%proposed in this paper is designed for steady-state SWEs \eqref{eq:GNeqsSTEEADY},
 can also be extended to solving  the nonlinear system
 arising from an implicit or semi-implicit scheme for the time-dependent SWEs \eqref{eq:GNeqs}, e.g.
%Specifically, we can use the Newton multigrid method to solve the nonlinear system (about the unknowns $\vec U^{p+1}_i$ if $\vec U^{p}_i$ are known)
\begin{equation*}
\frac{{  \overline{ \vec U}^{n + 1}_i
- \overline{\vec U}^{n}_i }}{{\Delta t_n }}
+
\sum\limits_{e_{ij}  \in \partial {\cal K}_i } {
\left| {e_{ij} } \right|
\left[
\beta \widehat{\vec {\cal F}}\left( {   \overline{\vec U}_i^{n} , \overline{\vec U}_j^{n} ,z_i ,z_j } \right)
+ (1- \beta) \widehat{\vec {\cal F}}\left( { \overline{\vec U}_i^{n+1} ,  \overline{\vec U_j}^{n+1} ,z_i ,z_j } \right)
\right]
} = 0,
\end{equation*}
%resulting from a time-implicit finite volume discretization of \eqref{eq:GNeqs}, and thus the evolution of the numerical solution
%from $t=t_p$ to $t=t_{p+1}$ becomes much efficient in the implicit or semi-implicit scheme.
where $\Delta t_n$ denotes the time step size, and the weight
$\beta \in [0,1)$.
\end{remark}

%%%%%%%%%%%%%%%%%%%%%Numerical Experiments

%%%%%%%%%%%%%%%%%%%%%%%%%%%%%%%%%%%%%%%%%%%%%%%%%%%%%%%%%%%%%%%%%%%%%%%%%%%%%%%%%%%
\section{Numerical Experiments}\label{sec:numerical-results}

The section presents several numerical examples
to demonstrate the robustness
and efficiency of  {\tt NMGM}
%(abbreviated by ``{\tt NMGM}'' for convenience)
for 1D and 2D steady-state
SWEs \eqref{eq:GNeqsSTEEADY}, and investigates
the relation between the convergence behavior of  {\tt NMGM}  and the distribution
of the eigenvalues of the iteration matrix detailedly.
Unless specifically stated,  the parameter $\epsilon_p$  in the initialization step of {\tt Algorithm 1}
is taken as $0.2$,
%the coefficient $\alpha$ used in the regularization is taken as $3$,
and the  multigrid iteration number  $N_{mg}$    is set to be $2$ (resp. $3$) for
1D (resp. 2D) problems. Moreover, the parameters $\epsilon$ in \eqref{eq:NumerDiff},
$\alpha$ in \eqref{eq:NewtonREG}, and $\tau_i$ in \eqref{eq:NewtonREG-12}
are   always taken as $10^{-8}$, $3$,  and 1, respectively.
All computations are carried out on the Linux environment of a personal
computer of Lenovo (Intel(R) Core(TM) i5 CPU 3.2GHZ 4GB RAM).

%For the sake of convenience, the partition $\cal T$ of the physical spatial domain $\Omega_p$ are chosen as uniform mesh for
%the 1D case and structured quadrilateral mesh for the 2D case, which are refined for several times from a
%given coarse mesh.
%

\subsection{1D case}

%%%%%%%%%%%%%%%%%%%%%%%%%%%%%%%%% Example 1D 02 %%%%%%%%%%%%%%%%%%%%%%%%%%%%%%%%%%%%%%%%

\begin{example}[Smooth subcritical flow] \label{example1Dsmooth}\rm
This problem has been studied in \cite{Xu2002} to check the dissipative and dispersive errors in the
kinetic schemes. The bottom shape of the river is
$$
z(x) = 0.2 {\rm{e}}^{ - \frac{(x+1)^2}{2} } + 0.3 {\rm {e}}^{- (x-1.5)^2 }, \quad x\in [-10, 10],
$$
and the boundary conditions at $x=\pm 10$ are specified as $h = 1$  and $hu= 1$.
%The computational domain is $\left[ -10, 10\right]$.
%The discharge $hu = 1$
%and the water height $h = 1$ are imposed
%at $x=-10$ and $10$ respectively.
 Fig. \ref{fig:1d_02_solu} shows the numerical
steady-state solutions obtained by {\tt NMGM} on the mesh of 512 uniform cells, in comparison
with the exact solutions  obtained   by solving the
algebraic system
\begin{align*}
%\begin{aligned}
& u^3 + ( 2gz - 2g - 1 )u + 2g=0, \ \
  hu = 1.
%\end{aligned}
\end{align*}
Fig. \ref{fig:1d_02_Re} displays the numerical residuals obtained by {\tt NMGM}
versus the {\tt NMGM} iteration number  $N_{\rm step}$
 and CPU time for three meshes of 512, 1024, and 2048 uniform cells respectively.
The results show that {\tt NMGM} is very efficient and  fast
to get the correct steady-state solutions.
Moreover, the convergence behaviors are similar
on those different meshes, and the {\tt NMGM} iteration number
 does not increase with refining the mesh.
In those computations, the HLL flux is used, and the grid level number in
the V-cycle multigrid is set to be $4$, i.e. $N_L=3$.

\begin{figure}[htbp]
	\centering
	\includegraphics[width=0.55\textwidth]{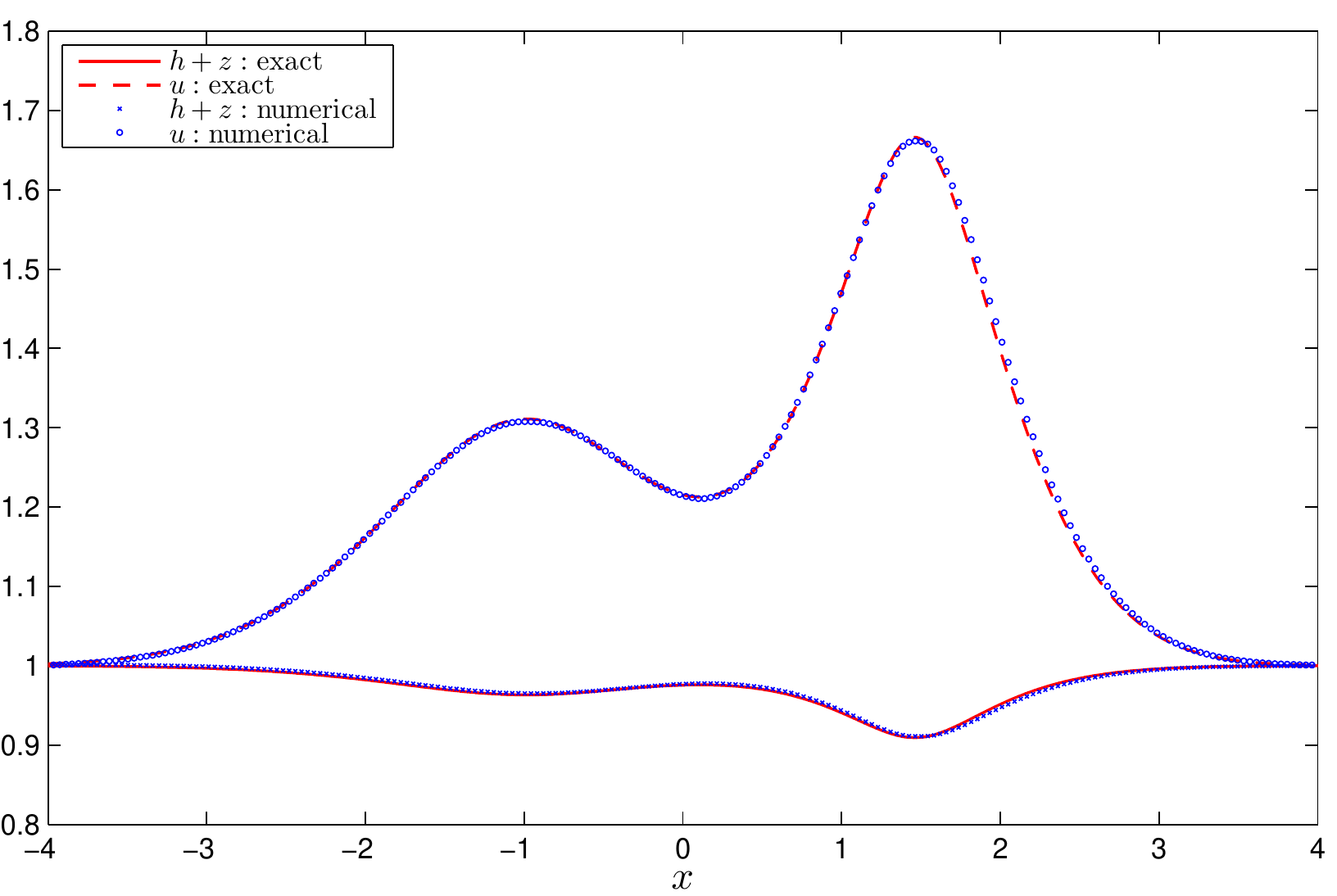}
	\caption{\small Example \ref{example1Dsmooth}: Close-up of the
		steady-state solutions $h+z$ and $u$ obtained by {\tt NMGM}   on the mesh of 512 uniform cells.
	}\label{fig:1d_02_solu}
\end{figure}
\begin{figure}[htbp]
	\centering
	\includegraphics[width=0.46\textwidth,height=5cm]{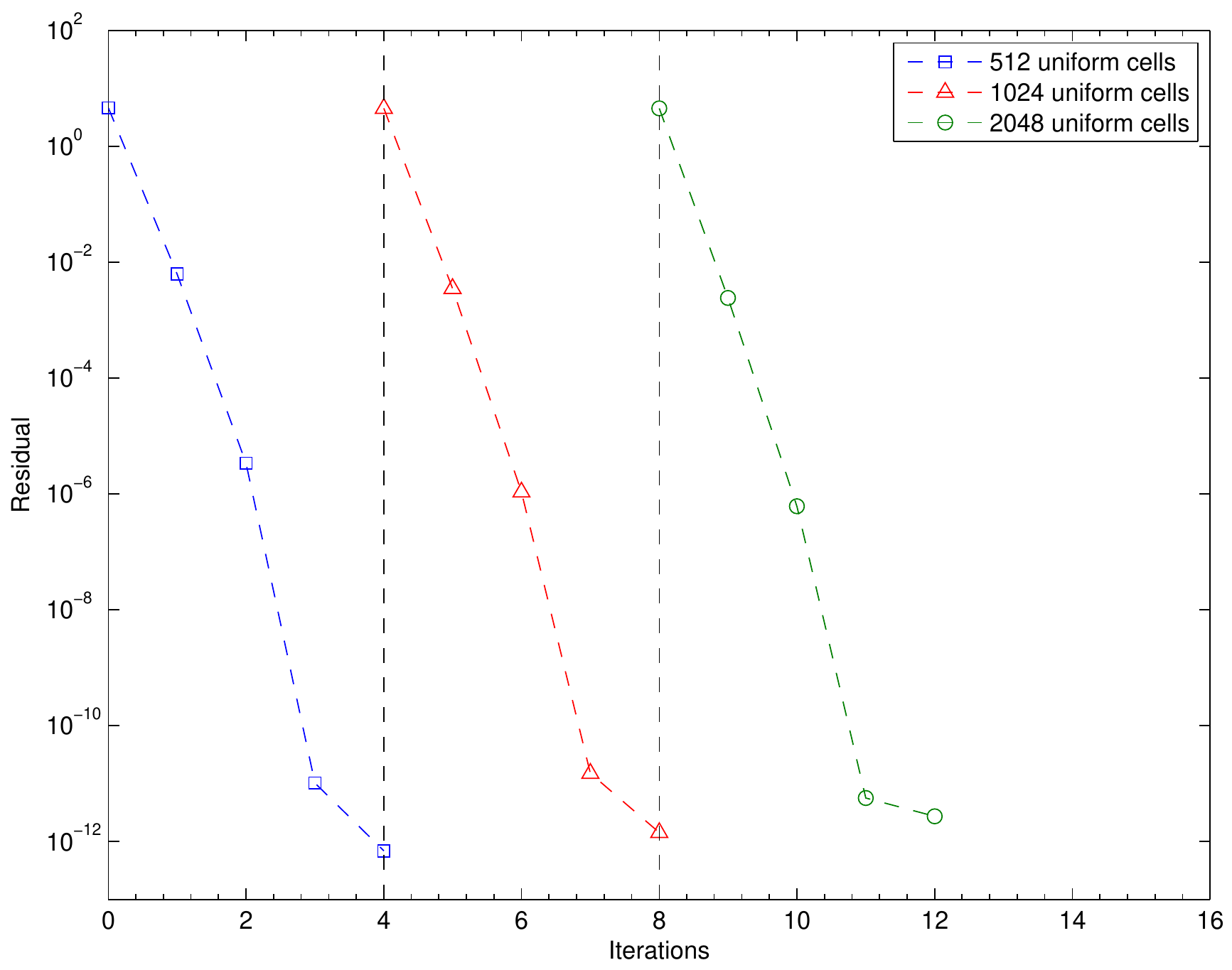}
	\includegraphics[width=0.46\textwidth,height=5cm]{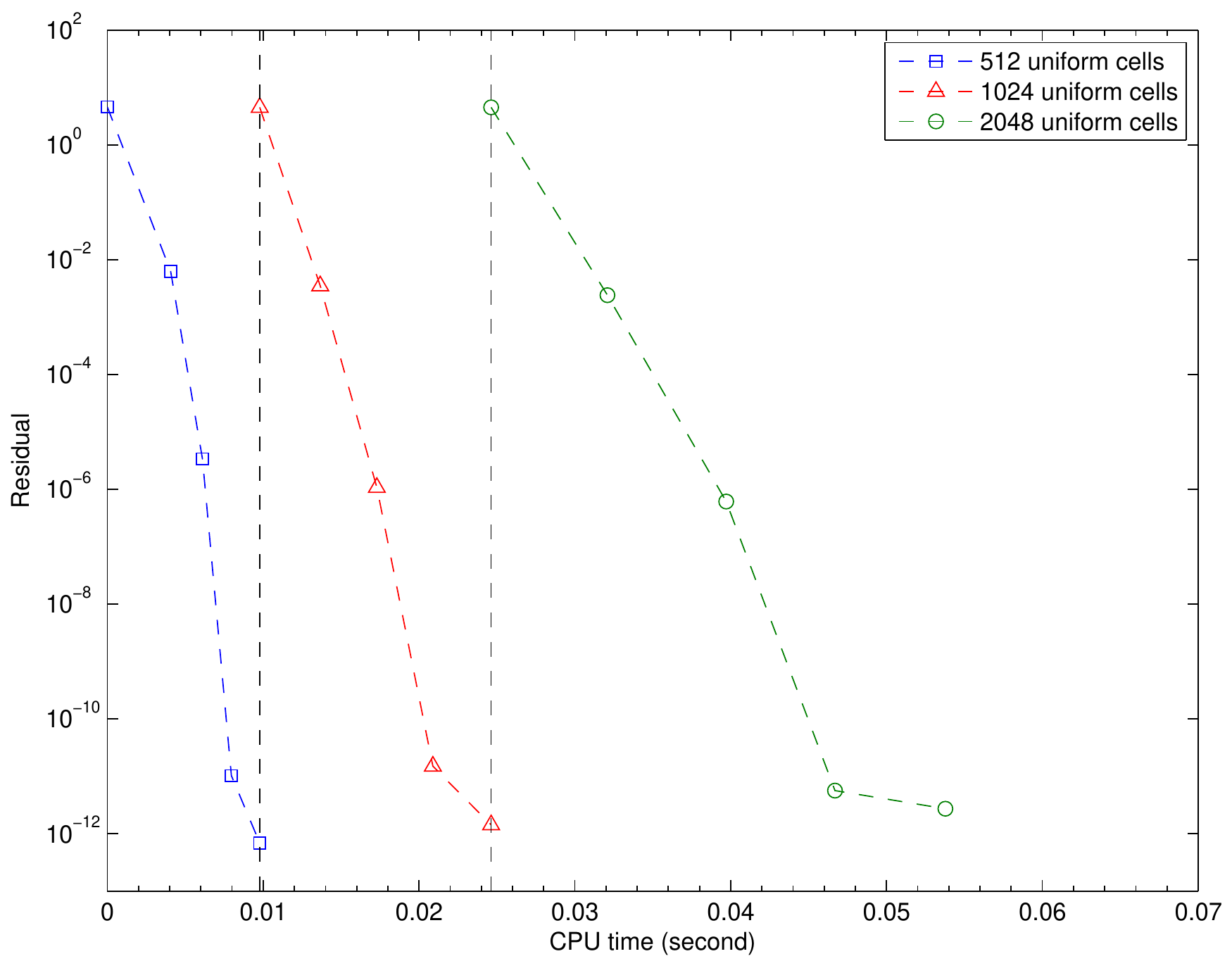}
	\caption{\small Example \ref{example1Dsmooth}:
		Residuals vs {\tt NMGM} iteration number (left) and CPU time (right) for  three uniform meshes.
	}\label{fig:1d_02_Re}
\end{figure}

\begin{table}[!h] %\multirow{3}{2pt}{$N$}
\centering
\caption{\small Example \ref{example1Dsmooth}:  Effect of the HLL, LLF, and Roe fluxes on
	convergence behaviors of
{\tt NMGM} and   the BLU-SGS iteration.
  }
\scriptsize
\begin{tabular}{|c|c|c|c|c|c|c|c|c|}
\hline
\multicolumn{3}{|c|}{$N$} & 64 & 128 & 256 & 512 & 1024 & 2048 \\
\hline
\multirow{8}{20pt}{HLL}
& Block & { $N_{\rm step}$} & 33 & 39  & 40 & 47 & 49 & 58  \\
& LU-SGS & { $T_{\rm cpu}$}  &6.57e-3 & 1.25e-2  & 1.68e-2 & 3.29e-2 & 4.87e-2 & 1.0e-1  \\
\cline{2-9}
& V-cycle & { $N_{\rm step}$}  & 4 & 4  & 3 & 4 & 4 & 4  \\
& {$N_L = 1$} & { $T_{\rm cpu}$}  & 2.51e-3 & 4.81e-3  & 5.37e-3 & 9.01e-3 & 1.37e-2 & 2.63e-2  \\
\cline{2-9}
& V-cycle & { $N_{\rm step}$}  & 4 & 4  & 3 & 4 & 4 & 4  \\
& {$N_L = 3$} & { $T_{\rm cpu}$}   & 2.75e-3 & 4.83-3  & 6.04e-3 & 9.78e-3 & 1.49e-2 & 2.91e-2  \\
%\cline{2-9}
%& V-cycle & { $N_{\rm step}$}  & x & x  & x & x & x & x  \\
%& {$N_L = 6$} & { $T_{\rm cpu}$}  & x & x  & x & x & x & x  \\
\cline{2-9}
& \multicolumn{2}{|c|}{$\rho$}   & 0.39347 & 0.43938  & 0.46526
 & 0.47874 & 0.48541 & 0.48877 \\
\cline{2-9}
& \multicolumn{2}{|c|}{$R_{\infty}$}  & 0.93276 & 0.82239  & 0.76517
 & 0.73660 & 0.72276 & 0.71586  \\
\hline
\multirow{12}{20pt}{LLF}
& Block & { $N_{\rm step}$} & 576 & 1370  & 3156 & 7116 & 15437 & 31855  \\
& LU-SGS & { $T_{\rm cpu}$}  & 4.18e-2 & 1.65e-1  & 5.44e-1 & 2.66e0 & 1.08e1 & 4.22e1  \\
\cline{2-9}
& V-cycle & { $N_{\rm step}$}  & 22 & 51  & 118 & 269 & 585 & 1224  \\
& {$N_L = 1$} & { $T_{\rm cpu}$}   & 6.21e-3 & 2.65e-2  & 8.56e-2 & 3.67e-1 & 1.55e0 & 6.33e0  \\
\cline{2-9}
& V-cycle & { $N_{\rm step}$}  & 8 & 9  & 23 & 55 & 125 & 267  \\
& {$N_L = 3$} & { $T_{\rm cpu}$}  & 4.93e-3 & 9.67e-3  & 1.90e-2 & 8.45e-2 & 3.69e-1 & 1.54e0  \\
\cline{2-9}
& V-cycle & { $N_{\rm step}$}  & 7 & 8  & 8 & 12 & 26 & 64  \\
& {$N_L = 5$} & { $T_{\rm cpu}$}   & 2.17e-3 & 6.23e-3  & 1.04e-2 & 1.91e-2 & 7.93e-2 & 3.63e-1  \\
\cline{2-9}
& W-cycle & { $N_{\rm step}$}  & 6 & 6  & 7 & 7 & 8 & 8  \\
& {$N_L = 5$} & { $T_{\rm cpu}$}  & 6.67e-3 & 9.87e-3  & 1.39e-2 & 2.38e-2 & 4.94e-2 & 9.95e-2  \\
\cline{2-9}
& \multicolumn{2}{|c|}{$\rho$}   & 0.96013 & 0.98273 & 0.99256 & 0.99673 & 0.99851 & 0.99929  \\
\cline{2-9}
& \multicolumn{2}{|c|}{$R_{\infty}$}   & 4.069e-2 & 1.742e-2  & 7.468e-3
 & 3.274e-3 & 1.493e-3 & 7.072e-4  \\
\hline
\multirow{8}{20pt}{ROE}
& Block & { $N_{\rm step}$}  & 35 &  38 & 39 & 43 & 45 & 51  \\
& LU-SGS & { $T_{\rm cpu}$}   & 6.07e-3 & 9.85e-3  & 2.18e-2 & 3.74e-2 & 5.58e-2 & 1.19e-1  \\
\cline{2-9}
& V-cycle & { $N_{\rm step}$}  & 4 & 5  & 5 & 5 & 5 & 5  \\
& {$N_L = 1$} & { $T_{\rm cpu}$}   & 3.06e-3 & 7.68e-3  & 9.57e-2 & 1.31e-2 & 2.06e-2 & 4.12e-2  \\
\cline{2-9}
& V-cycle & { $N_{\rm step}$}  & 4 & 5  & 5 & 5 & 5 & 5  \\
& {$N_L = 3$} & { $T_{\rm cpu}$}   & 3.23e-3 & 7.79e-3  & 1.05e-1 & 1.42e-2 & 2.20e-2 & 4.29e-2  \\
%\cline{2-9}
%& V-cycle & { $N_{\rm step}$}  & x & x  & x & x & x & x  \\
%& {$N_L = 5$} & { $T_{\rm cpu}$}    & x & x  & x & x & x & x  \\
\cline{2-9}
& \multicolumn{2}{|c|}{$\rho$}  & 0.48516 & 0.49799  & 0.50571
 & 0.50827 & 0.50655 & 0.49969  \\
\cline{2-9}
& \multicolumn{2}{|c|}{$R_{\infty}$}   & 0.72327 & 0.69717
& 0.68178 & 0.67673 & 0.68012 & 0.69377  \\
\hline
\end{tabular}\label{tab:1Dsmooth}
\end{table}

\begin{figure}[!htbp]
  \centering
\subfigure[]
      { \includegraphics[width=0.46\textwidth,height=5cm]{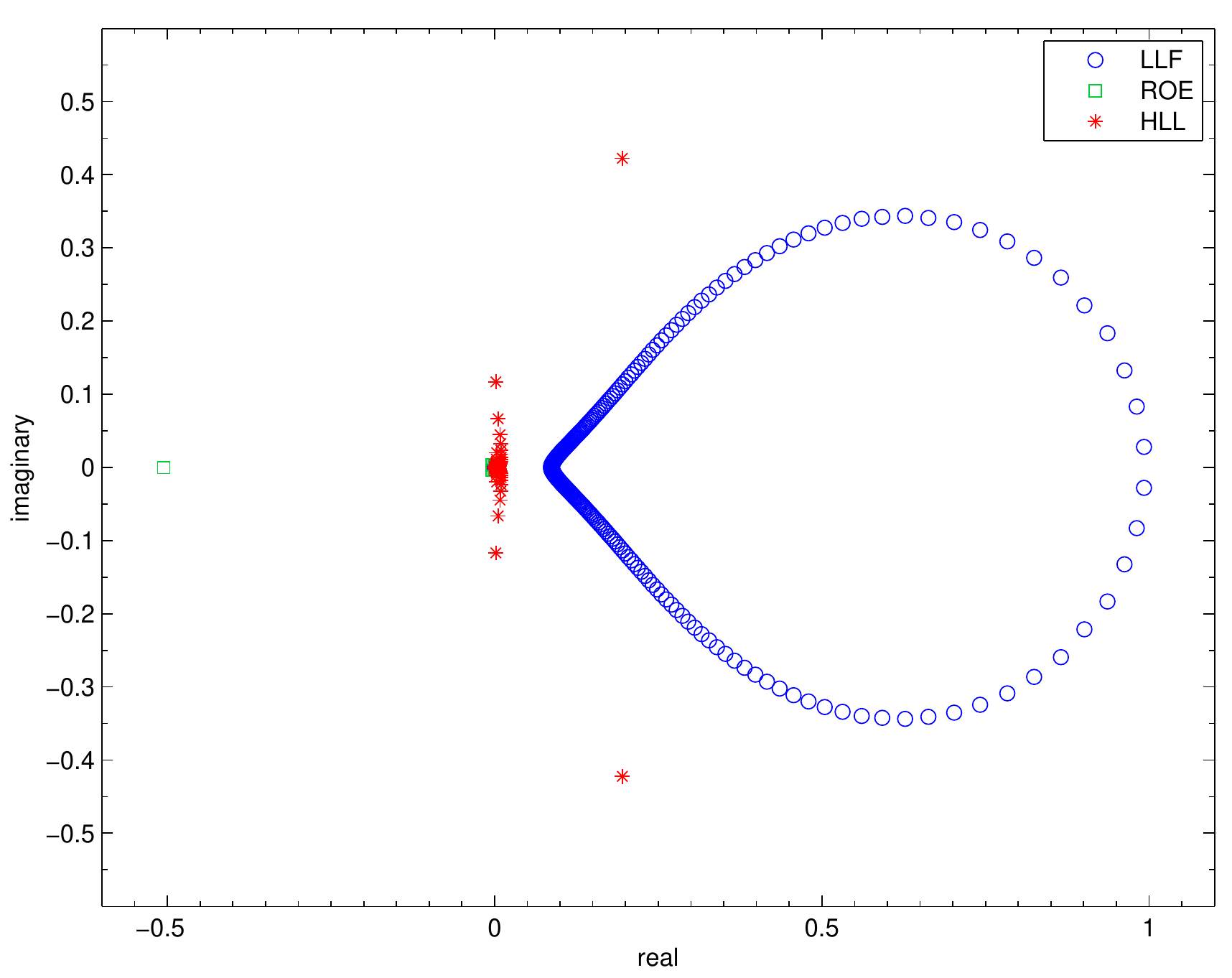} }
\subfigure[]
      { \includegraphics[width=0.46\textwidth,height=5cm]{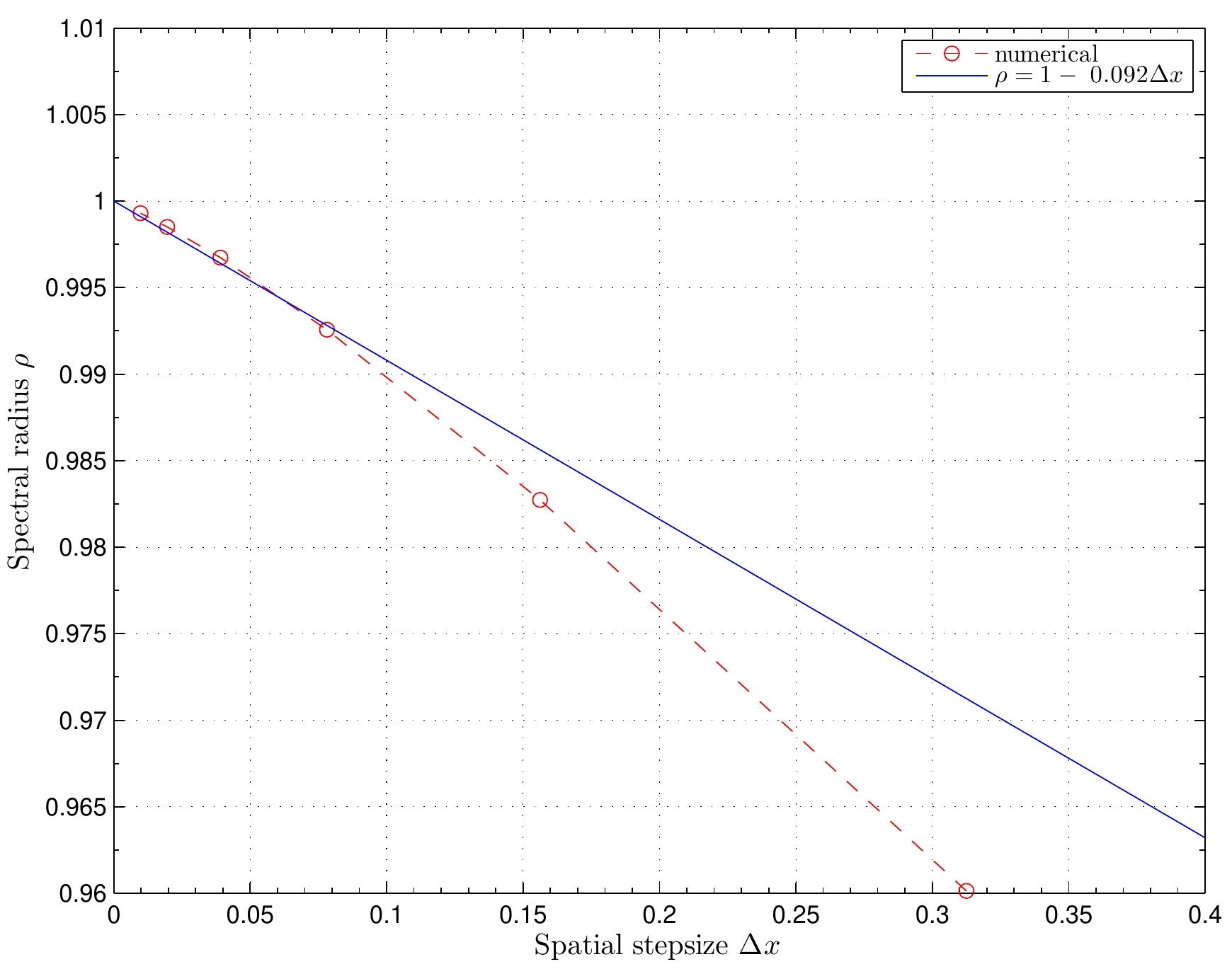} }
  \caption{\small Example \ref{example1Dsmooth}: (a). Distribution of the eigenvalues
of the block SGS iteration matrix
for three numerical fluxes over the mesh of 256 uniform cells.
(b). Asymptotic relation of the spectral radius $\rho$
of the block SGS iteration matrix
with respect to the spatial step size $\Delta x$ for the LLF flux.
}\label{fig:1d_02_eig}
\end{figure}

%\begin{figure}[!htbp]
%  \centering
%\includegraphics[width=0.46\textwidth,height=5cm]{fig1D/1D02/EigEx1D02}
%\includegraphics[width=0.46\textwidth,height=5cm]{fig1D/1D02/EigEx1D02-512}
%\caption{\small Example \ref{example1Dsmooth}: Eigenvalues
%(in the complex plane) of the block SGS iteration matrix
%for the first-order accurate scheme with three different numerical fluxes on 256 (left) and 512 (right) uniform cells,
%respectively.
%}\label{fig:1d_02_eig2}
%\end{figure}
%
%\begin{figure}[!htbp]
%  \centering
% \includegraphics[width=0.46\textwidth,height=5cm]{fig1D/1D02/2ndEigEx1D02-256}
% \includegraphics[width=0.46\textwidth,height=5cm]{fig1D/1D02/2ndEigEx1D02-512}
%  \caption{\small Same as Fig. \ref{fig:1d_02_eig2} except for the
%  second-order accurate scheme.}
%  \label{fig:1d_02_eig2b}
%\end{figure}

Table \ref{tab:1Dsmooth} investigates the effect of  the HLL, LLF,  and Roe fluxes shown in Remark \ref{Rem2.2}
on the convergence behavior  of  {\tt NMGM},  in comparison to {the BLU-SGS iteration}
 in solving \eqref{eq:NewtonREG},
where $N$, $N_L$, $N_{\rm step}$,  and   $T_{\rm cpu}$ denote the cell number,
 the {coarse mesh} level number of the geometric MG, the total Newton iteration number,
and the   CPU time $T_{\rm cpu}$, respectively,
$\rho$ denotes the spectral radius of the iteration matrix % $\vec M_{\rm SGS}(\vec U)$
of the block SGS method around the steady-state solution, %i.e. $\vec U \approx \vec U^* $,
and $R_{\infty} := - \ln  \rho $  corresponds to the asymptotic convergence rate.
The results show that  {\tt NMGM} is a little more efficient
than  the BLU-SGS iteration for the HLL and Roe fluxes,
and exhibits great advantages for the LLF flux,
specially, in the case of the W-cycle multigrid with  $N_L=5$.
The W-cycle multigrid is also tested for the HLL  and Roe fluxes,
and its performance is almost the same as the V-cycle multigrid.
%  it shows no improvement in efficiency for this problem.  %substantial
%
The BLU-SGS iteration with the HLL and Roe fluxes are more efficient
than the  LLF flux. For the former,  $N_{\rm step}$ increases very slowly with refining the mesh,
but the latter  does nearly linearly increase in terms of the cell number $N$.
Such phenomenon  could be explained by comparing  the spectral radius $\rho$ of the iteration matrix,
whose values are around  $0.5$ and larger than 0.9, respectively.
%For HLL flux and Roe  flux, the corresponding spectral radius $\rho$ keeps around $0.5$ with the increase of $N$,
%but for LLF flux, it keeps larger than $0.9$ and increases asymptotically to $1$ as $N$ increases.
%
%
Fig. \ref{fig:1d_02_eig}(a) displays
the distribution (in the complex plane) of the eigenvalues of the block SGS iteration matrix
on the mesh with 256 uniform cells for the HLL, LLF, and Roe fluxes.
We see that the eigenvalues are almost around $0$ for the HLL flux
(except for two eigenvalues with relatively large imaginary parts of $\pm 0.4$ respectively)
and the Roe flux (except for
only one eigenvalue located around $-0.5$), while those are widely distributed for the LLF flux
({some of them are near} $1$).
Thus   the ``high frequency'' eigenvalues for the LLF flux are much more
 than   for the HLL and Roe fluxes.
Fig. \ref{fig:1d_02_eig}(b) plots an
asymptotic relation of the spectral radius $\rho$
of the iteration matrix
with respect to the spatial step size $\Delta x$ for the LLF flux as follows
%
%the following asymptotic relation of $\rho$ with respect to the spatial stepsize $\Delta x$ can be observed
\begin{equation}\label{eq-04-01}
\rho \sim  1 - C \Delta x , \quad {\rm as}~ \Delta x \to 0,
\end{equation}
where $C$ is about $0.092$.

\end{example}

%%%%%%%%%%%%%%%%%%%%%%%%%%%%%%%%% Example 1D 01 A %%%%%%%%%%%%%%%%%%%%%%%%%%%%%%%%%%%%%%%%

\begin{example}[Two transcritical flows over a bump]
\label{example1Dtrans}\rm
The example simulates two transcritical flows over a bump,
which  have been widely used to test  the {SWEs} solvers
\cite{LiChen2006,Xu2002,Vazquez1999}. The bed topography is
\begin{equation}\label{eq:bottom1}
z(x)=
\begin{cases}
0.2 - 0.05(x-10)^2, &  8 < x < 12,\\
0,&    \rm{otherwise},\end{cases}
\end{equation}
for a channel of length 25.
%Depending on different boundary conditions, the flow can be subcritical, transcritical
%with or without a steady-state shock, or supercritical.
%The transcritical cases are only considered here since these two cases are more
%difficult.

{\bf (I). Transcritical flow without a shock}:
The discharge $hu = 1.53$ is imposed
at the upstream boundary condition $x = 0$£¬ while
the water height $h = 0.66$ is imposed at the downstream
end of the channel $x = 25$ when the flow is sub-critical.

 Fig. \ref{fig:1d_01a_solu} shows the numerical
 steady-state solution $h+z$  obtained by using  the HLL flux and the V-cycle multigrid with $N_L=3$
 on the mesh of 512 uniform cells in comparison to the exact solution
 given by {\tt SWASHES} \cite{Delestre2013}.
Fig. \ref{fig:1d_01a_Re} gives the convergence history of {\tt NMGM}
 in terms of the  {\tt NMGM} iteration number  $N_{\rm step}$  and CPU time on three
meshes of 512, 1024, and 2048 uniform cells respectively.
The results show that
the steady-state solutions can be correctly and
fast obtained  by using {\tt NMGM}, the convergence behaviors are
similar to Example \ref{example1Dsmooth}, and the  {\tt NMGM} iteration number does not increase
with the mesh refinement.

\begin{figure}[htbp]
	\centering
	\includegraphics[width=0.55\textwidth]{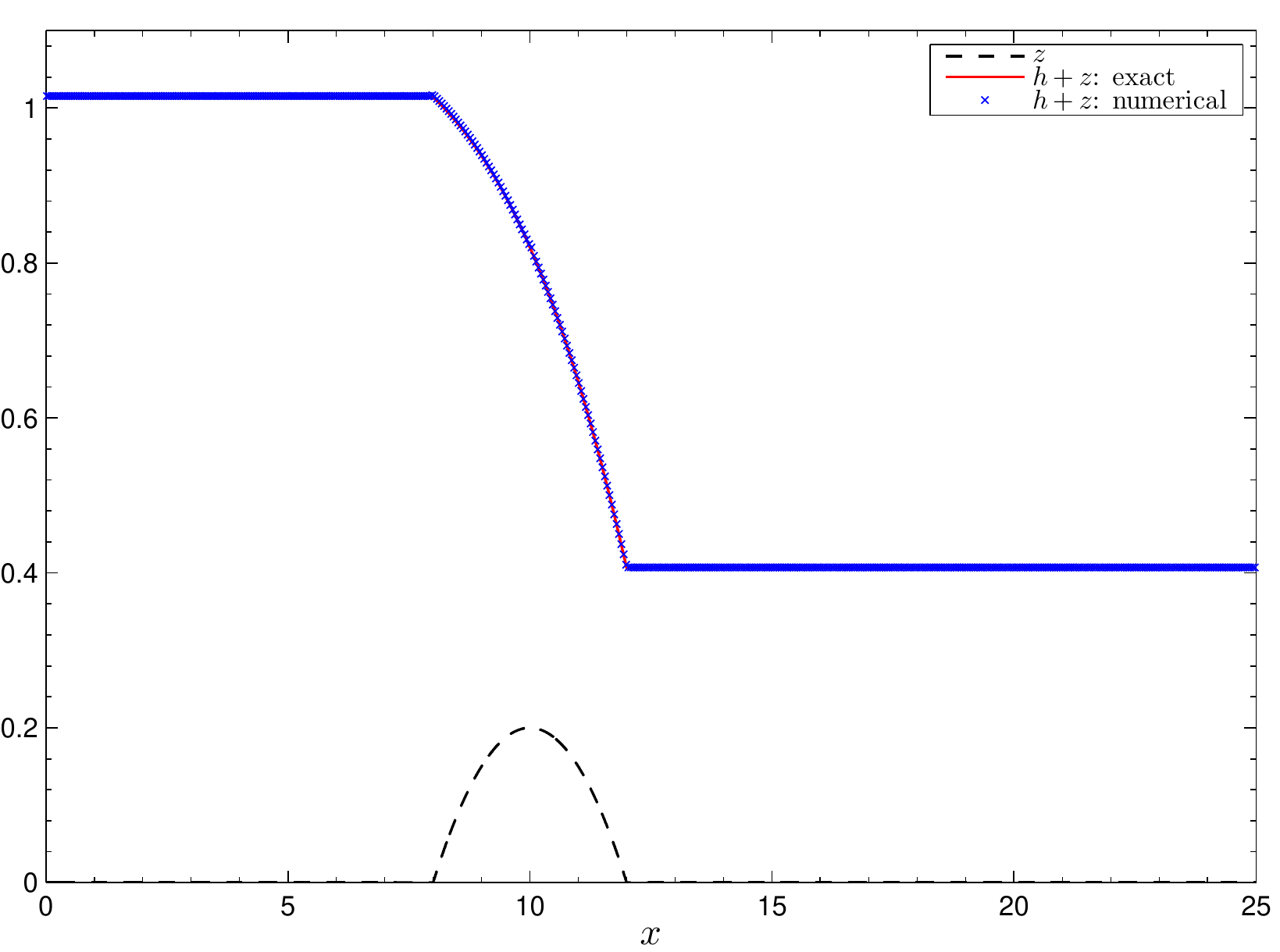}
	\caption{\small Example \ref{example1Dtrans}(I):
		Steady-state solution $h+z$ obtained by {\tt NMGM} on the mesh of 512 uniform cells.
	}\label{fig:1d_01a_solu}
\end{figure}

\begin{figure}[htbp]
  \centering
  \includegraphics[width=0.46\textwidth,height=5cm]{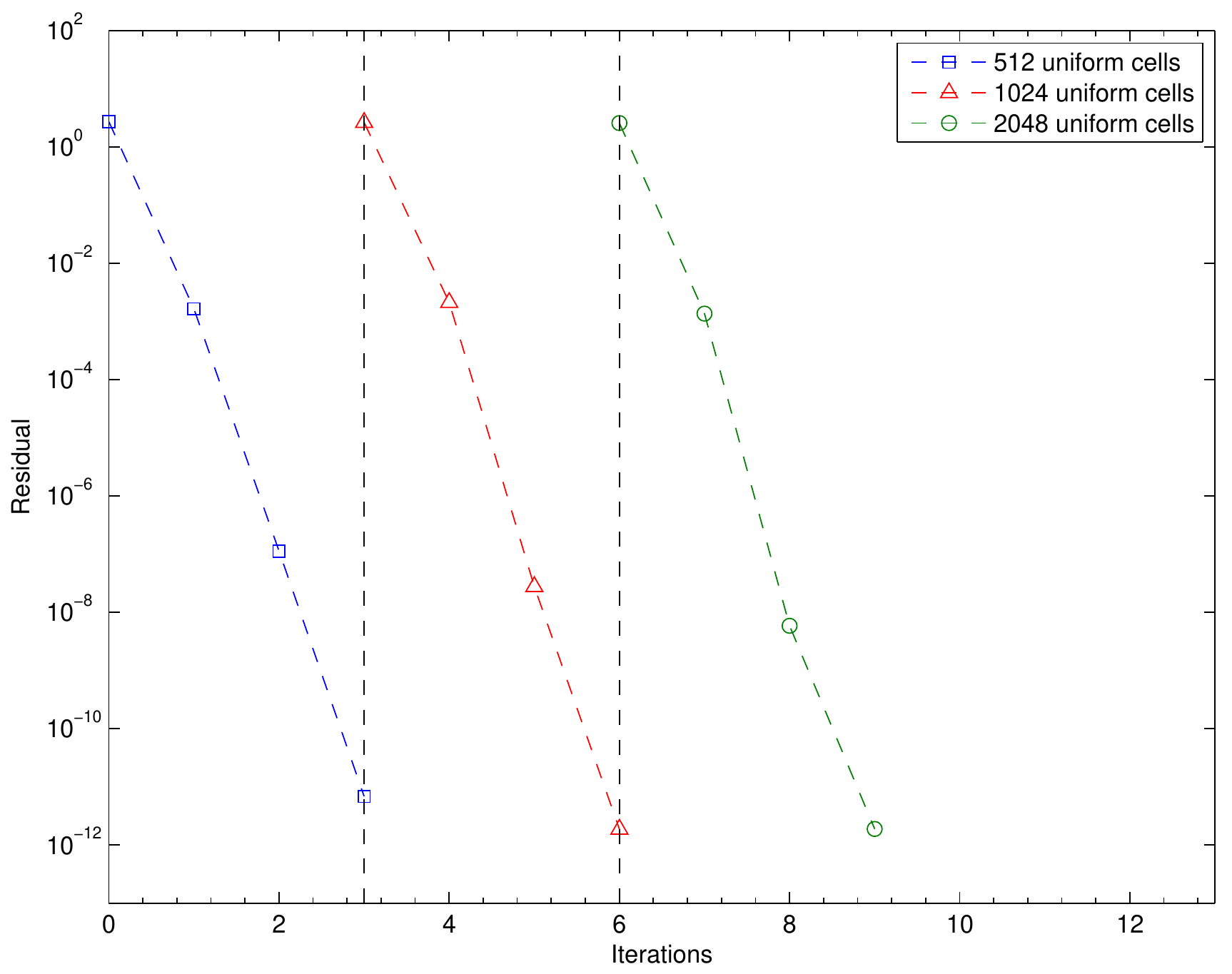}
  \includegraphics[width=0.46\textwidth,height=5cm]{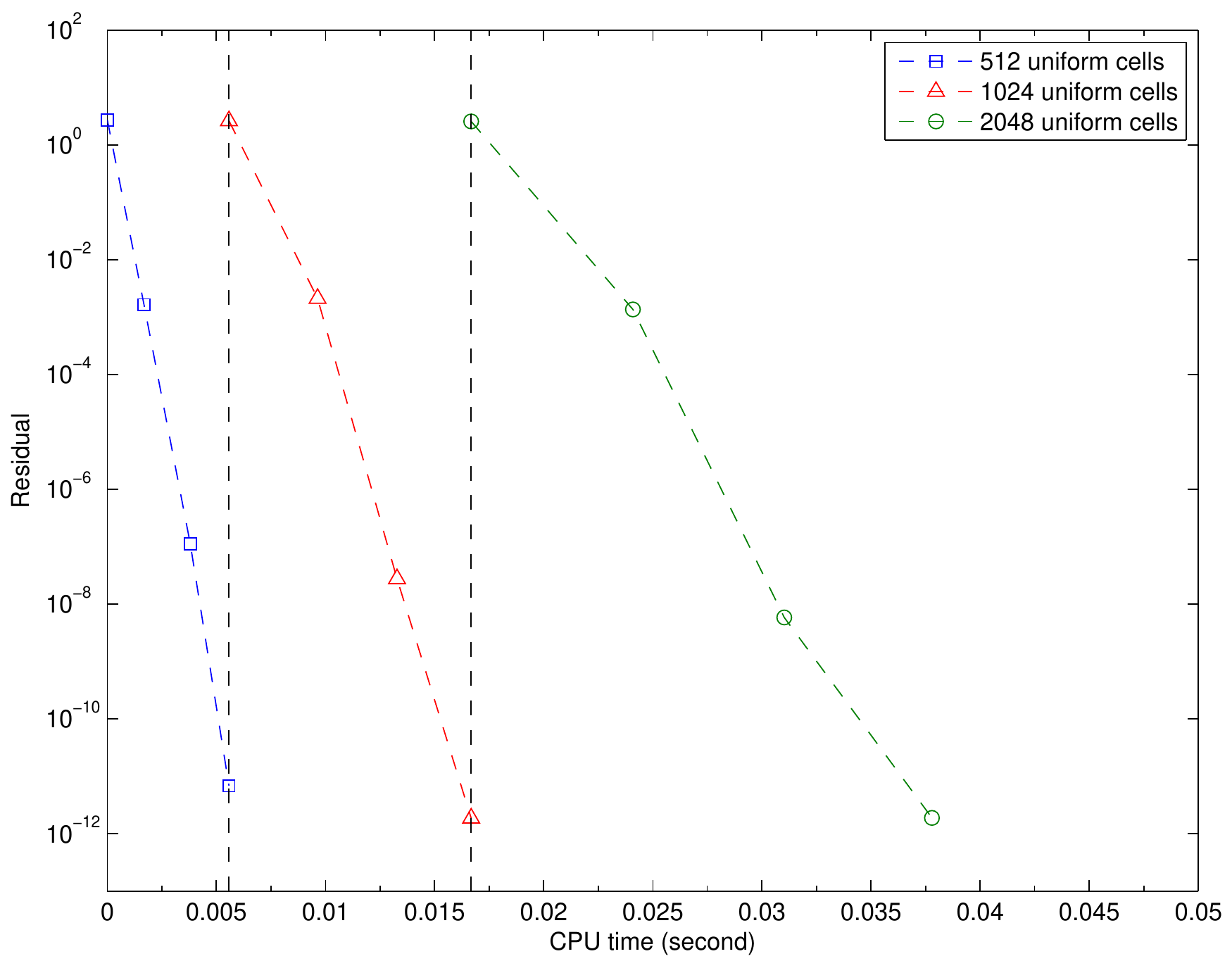}
  \caption{\small Example \ref{example1Dtrans}(I):
Convergence history in terms of the  {\tt NMGM} iteration number  $N_{\rm step}$  (left) and CPU time (right) on three uniform meshes.
}\label{fig:1d_01a_Re}
\end{figure}

\begin{table}[htbp] %\multirow{3}{2pt}{$N$}
\centering
\caption{\small Example \ref{example1Dtrans}(I):
%Convergence behaviors of
%{\tt NMGM} and block SGS method with different numerical fluxes.
 Effect of the HLL, LLF, and Roe fluxes on
 convergence behaviors of
 {\tt NMGM} and   the BLU-SGS iteration.
  }
\scriptsize
\begin{tabular}{|c|c|c|c|c|c|c|c|c|}
\hline
\multicolumn{3}{|c|}{$N$} & 64 & 128 & 256 & 512 & 1024 & 2048 \\
\hline
\multirow{8}{20pt}{HLL}
& Block & { $N_{\rm step}$} & 25 & 28  & 31 & 31 & 32 & 33  \\
& LU-SGS & { $T_{\rm cpu}$}   & 8.87e-3 & 1.11e-2  & 1.43e-2 & 1.76e-2 & 2.85e-2 & 5.21e-2  \\
\cline{2-9}
& V-cycle & { $N_{\rm step}$}  & 3 & 3  & 3 & 3 & 3 & 3  \\
& {$N_L = 1$} & { $T_{\rm cpu}$}  & 1.88e-3 & 2.24e-3  & 4.52e-3 & 8.54e-3 & 1.24e-2 & 2.07e-2  \\
\cline{2-9}
& V-cycle & { $N_{\rm step}$}  & 3 & 3  & 3 & 3 & 3 & 3  \\
& {$N_L = 3$} & { $T_{\rm cpu}$}   & 1.96e-3 & 3.31e-3  & 6.44e-3 & 9.16e-2 & 1.30e-2 & 2.09e-2  \\
\cline{2-9}
& \multicolumn{2}{|c|}{$\rho$}   & 0.36309 & 0.38504  & 0.41148 & 0.47183 & 0.4518 & 0.46042  \\
\cline{2-9}
& \multicolumn{2}{|c|}{$R_{\infty}$}   & 1.0131 & 0.9544  & 0.8879 & 0.7511 & 0.7945 & 0.7756  \\
\hline
\multirow{12}{20pt}{LLF}
& Block & { $N_{\rm step}$}  & 373 & 481  & 878 & 1667 & 3277 & 6462 \\
& LU-SGS & { $T_{\rm cpu}$}   & 3.67e-2 & 5.65e-2  & 1.74e-1 & 6.17e-1 & 2.31e0 & 8.73e0  \\
\cline{2-9}
& V-cycle & { $N_{\rm step}$}  & 13 & 20 & 35 & 66 & 128 & 249  \\
& {$N_L = 1$} & { $T_{\rm cpu}$} & 7.89e-3 & 1.14e-2 & 2.48e-1 & 9.09e-1 & 3.40e-1 & 1.27e0  \\
\cline{2-9}
& V-cycle & { $N_{\rm step}$}  & 14 & 16  & 17 & 18 & 31 & 59  \\
& {$N_L = 3$} & { $T_{\rm cpu}$}   & 9.42e-3 & 1.02e-2  & 1.35e-2 & 2.78e-2 & 9.12e-2 & 3.43e-1  \\
\cline{2-9}
& V-cycle & { $N_{\rm step}$}  & 14 & 17  & 18 & 19 & 20 & 22  \\
& {$N_L = 5$} & { $T_{\rm cpu}$}  & 9.67e-3 & 1.15e-2  & 1.51e-2 & 3.02e-2 & 6.11e-2 & 1.31e-1  \\
\cline{2-9}
& W-cycle & { $N_{\rm step}$}  & 11 & 12  & 9 & 8 & 8 & 8  \\
& {$N_L = 5$} & { $T_{\rm cpu}$}  & 6.73e-3 & 1.06e-2  & 1.48e-2 & 2.50e-2 & 4.95e-2 & 9.63e-2  \\
\cline{2-9}
& \multicolumn{2}{|c|}{$\rho$}   & 0.92458 & 0.94598  & 0.97116 & 0.98515 & 0.99246 & 0.9962  \\
\cline{2-9}
& \multicolumn{2}{|c|}{$R_{\infty}$}   & 7.842e-2 & 5.553e-2  & 2.936e-2 & 1.496e-2 & 7.566e-3 & 3.805e-3  \\
\hline
\multirow{8}{20pt}{ROE}
& Block & { $N_{\rm step}$}  & 15 & 21  & 31 & 49 & 88 & 172  \\
& LU-SGS & { $T_{\rm cpu}$}   & 4.42e-3 & 1.06e-2  & 1.68e-2 & 3.82e-2 & 1.08e-1 & 4.10e-1  \\
\cline{2-9}
& V-cycle & { $N_{\rm step}$}  & 3 & 3  & 3 & 3 & 4 & 5  \\
& {$N_L = 1$} & { $T_{\rm cpu}$}   & 2.68e-3 & 5.01e-3  & 8.33e-3 & 1.16e-2 & 1.82e-2 & 4.52e-2  \\
\cline{2-9}
& V-cycle & { $N_{\rm step}$}  & 3 & 3  & 3 & 3 & 4 & 4  \\
& {$N_L = 3$} & { $T_{\rm cpu}$}   & 2.70e-3 & 5.10e-3  & 8.54e-3 & 1.36e-2 & 1.85e-2 & 3.63e-2  \\
%\cline{2-9}
%& V-cycle & { $N_{\rm step}$}  & x & x  & x & x & x & x  \\
%& {$N_L = 5$} & { $T_{\rm cpu}$}    & x & x  & x & x & x & x  \\
\cline{2-9}
& \multicolumn{2}{|c|}{$\rho$}  & 0.16498 & 0.33576  & 0.53079 & 0.70641
 & 0.83135 & 0.90858  \\
\cline{2-9}
& \multicolumn{2}{|c|}{$R_{\infty}$}   & 1.802e0 & 1.091e0  & 6.334e-1
& 3.476e-1 & 1.847e-1 & 9.587e-2  \\
\hline
\end{tabular}\label{tab:1Dtrans:a}
\end{table}

\begin{figure}[htbp]
  \centering
\subfigure[]
      { \includegraphics[width=0.46\textwidth,height=5cm]{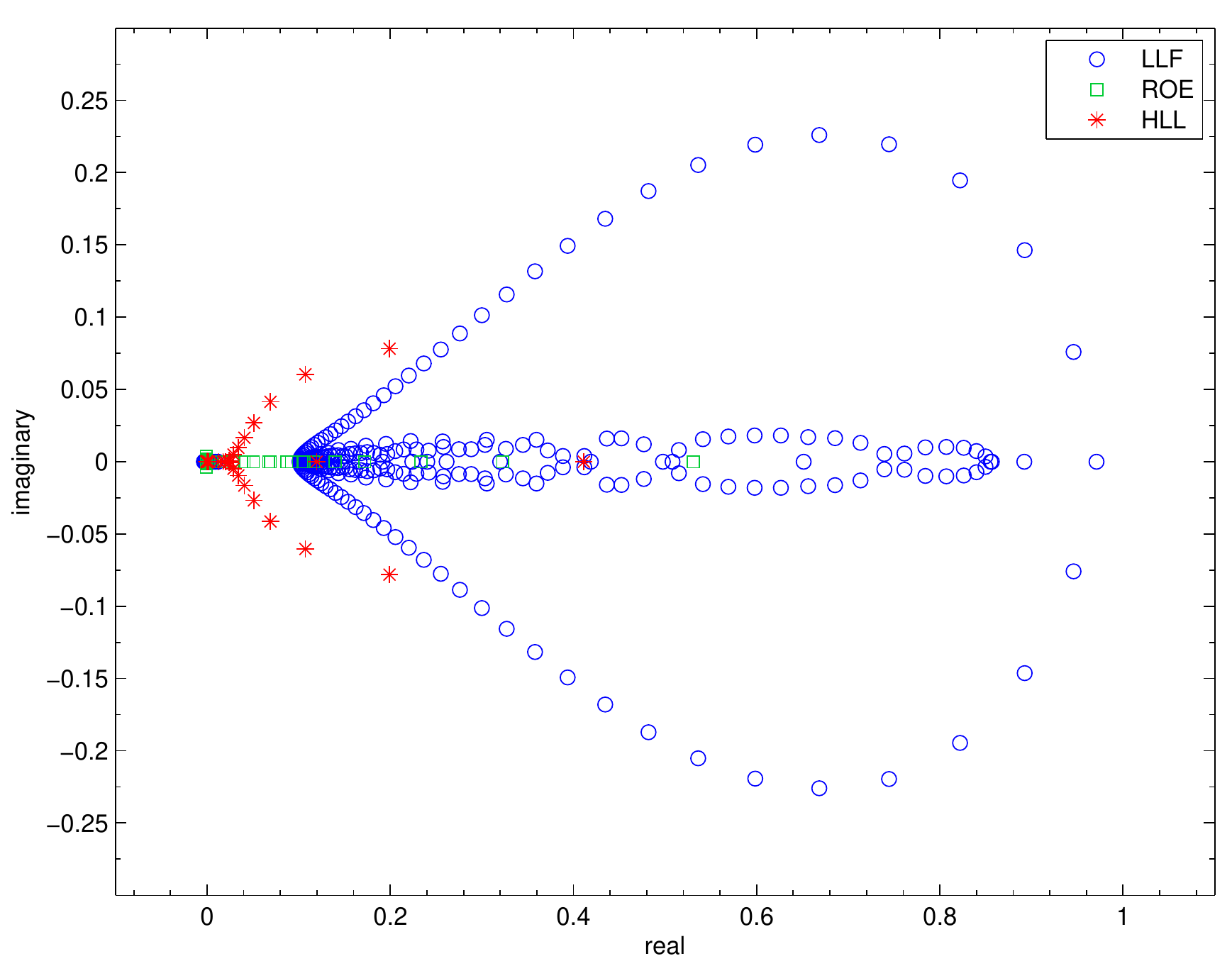} }
\subfigure[]
      { \includegraphics[width=0.46\textwidth,height=5cm]{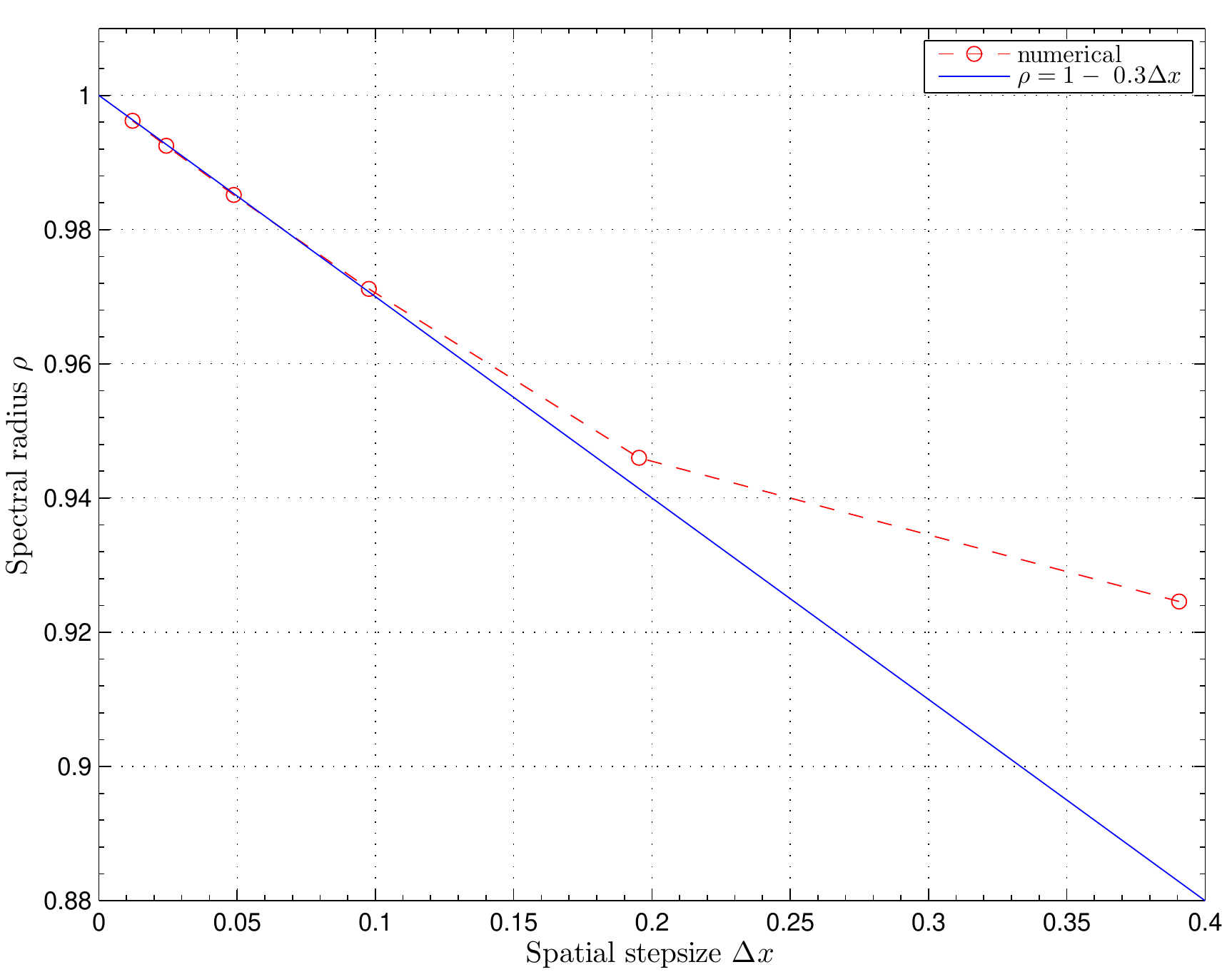} }
  \caption{\small Example \ref{example1Dtrans}(I): (a). Distribution of the eigenvalues
  of the block SGS iteration matrix
with the HLL, LLF, Roe fluxes over the mesh of 256 uniform cells.
(b). Asymptotic relation of the spectral radius $\rho$ of the  block SGS iteration matrix
with respect to the spatial stepsize $\Delta x$ for LLF flux.
}\label{fig:1d_01a_eig}
\end{figure}

%\begin{figure}[!htbp]
%  \centering
%\subfigure[1st-order accurate scheme with 256 uniform cells]
%      { \includegraphics[width=0.46\textwidth,height=5cm]{fig1D/1D01a/EigEx1D01a} }
%\subfigure[1st-order accurate scheme with 512 uniform cells]
%      { \includegraphics[width=0.46\textwidth,height=5cm]{fig1D/1D01a/EigEx1D01a-512} }
%\subfigure[2nd-order accurate scheme with 256 uniform cells]
%      { \includegraphics[width=0.46\textwidth,height=5cm]{fig1D/1D01a/2ndEigEx1D01a-256} }
%\subfigure[2nd-order accurate scheme with 512 uniform cells]
%      { \includegraphics[width=0.46\textwidth,height=5cm]{fig1D/1D01a/2ndEigEx1D01a-512} }
%  \caption{\small Example \ref{example1Dtrans}(I): Distribution of the eigenvalues
%(on the complex plane) of the iteration matrix for the block SGS method
%with different numerical fluxes.
%}\label{fig:1d_01a_eig2}
%\end{figure}

%Similar to Example \ref{example1Dsmooth}, using different numerical fluxes gives
%different convergence behaviors of  {\tt NMGM}.
Similar to Table \ref{tab:1Dsmooth},
Table \ref{tab:1Dtrans:a} investigates the effect of  the HLL, LLF,  and Roe fluxes
on the convergence behavior  of  {\tt NMGM},  in comparison to  the BLU-SGS iteration.
% gives the number of iterations $N_{\rm step}$ and
%the estimated CPU time $T_{\rm cpu}$ of {\tt NMGM} on different
%meshes by using HLL flux, LLF flux and Roe  flux, and the performance is shown
%with different number of coarser levels
%in comparing to the block SGS method without multigrid iterations.
%
The convergence behaviors of {\tt NMGM} and  the BLU-SGS iteration are almost the same as those in Example \ref{example1Dsmooth}.
The distribution (in the complex plane) of the eigenvalues of the block SGS iteration
matrix  in Fig. \ref{fig:1d_01a_eig}(a) shows that
 the ``high frequency'' eigenvalues  for the LLF flux  are much more
 than  for HLL or Roe  flux,
%given by the computing over mesh of 256 uniform cells.
%{\tt NMGM} shows a little more efficient performance than block SGS method for HLL flux and Roe  flux,
%and exhibits great advantages for LLF flux, especially when the number of levels of multigrid increases to $6$
%and the multigrid strategy of W-cycle is adopted.
%
%W-cycle is also tried for HLL flux and Roe  flux in our experiments,
%and it shows no improvement in efficiency for this problem.  %substantial
 Fig. \ref{fig:1d_01a_eig}(b) shows that
the spectral radius $\rho$ of the iteration matrix has asymptotic relation \eqref{eq-04-01} with $C\sim 0.3$ for the  LLF flux
in terms of  the spatial step size $\Delta x$.
%$$\rho \sim  1 - C \Delta x , \quad {\rm as}~ \Delta x \to 0, $$

%%%%%%%%%%%%%%%%%%%%%%%%%%%%%%%%% Example 1D 01 B %%%%%%%%%%%%%%%%%%%%%%%%%%%%%%%%%%%%%%%%

{\bf(II). Transcritical flow with a shock}:
The discharge $hu$ is taken as $0.18$ on the upstream boundary, and $h=0.33$
is specified on the downstream boundary condition.
In this case, the Froude number $F_r=u/\sqrt{gh}$
increases to a value larger than 1 above the bump, and then
decreases to less than 1.

Fig. \ref{fig:1d_01b_solu} shows the numerical
steady-state solution   on the mesh of 512 uniform cells in comparison to the exact solution
  \cite{Delestre2013}£¬ where  a stationary shock appears and is well resolved.
Fig. \ref{fig:1d_01b_Re} gives the convergence history of {\tt NMGM}
versus the  {\tt NMGM} iteration number  $N_{\rm step}$   and CPU time on three
meshes of 512, 1024, and 2048 uniform cells respectively.
The results show that
  {\tt NMGM} is very efficient and fast
to get the correct steady-state solution with a shock,
  the convergence behaviors are almost similar
on three uniform meshes, and the iteration number  $N_{\rm step}$  scarcely changes with the
mesh refinement.
The HLL flux and V-cycle multigrid  with $N_L=3$ have been adopted in those computations.
\begin{figure}[htbp]
	\centering
	\includegraphics[width=0.5\textwidth]{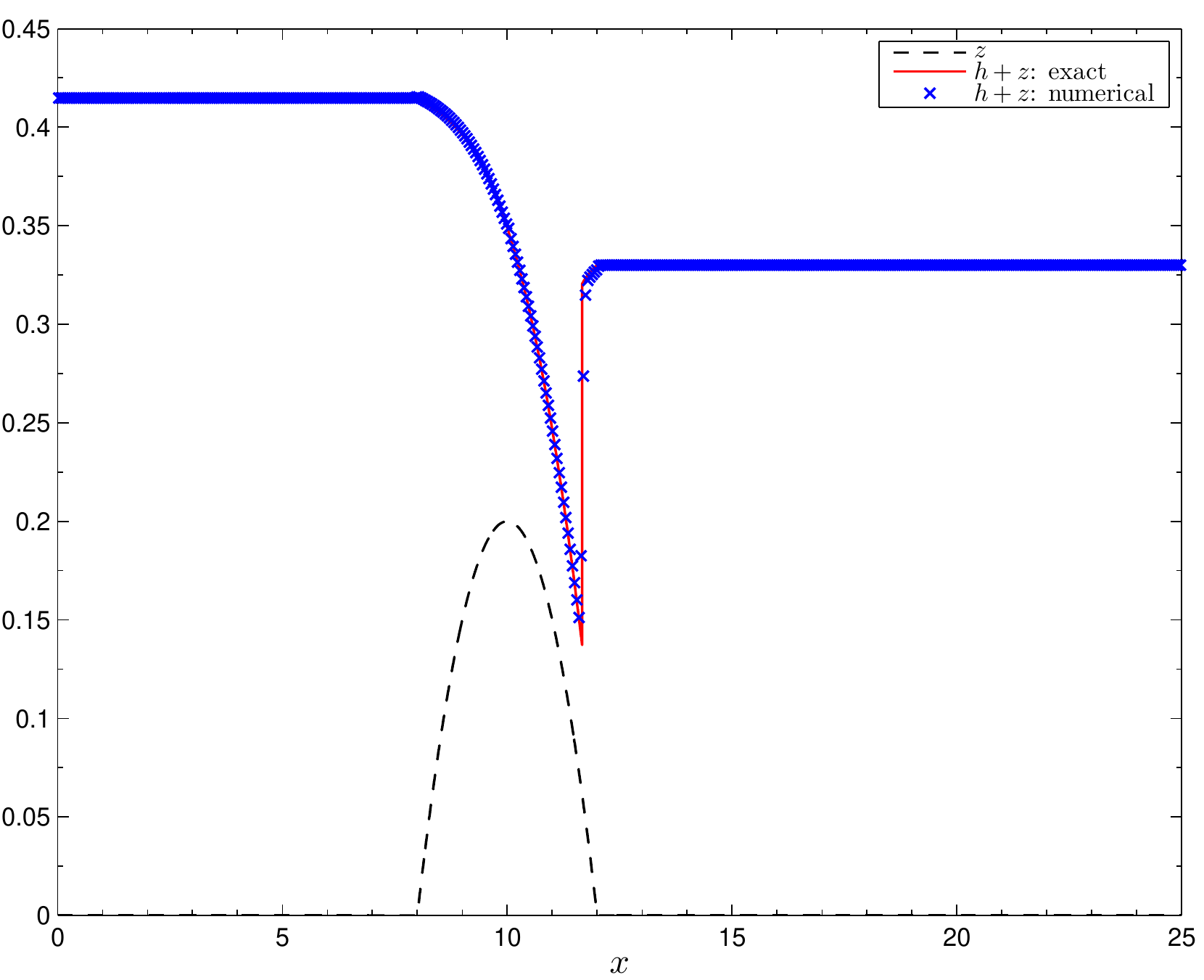}
	\caption{\small Example \ref{example1Dtrans}(II):
		Steady-state solution $h+z$ obtained by {\tt NMGM} on the mesh of 512 uniform cells.
	}\label{fig:1d_01b_solu}
\end{figure}

\begin{figure}[htbp]
  \centering
  \includegraphics[width=0.46\textwidth,height=5cm]{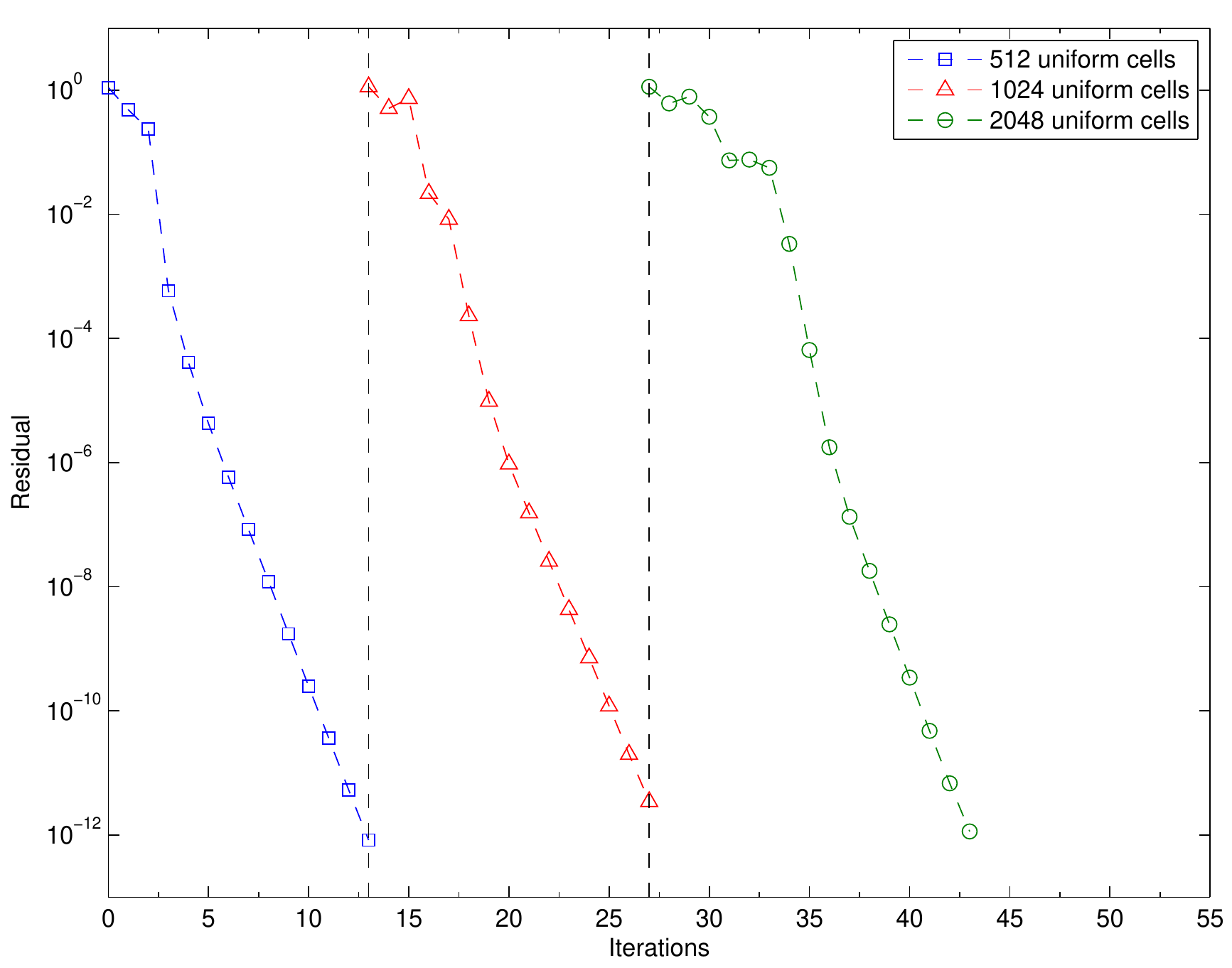}
  \includegraphics[width=0.46\textwidth,height=5cm]{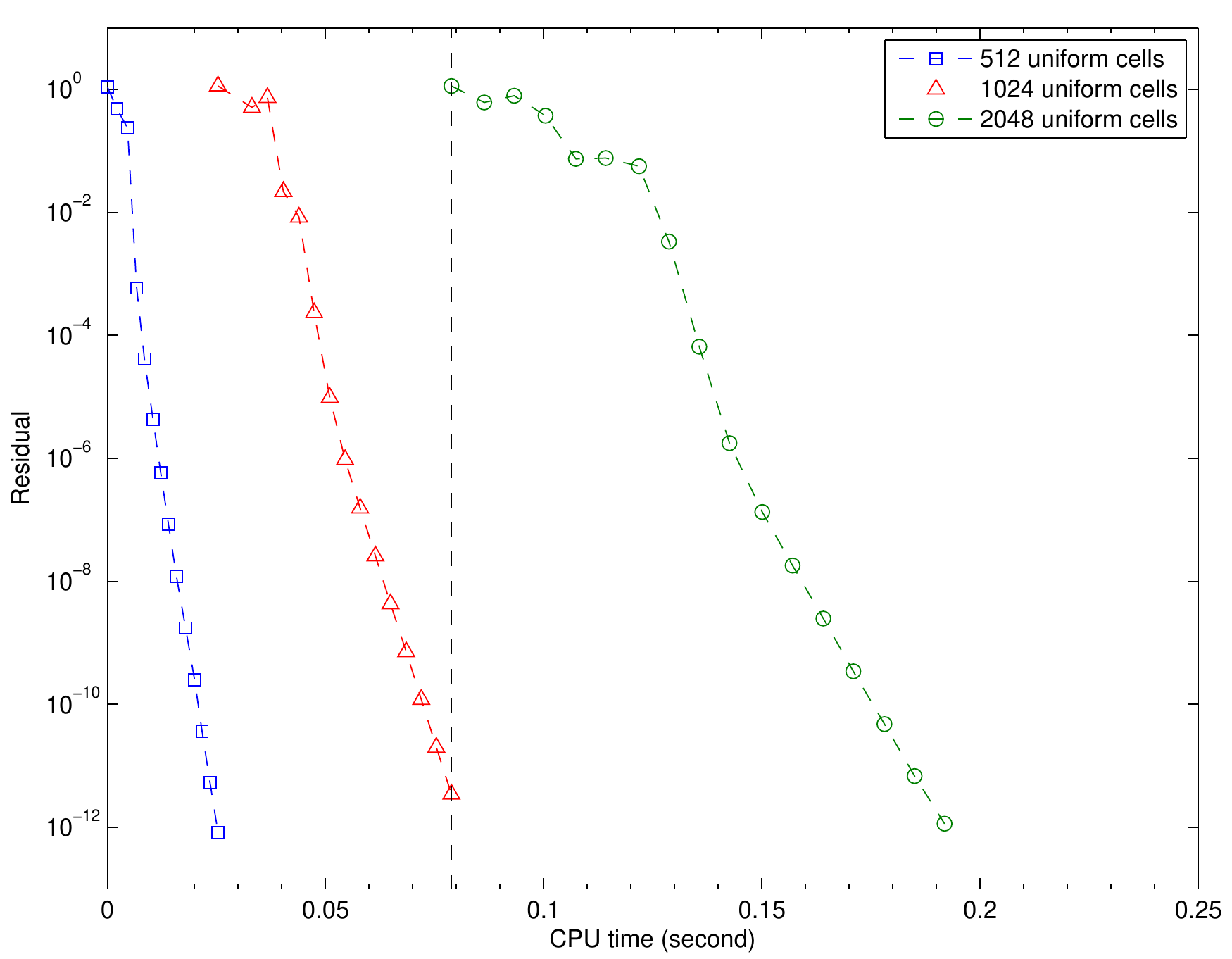}
  \caption{\small Example \ref{example1Dtrans}(II):
Convergence history in terms of  the  {\tt NMGM} iteration number  $N_{\rm step}$ (left) and CPU time (right) on three uniform meshes.
}\label{fig:1d_01b_Re}
\end{figure}

\begin{table}[htbp] %\multirow{3}{2pt}{$N$}
\centering
\caption{\small Example \ref{example1Dtrans}(II):
	%Convergence behaviors of
%{\tt NMGM} and BLU-SGS method with different numerical fluxes.
 Effect of the HLL, LLF, and Roe fluxes on
 convergence behaviors of
 {\tt NMGM} and   the BLU-SGS iteration.
  }
\scriptsize
\begin{tabular}{|c|c|c|c|c|c|c|c|c|}
\hline
\multicolumn{3}{|c|}{$N$} & 64 & 128 & 256 & 512 & 1024 & 2048 \\
\hline
\multirow{10}{20pt}{HLL}
& Block & { $N_{\rm step}$} & 28 & 37  & 82 & 161 & 328 & 658  \\
& LU-SGS & { $T_{\rm cpu}$}   & 2.87e-3 & 7.07e-3 & 2.94e-2 & 1.07e-1 & 4.88e-1 & 1.92e0 \\
\cline{2-9}
& V-cycle & { $N_{\rm step}$} & 4 & 4  & 4 & 13  & 14 & 28  \\
& {$N_L = 1$} & { $T_{\rm cpu}$}  & 2.62e-3 & 4.49e-3 & 7.99e-3 & 2.52e-2 & 4.68e-2 & 1.73e-1 \\
\cline{2-9}
& V-cycle & { $N_{\rm step}$} & 4 & 4 & 5 & 13 & 14 & 16  \\
& {$N_L = 3$} & { $T_{\rm cpu}$}   & 2.87e-3 & 5.27e-3 & 8.64e-3 & 2.54e-2 & 5.35e-2 & 1.13e-1 \\
\cline{2-9}
& V-cycle & { $N_{\rm step}$} & 4 & 4  & 5 & 13 & 14 & 16  \\
& {$N_L = 5$} & { $T_{\rm cpu}$} & 2.96e-3 & 5.32e-3 & 8.89e-3 & 3.76e-2 & 5.44e-2 & 1.14e-1 \\
\cline{2-9}
& \multicolumn{2}{|c|}{$\rho$}  & 0.39173 & 0.52614 & 0.74702 &
0.86591 & 0.93346 & 0.96696  \\
\cline{2-9}
& \multicolumn{2}{|c|}{$R_{\infty}$}  & 9.372e-1 & 6.422e-1 & 2.917e-1
 & 1.44e-1 & 6.885e-2 & 3.36e-2  \\
\hline
\multirow{12}{20pt}{LLF}
& Block & { $N_{\rm step}$} & 288 & 662  & 1424 & 3084 & 6518 & 13179  \\
& LU-SGS & { $T_{\rm cpu}$}  & 1.64e-2 & 9.59e-2 & 2.61e-1 & 1.91e0 & 6.01e0 & 2.60e1 \\
\cline{2-9}
& V-cycle & { $N_{\rm step}$} & 12 & 25 & 52 & 113 & 242 & 488  \\
& {$N_L = 1$} & { $T_{\rm cpu}$}  & 6.29e-3 & 2.37e-2 & 4.93e-2 & 1.55e-1 & 6.43e-1 & 2.49e0  \\
\cline{2-9}
& V-cycle & { $N_{\rm step}$} & 8 & 8 & 14 & 28 & 53 & 103  \\
& {$N_L = 3$} & { $T_{\rm cpu}$}  & 4.98e-3 & 9.84e-3 & 1.51e-2 & 4.73e-2 & 1.56e-1 & 5.99e-1  \\
\cline{2-9}
& V-cycle & { $N_{\rm step}$} & 8 & 9  & 10 & 10  & 17 & 47  \\
& {$N_L = 5$} & { $T_{\rm cpu}$}  & 5.13e-3 & 1.16e-2 & 1.71e-2 & 2.28e-2 & 5.27e-2 & 2.76e-1 \\
\cline{2-9}
& W-cycle & { $N_{\rm step}$} & 7 & 6  & 6 & 7  & 7 & 10  \\
& {$N_L = 5$} & { $T_{\rm cpu}$} & 7.62e-3 & 1.27e-2 & 1.73e-2 & 2.23e-2 & 4.49e-2 & 1.24e-1 \\
\cline{2-9}
& \multicolumn{2}{|c|}{$\rho$}  & 0.92518 & 0.96734 & 0.98529 & 0.99604 & 0.99663 & 0.99834  \\
\cline{2-9}
& \multicolumn{2}{|c|}{$R_{\infty}$}  & 7.777e-2 & 3.321e-2 & 1.482e-2 & 6.987e-3 & 3.378e-3 & 1.659e-3  \\
\hline
\multirow{10}{20pt}{ROE}
& Block & { $N_{\rm step}$} & 19 & 26  & 36 & 102 & 204 &  378\\
& LU-SGS & { $T_{\rm cpu}$}  & 1.52e-3 & 4.13e-3 & 1.32e-2 & 6.37e-2 & 2.45e-1 & 9.04e-1 \\
\cline{2-9}
& V-cycle & { $N_{\rm step}$} & 3 & 4  & 4 & 12 & 10 & 18  \\
& {$N_L = 1$} & { $T_{\rm cpu}$}  & 2.26e-3 & 5.83e-3 & 1.09e-2 & 2.98e-2 & 3.88e-2 & 1.34e-1 \\
\cline{2-9}
& V-cycle & { $N_{\rm step}$} & 3 & 4  & 4 & 12  & 11 & 12  \\
& {$N_L = 3$} & { $T_{\rm cpu}$}  & 2.42e-3 & 6.37e-3 & 1.22e-2 & 3.05e-2 & 4.59e-2 & 9.86e-2 \\
\cline{2-9}
& V-cycle & { $N_{\rm step}$} & 3 & 4  & 4 & 12 & 11 & 12  \\
& {$N_L = 5$} & { $T_{\rm cpu}$}   & 2.61e-3 & 6.64e-3 & 1.26e-2 & 3.10e-2 & 4.79e-2 & 1.01e-1 \\
\cline{2-9}
& \multicolumn{2}{|c|}{$\rho$}  & 0.25371 & 0.36652 & 0.54545 & 0.79807 & 0.89694 & 0.9425  \\
\cline{2-9}
& \multicolumn{2}{|c|}{$R_{\infty}$}  & 1.372e0 & 1.004e0 & 6.062e-1 & 2.256e-1 & 1.088e-1 & 5.922e-2  \\
\hline
\end{tabular}\label{tab:1Dtrans:b}
\end{table}

\begin{figure}[htbp]
  \centering
\subfigure[Eigenvalues in complex plane]
      { \includegraphics[width=0.46\textwidth,height=5cm]{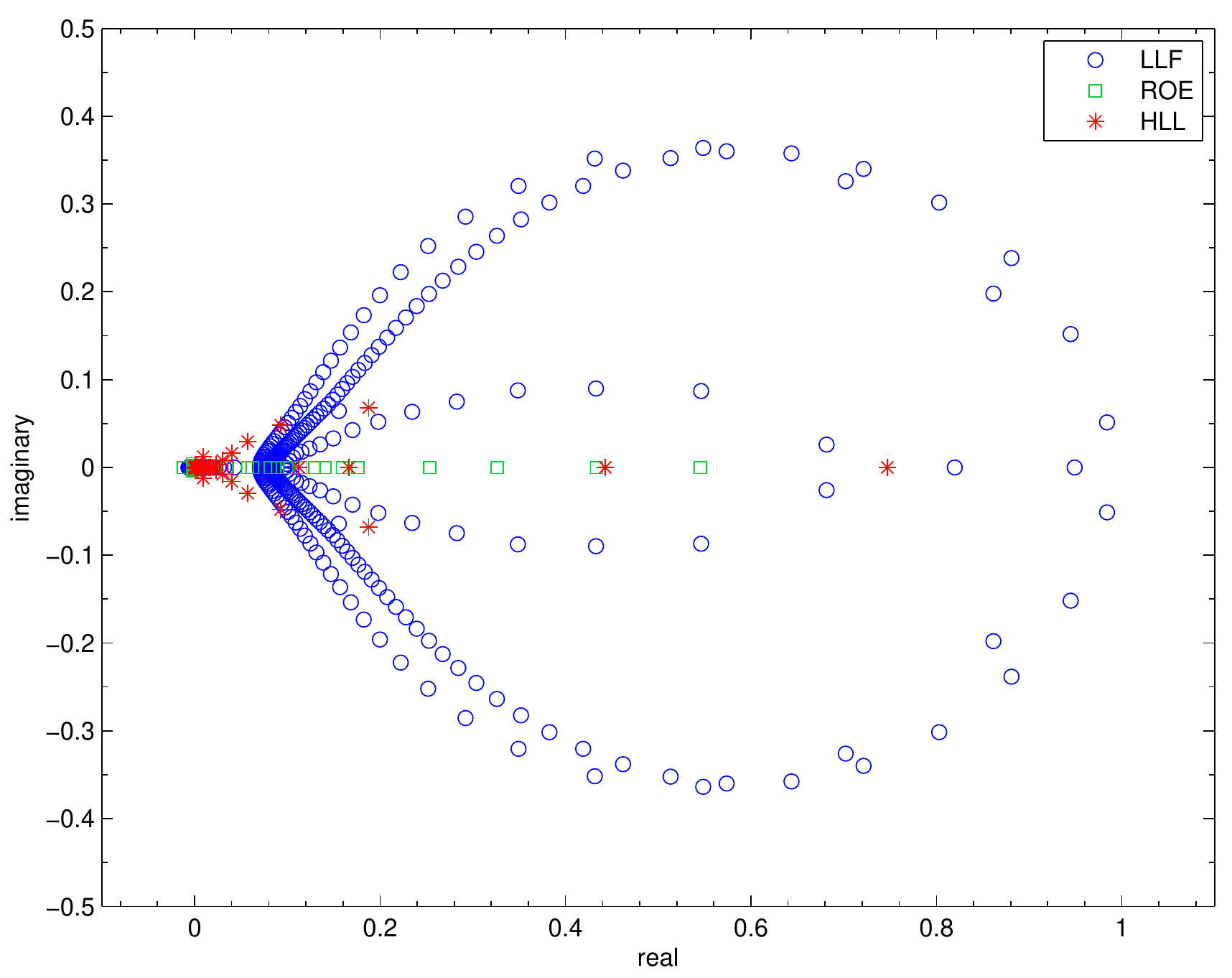} }
\subfigure[HLL flux]
      { \includegraphics[width=0.46\textwidth,height=5cm]{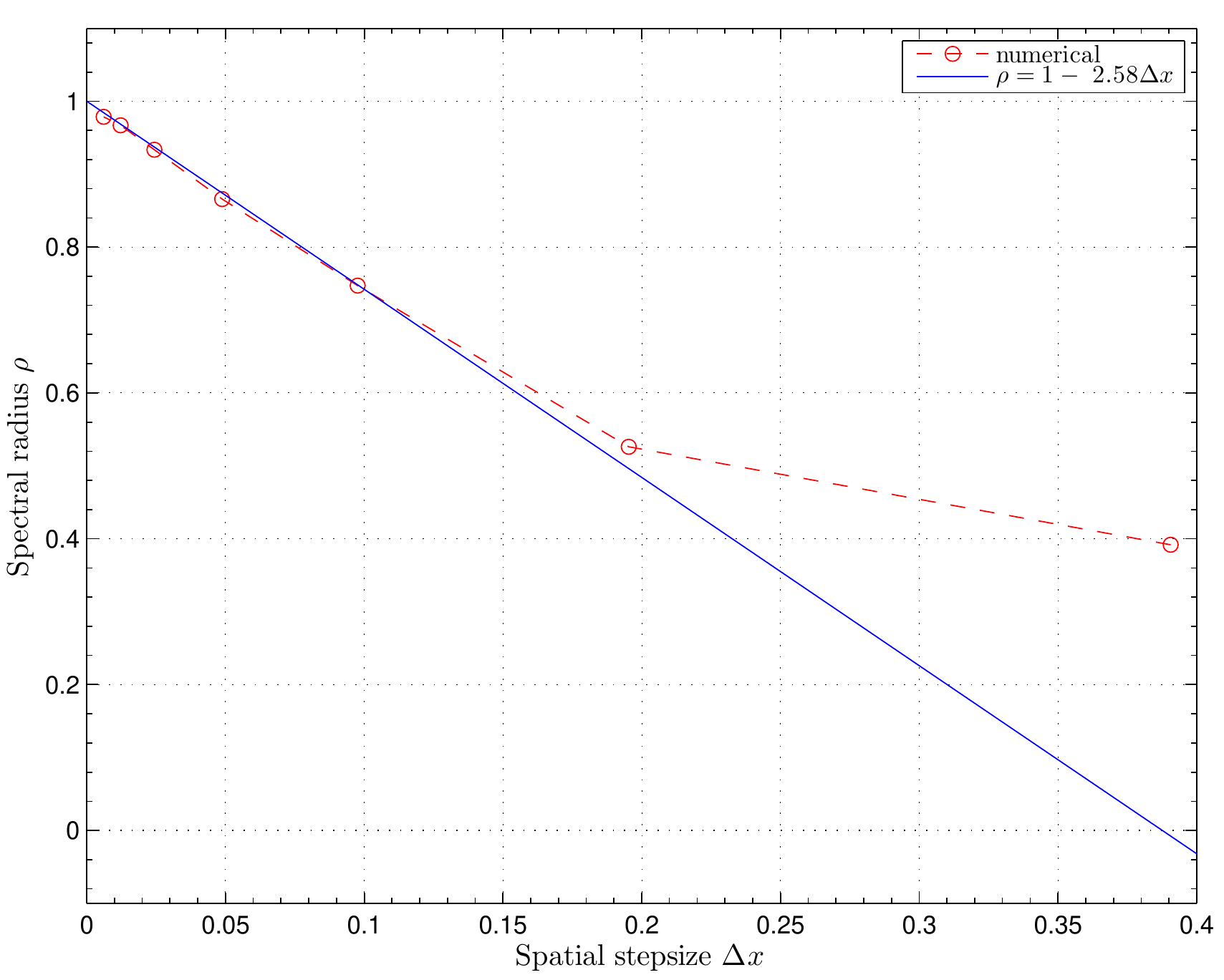} }
\subfigure[LLF flux]
      { \includegraphics[width=0.46\textwidth,height=5cm]{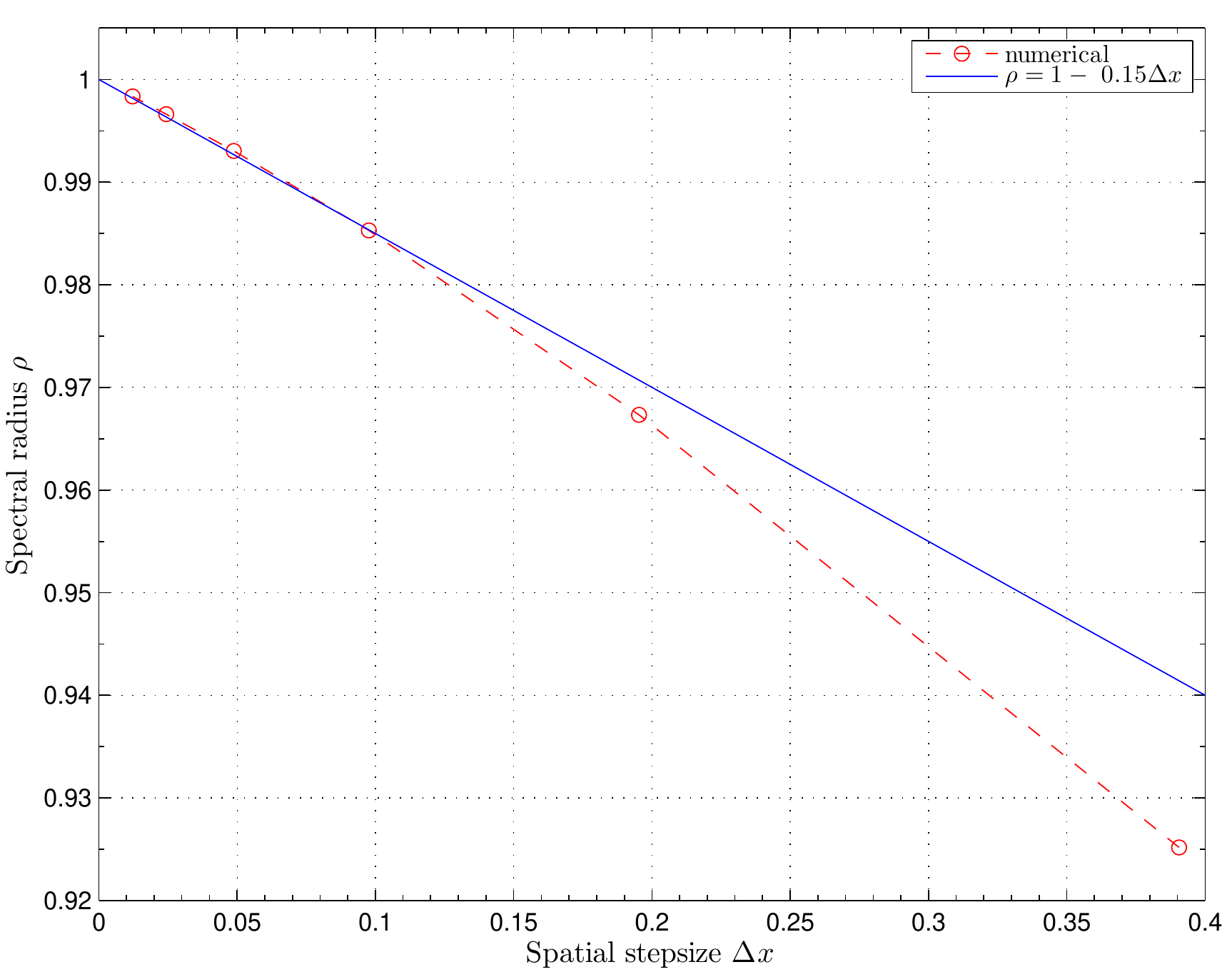} }
\subfigure[Roe flux]
      { \includegraphics[width=0.46\textwidth,height=5cm]{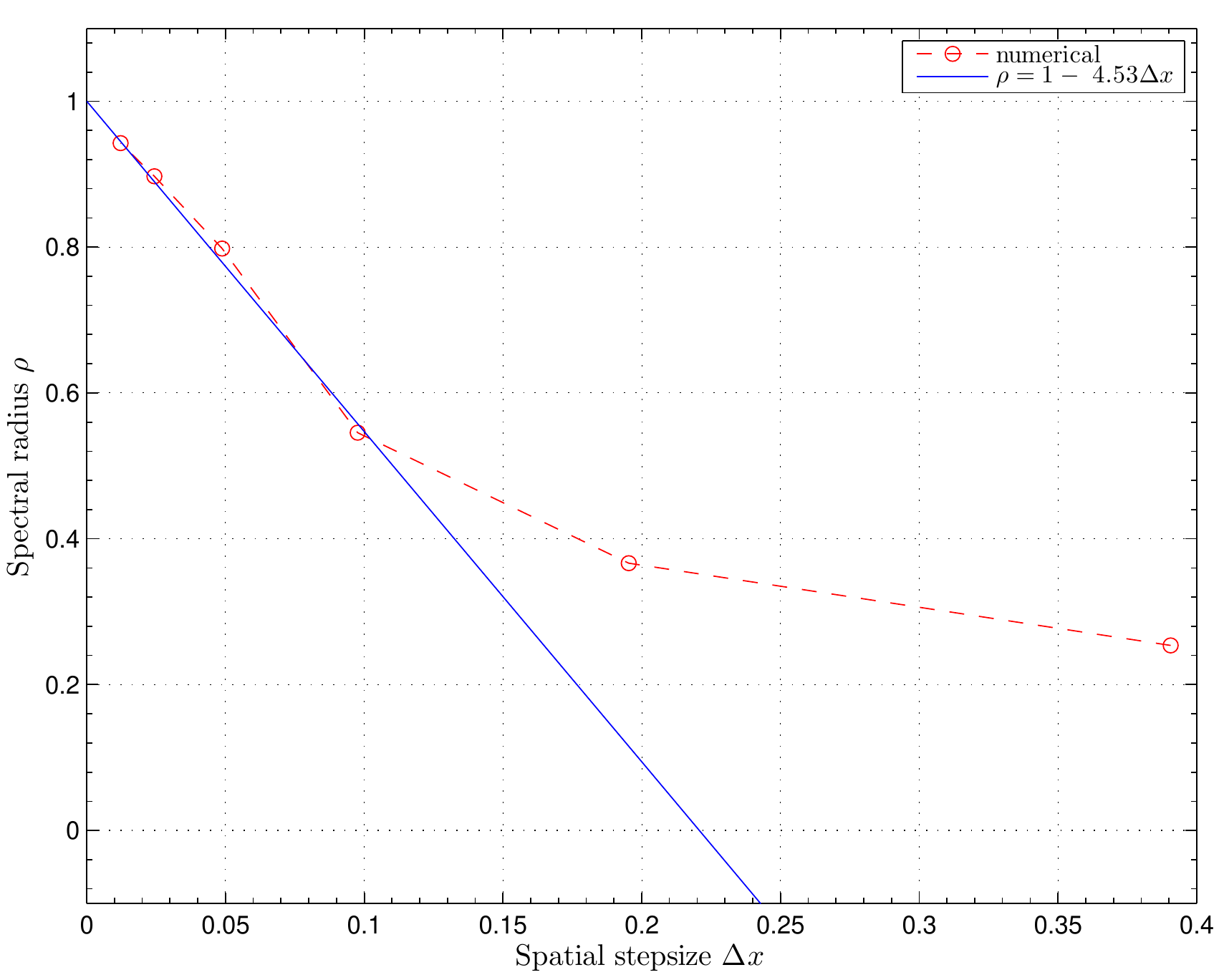} }
  \caption{\small Example \ref{example1Dtrans}(II): (a). Distribution of the eigenvalues
 of the block SGS iteration matrix
with the HLL, LLF, and Roe  fluxes on the mesh of 256 uniform cells.
(b-d). Asymptotic relation of the spectral radius $\rho$ of the block SGS iteration matrix
with respect to the spatial step size $\Delta x$.
}\label{fig:1d_01b_eig}
\end{figure}

%\begin{figure}[!htbp]
%  \centering
%\subfigure[1st-order accurate scheme with 256 uniform cells]
%      { \includegraphics[width=0.46\textwidth,height=5cm]{fig1D/1D01b/EigEx1D01b} }
%\subfigure[1st-order accurate scheme with 512 uniform cells]
%      { \includegraphics[width=0.46\textwidth,height=5cm]{fig1D/1D01b/EigEx1D01b-512} }
%\subfigure[2nd-order accurate scheme with 256 uniform cells]
%      { \includegraphics[width=0.46\textwidth,height=5cm]{fig1D/1D01b/2ndEigEx1D01b-256} }
%\subfigure[2nd-order accurate scheme with 512 uniform cells]
%      { \includegraphics[width=0.46\textwidth,height=5cm]{fig1D/1D01b/2ndEigEx1D01b-512} }
%  \caption{\small Example \ref{example1Dtrans}(II): Distribution of the eigenvalues
%(on the complex plane) of the iteration matrix for the block SGS method
%with different numerical fluxes.
%}\label{fig:1d_01b_eig2}
%\end{figure}

%Similar to Example \ref{example1Dsmooth}, using different numerical fluxes gives
%different convergence behaviors of  {\tt NMGM}.
The effect of  the HLL, LLF,  and Roe fluxes
on the convergence behavior  of  {\tt NMGM} and  the BLU-SGS iteration.
is investigated in Table \ref{tab:1Dtrans:b}, which shows that
the {\tt NMGM} performs  much more efficiently than the BLU-SGS iteration,
and it becomes obvious as $N_L$ increases in the multigrid
and the W-cycle multigrid is adopted with the LLF flux. Moreover,
convergence behaviors depend on the numerical flux, and
the block LU-SGS method with the HLL or Roe  flux is more efficient
than the LLF flux. Such observation is consistent with the previous.
However, different from the results of Example \ref{example1Dsmooth} and \ref{example1Dtrans}(I),
the results in Table \ref{tab:1Dtrans:b} show that
$N_{\rm step}$ increases almost linearly with the mesh refinement for the BLU-SGS iteration with the HLL, LLF,  and Roe fluxes, and the spectral radius $\rho$ of the iteration matrix increases asymptotically to $1$ as $N$ increases.
The distribution in the complex plane of the eigenvalues of the block SGS iteration
matrix  in Fig. \ref{fig:1d_01b_eig}(a) still shows that
the ``high frequency'' eigenvalues  for the LLF flux  are much more
than  for the HLL or Roe  flux.
Figs. \ref{fig:1d_01b_eig}(b-d) show
the asymptotic relation  \eqref{eq-04-01} of $\rho$ with respect to the spatial step size $\Delta x$,
%$$\rho \sim  1 - C \Delta x , \quad {\rm as}~ \Delta x \to 0 $$
where  $C$ is about  $2.58, 0.15$, and $4.53$ for the HLL, LLF, and Roe fluxes, respectively.

\end{example}

%%%%%%%%%%%%%%%%%%%%%%%%%%%%%%%%%%%%%% Example 1D 04 %%%%%%%%%%%%%%%%%%%%%%%%%%%%%%%

\begin{example}[Wet-dry boundary problem] \label{example1Ddry}\rm

The bed topography  for this problem is the same  as one in Example \ref{example1Dtrans},
%i.e.
%\begin{equation*}
%z(x)=
%\begin{cases}
%0.2 - 0.05(x-10)^2, &  8 < x < 12,\\
%0,&    \rm{otherwise},\end{cases}
%\end{equation*}
but for a channel of length 20.
Initially, the flow with the water height $h^{(0)}(x)=0.22-z(x)$ is static in the channel, i.e., $u(x)=0$.
A steady state with the dry bed
\begin{align*}
\left\{
\begin{aligned}
& h(x) = \max \{0.1 - z(x),0\}, \\
& u(x) = 0,
\end{aligned}
\right.
\end{align*}
will  be reached if  imposing $h=0.1$  at   the interval ends  $x=0$ and
$20$. %setting two ``dam breaking'' boundaries
% The steady solution called ``lake at rest with an emerged bump''
%was used to verify the well-balanced property of a numerical scheme for shallow water equations \cite{Delestre2013}.
Since the steady solution involves wet/dry transition,
the Jacobian matrix of these numerical fluxes near the
wet-dry boundary becomes singular so that
it is much more difficult and challenging.

 Fig. \ref{fig:1d_04_solu} shows the
 steady solutions $h+z$ and the error in the rest  water surface
 obtained on the mesh of 512 uniform cells, where the LLF flux and  W-cycle multigrid with $N_L=3$ are used.
 It is seen that
 the correct steady solutions are efficiently obtained by {\tt NMGM}, and
 the rest water surface is preserved exactly up to the machine precision.
 %which verifies the ``well-balanced property'' of our scheme.
Fig. \ref{fig:1d_04_Re} gives the convergence history of {\tt NMGM} in terms of  the  {\tt NMGM} iteration number  $N_{\rm step}$ and CPU time on four
meshes of 512, 1024, 2048, and 4096 uniform cells, respectively.
Convergence behaviors of {\tt NMGM} and $N_{\rm step}$ do not depend
on  the uniform  mesh number $N$.

  Because  {\tt NMGM}  and the BLU-SGS iteration with the HLL or Roe flux {fail} to work now,
Table \ref{tab:1Ddry} only {investigates} the convergence behaviors of {\tt NMGM} and the BLU-SGS iteration
with the LLF flux.
We see that {\tt NMGM} is much more efficient
than the BLU-SGS iteration and it is obvious  as $N_L$ increases
and the W-cycle  multigrid is adopted.
In this problem, the big solution error is mainly  introduced around  the wet-dry boundary and can be
fast reduced by using  the W-cycle
multigrid with properly increasing the coarse level number $N_L$.

\begin{figure}[htbp]
	\centering
	\includegraphics[width=0.46\textwidth,height=5cm]{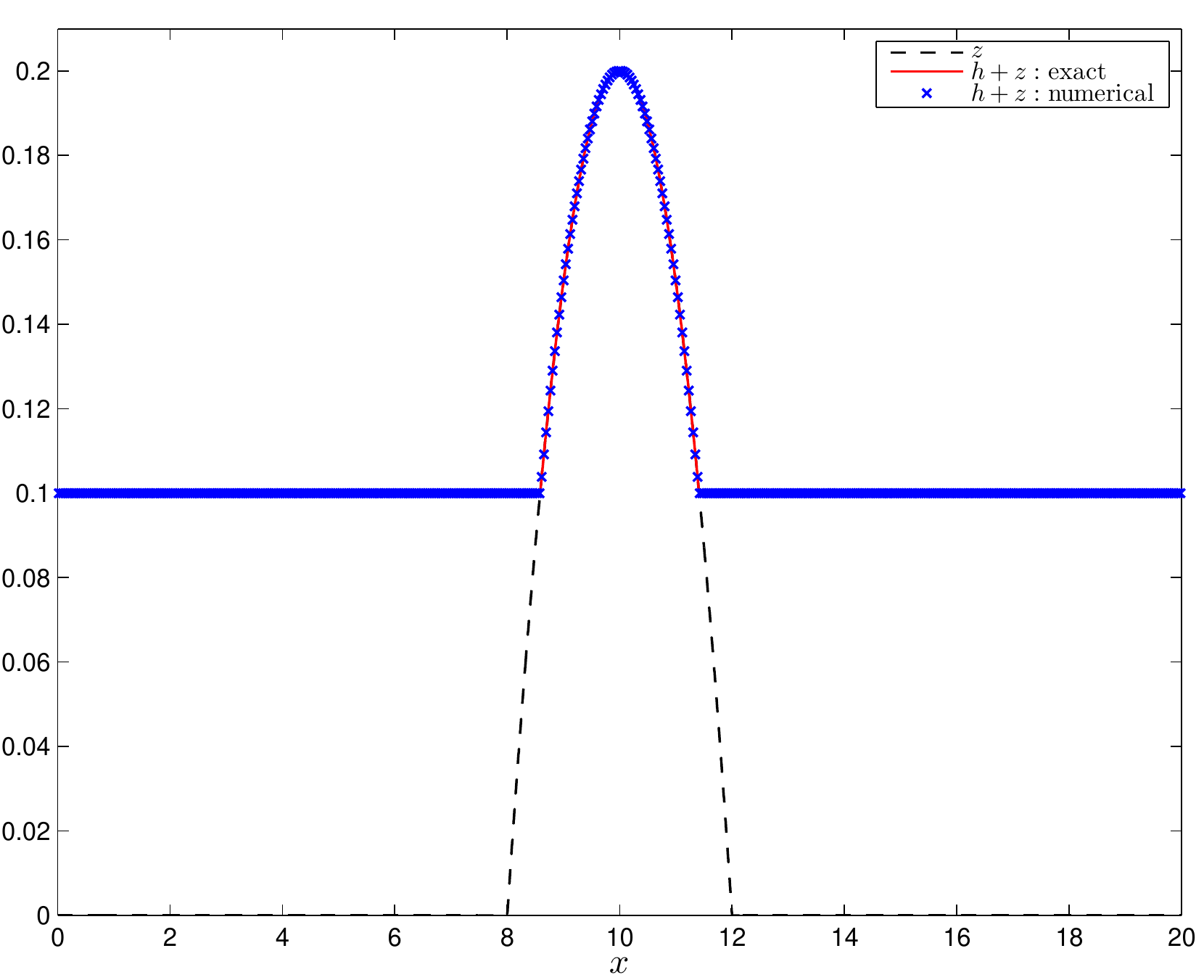}
	\includegraphics[width=0.46\textwidth,height=5cm]{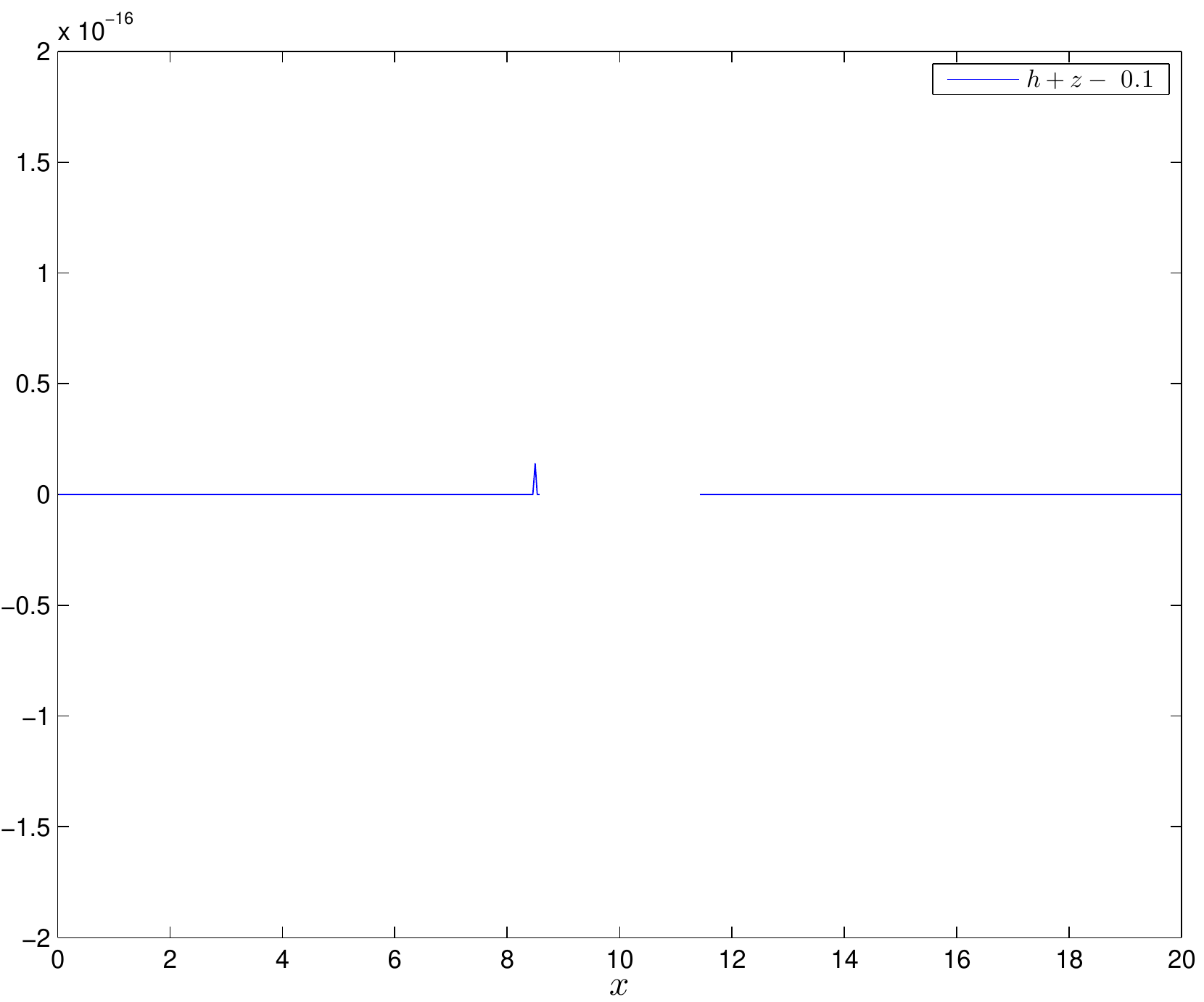}
	\caption{\small Example \ref{example1Ddry}:
		Steady-state solution (left) and the error in the water surface (right) obtained by {\tt NMGM} on the mesh of 512 uniform cells.
	}\label{fig:1d_04_solu}
\end{figure}

\begin{figure}[htbp]
	\centering
	\includegraphics[width=0.46\textwidth,height=5cm]{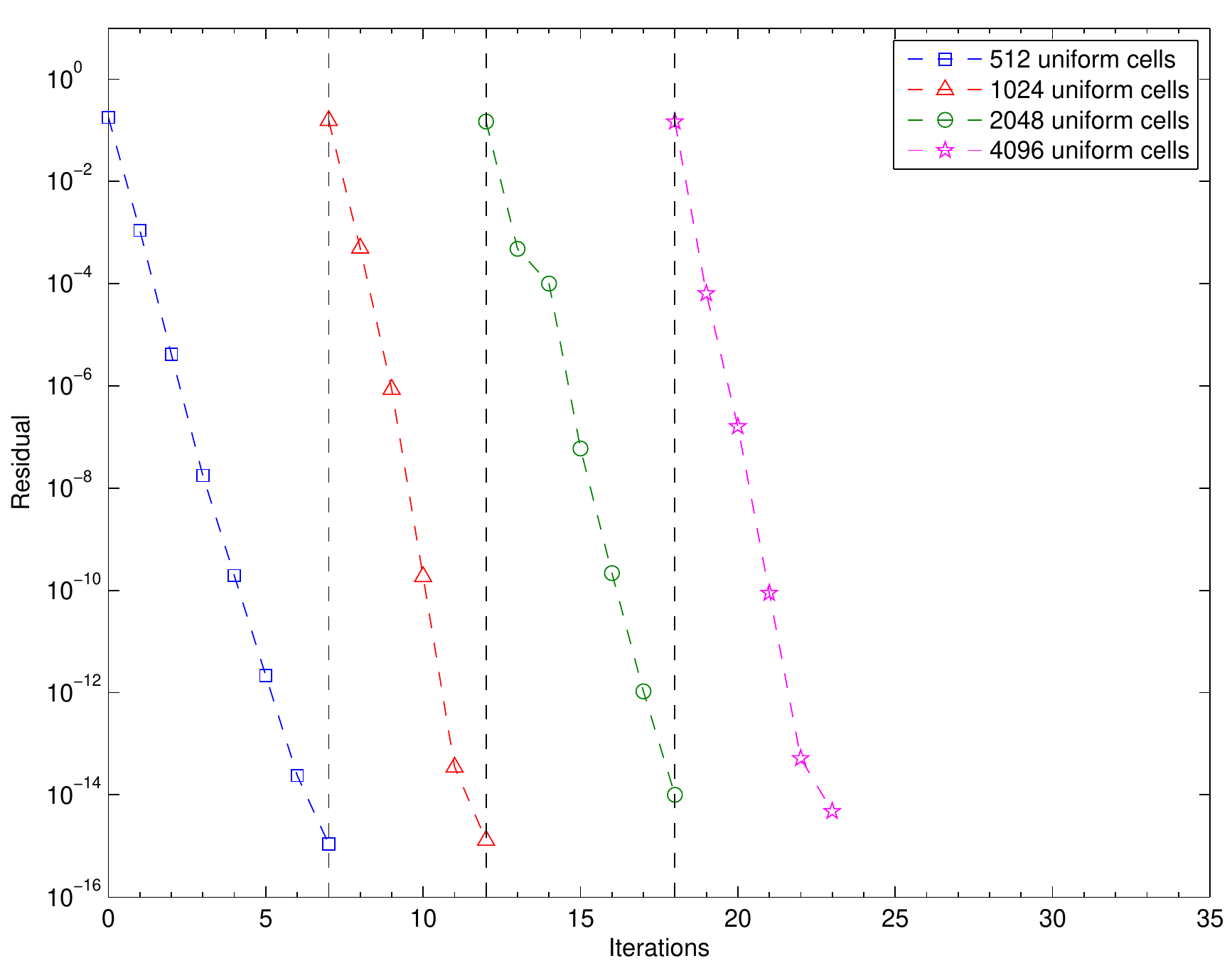}
	\includegraphics[width=0.46\textwidth,height=5cm]{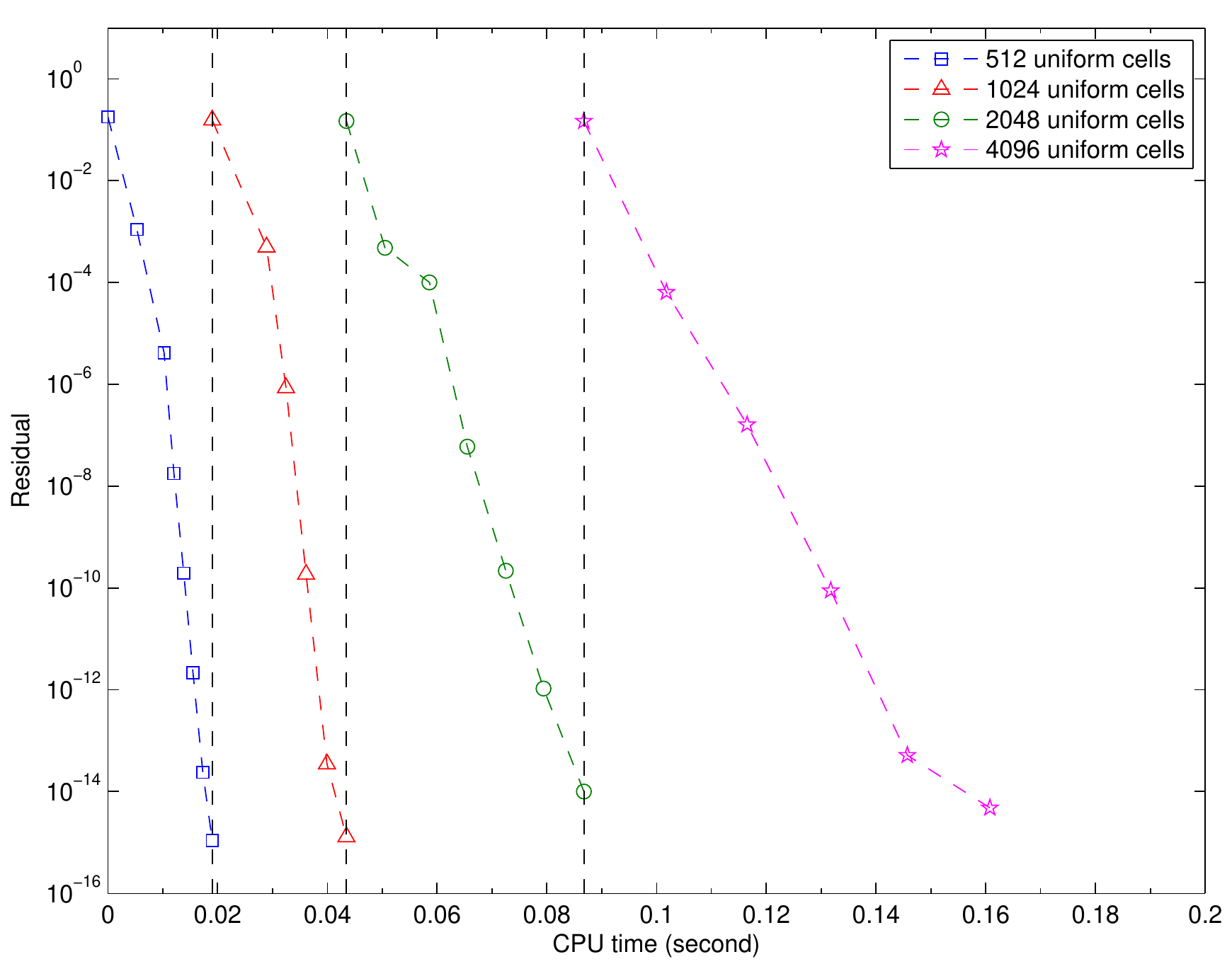}
	\caption{\small Example \ref{example1Ddry}:
		Convergence history in terms of  the  {\tt NMGM} iteration number  $N_{\rm step}$  (left) and CPU time (right) on four uniform meshes.
	}\label{fig:1d_04_Re}
\end{figure}

\begin{table}[htbp] %\multirow{3}{2pt}{$N$}
\centering
\caption{\small Example \ref{example1Ddry}: Convergence behaviors of
{\tt NMGM} and the BLU-SGS iteration with different $N_L$,  V or
W cycle  multigrid, and   LLF flux.
  }
\scriptsize
\begin{tabular}{|c|c|c|c|c|c|c|c|c|}
\hline
\multicolumn{2}{|c|}{$N$} & 64 & 128 & 256 & 512 & 1024 & 2048 & 4096 \\
\hline
 Block & { $N_{\rm step}$}  & 42 & 69  & 129 & 364 & 1122 & 4009 & 14459 \\
 LU-SGS & { $T_{\rm cpu}$}   & 8.16e-3 & 1.78e-2  & 3.41e-2 & 1.37e-1 & 7.71e-1 & 5.17e0 & 3.63e1 \\
\cline{1-9}
V-cycle & { $N_{\rm step}$}  & 4 & 5  & 6 & 11 & 32 & 110  & 356\\
{$N_L = 1$} & { $T_{\rm cpu}$}  & 1.73e-3 & 4.12e-3  & 1.03e-2 & 2.30e-2 & 8.13e-2 & 4.91e-1 & 3.08e0 \\
\cline{1-9}
V-cycle & { $N_{\rm step}$}  & 4 & 5  & 5 & 7 & 9 & 11  & 44\\
{$N_L = 3$} & { $T_{\rm cpu}$}   & 2.02e-3 & 4.58e-3  & 9.82e-3 & 1.63e-2 & 2.55e-2 & 5.92e-2  & 4.23e-1\\
\cline{1-9}
V-cycle & { $N_{\rm step}$}  & 4 & 5  & 5 & 7 & 9 & 12  & 12\\
{$N_L = 5$} & { $T_{\rm cpu}$}  & 2.07e-3 & 4.76e-3  & 1.01e-2 & 1.76e-2 & 2.64e-2 & 6.06e-2  & 1.18e-1\\
\cline{1-9}
W-cycle & { $N_{\rm step}$}  & 4 & 5  & 6 & 8 & 18 & 58  & 198\\
{$N_L = 1$} & { $T_{\rm cpu}$}  & 1.91e-3 & 4.62e-3  & 1.16e-2 & 1.72e-2 & 4.82e-2 & 2.85e-1  & 1.94e0\\
\cline{1-9}
W-cycle & { $N_{\rm step}$}  & 4 & 5  & 5 & 7 & 5  & 6  & 5\\
{$N_L = 3$} & { $T_{\rm cpu}$}  & 2.58e-3 & 5.78e-3  & 1.90e-2 & 2.45e-2 & 2.93e-2 & 4.33e-2  & 7.40e-2\\
\cline{1-9}
\multicolumn{2}{|c|}{$\rho$}   & 0.45675 & 0.63258  & 0.80392
 & 0.92105 & 0.97479 & 0.99296  & 0.99816\\
\cline{1-9}
\multicolumn{2}{|c|}{$R_{\infty}$}   & 7.836e-1 & 4.579e-1  & 2.183e-1
 & 8.224e-2 & 2.553e-2 & 7.062e-3 & 1.846e-3 \\
\hline
\end{tabular}\label{tab:1Ddry}
\end{table}

\end{example}

%%%%%%%%%%%%%%%%%%%%%%%%%%%%%%%%%%%%%%%%%%%%%%% 2D Examples %%%%%%%%%%%%%%%%%%%%%%%%%%%%%%%%%%%%%%%

\subsection{2D case}
Unless specifically stated, in the following,
the HLLC flux and  V-cycle multigrid are adopted, and the coarse  mesh level number
$N_L$ is set to be $3$.

\begin{figure}[htbp]
  \centering
  \includegraphics[width=0.47\textwidth]{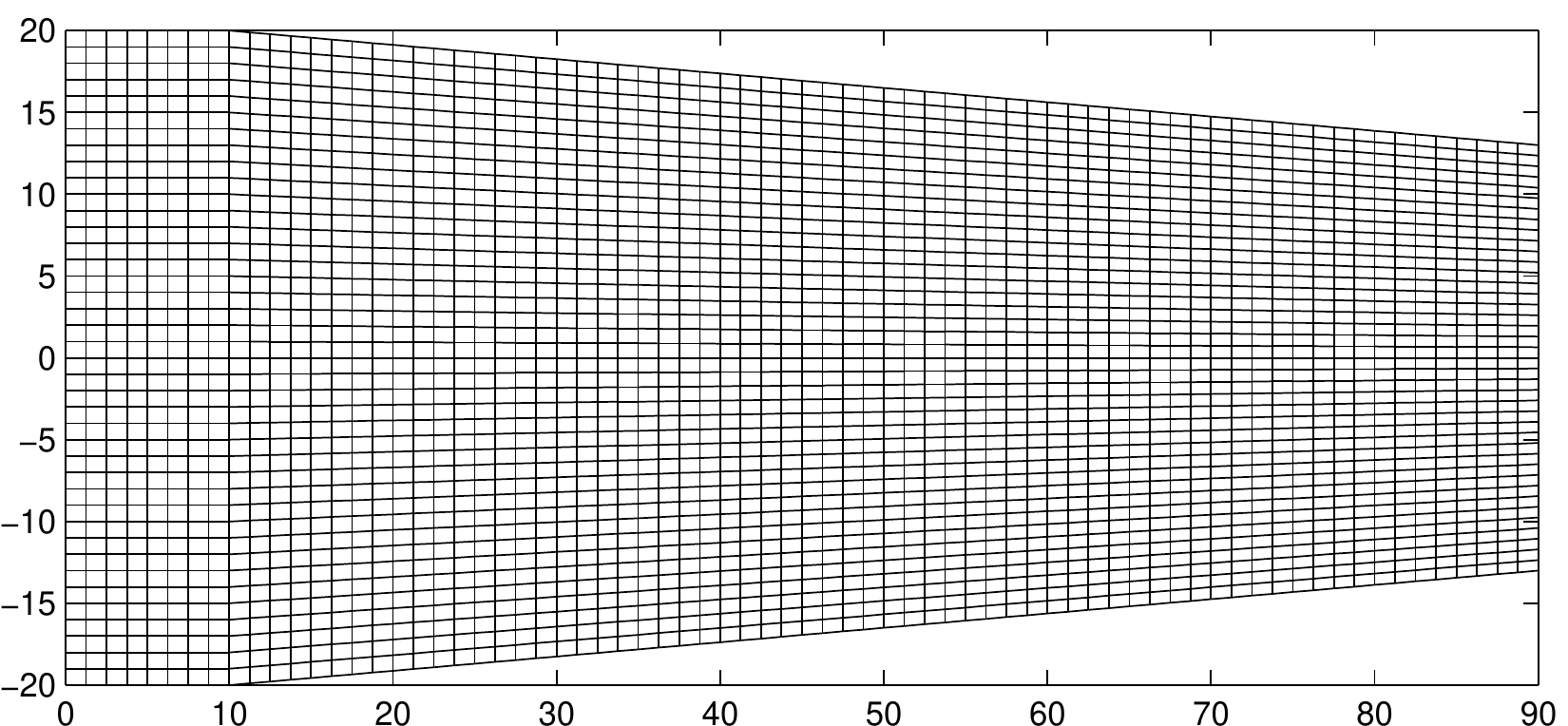}
  \includegraphics[width=0.47\textwidth]{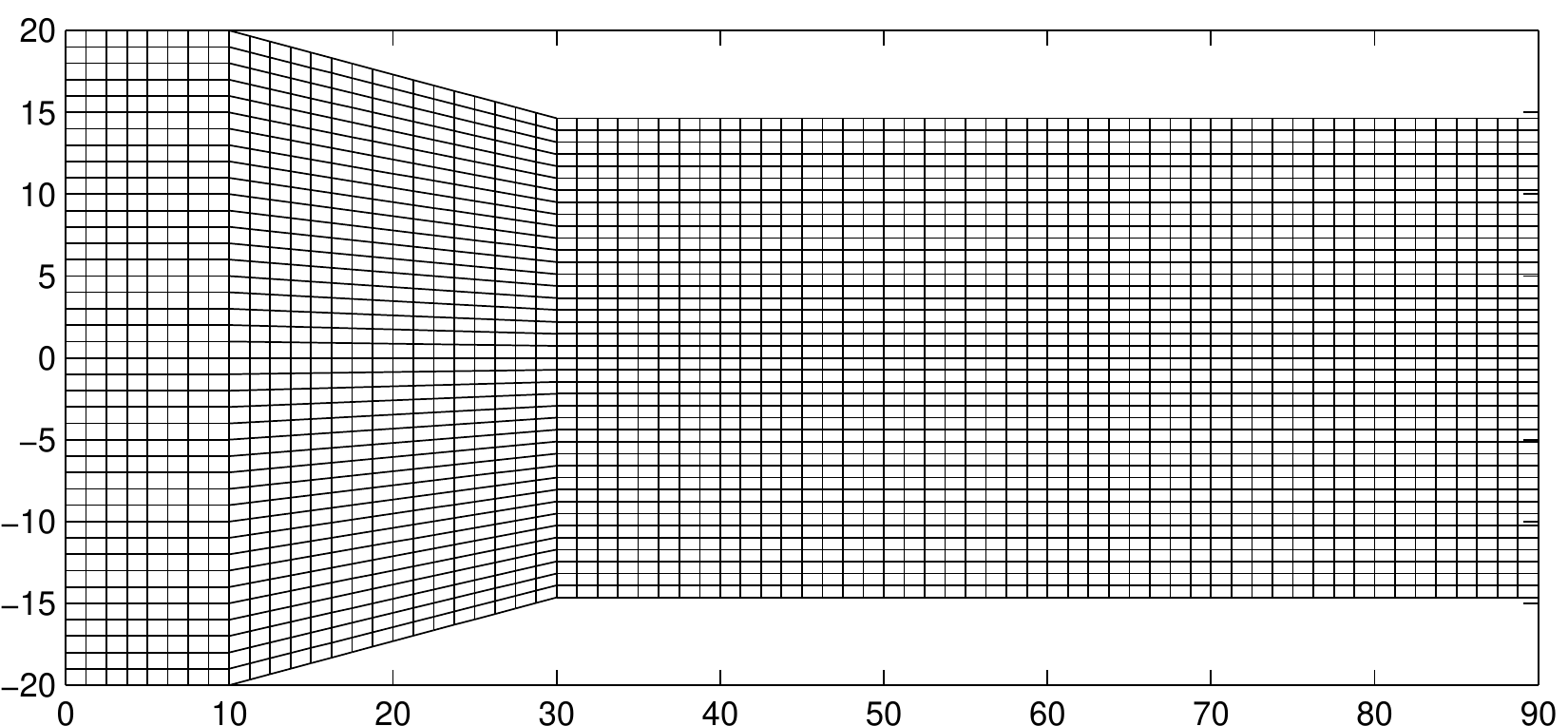}
  \caption{\small Example \ref{example2Dsym}:
 The  geometries of {\em channel I} (left) and {\em II} (right)
and  structured mesh of $72 \times 40$ cells.
}\label{fig:2d_01_mesh}
\end{figure}

\begin{example} \label{example2Dsym}\rm
The first 2D example is to simulate  fluid flows
in two symmetric channels with flat
bottom constricted from both side in the $y$--direction, see
\cite{Tang2004}. %\cite{Liska1999,Tang2004}.
Specifically,
 {\em channel I} is with a constriction angle
$\alpha = 5^{\circ}$ started at $x=10$, and {\em channel II} with a constriction
angle $\alpha = 15^{\circ}$ started at $x=10$, ended at $x=30$,  and then
followed  by a straight narrower channel,
see Fig. \ref{fig:2d_01_mesh} for their geometries and corresponding  structured mesh of $72 \times 40$ cells.
The inflow conditions  $h=1$ and $v=0$ and the Froude number
$F_r = |\vec u|/\sqrt{gh}=2.5$ are imposed at $x=0$,  the outflow boundary
condition is specified at $x=90$,  while slip boundary
conditions are employed on the top and bottom boundaries.

\begin{figure}[htbp]
  \centering
  \includegraphics[width=0.47\textwidth]{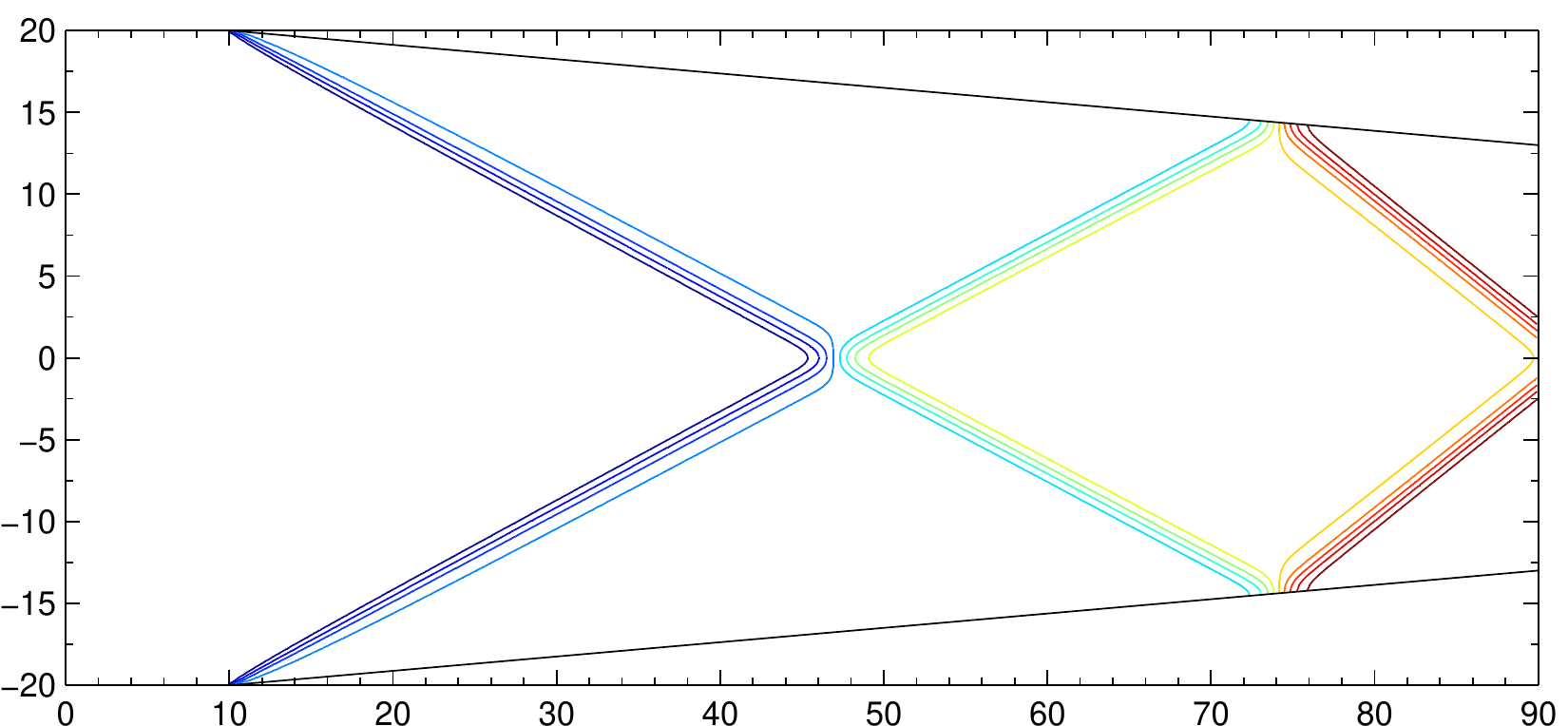}
  \includegraphics[width=0.47\textwidth]{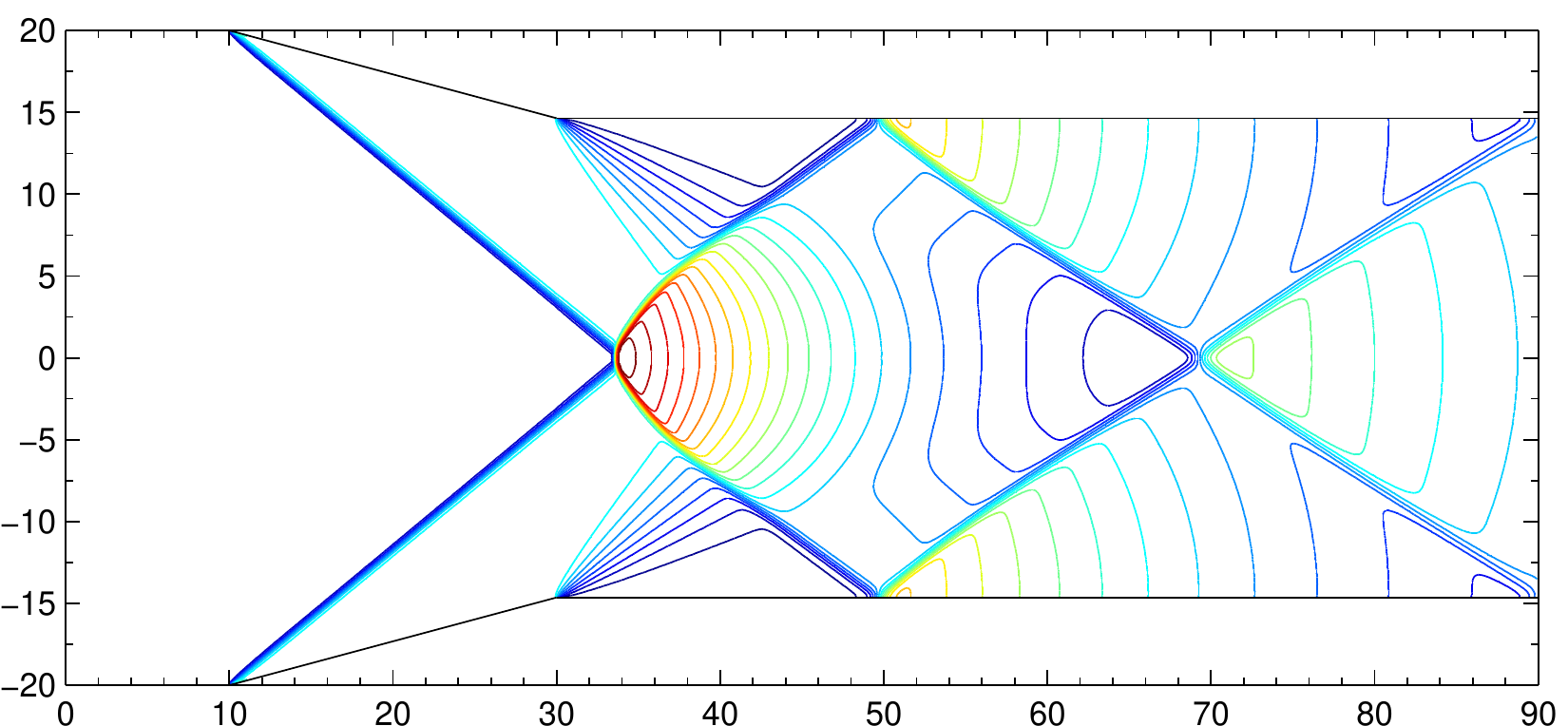}
  \caption{\small Example \ref{example2Dsym}:
The steady water depth $h$ with 13  and 20  equally spaced
contour lines for {\em channel I}  and {\em   II}  obtained by
{\tt NMGM} on the mesh of $576\times 320$ cells,  respectively.
}\label{fig:2d_01_solu}
\end{figure}

Fig. \ref{fig:2d_01_solu} displays
the contour plots of the steady solutions for  two channels
obtained by {\tt NMGM} on the mesh of $576\times 320$ cells.
For the  steady state flow in {\em channel I},  two  bore waves starting at the points
$(10,\pm 20)$  interact at the point (45,0) and then two regular
{reflections happen} around $x=74$ due to the constriction. The present result is well  comparable to the numerical result
 given by using the adaptive moving mesh method to solve the time-dependent SWEs \cite{Tang2004} and
analytical ones \cite{Zienkiewicz1995}, especially,
numerical values of the water heights
of the first and second plateau are $1.25$ and $1.5271$, respectively, which agree well with those %numerical results $h_2 = 1.25$, $h_3 = 1.527$ given
in \cite{Tang2004,Zienkiewicz1995}.
%and analytical results $h_2 = 1.254$, $h_3 = 1.55$
%\cite{Zienkiewicz1995}.
For the steady flows in {\em channel II}, two shock waves
started at $(10,\pm20)$ meet  around (33,0) each other,
then interact with two expansion waves generated at the points $(30,\pm15)$, and then two regular reflections {happen} at
$(50,\pm15)$. Finally, two reflected shock waves meet around (69,0).
From the contour plots in Fig. \ref{fig:2d_01_solu}, we see
that {\tt NMGM} can  well capture those waves and their interactions with high resolution and non-oscillation.

Figs. \ref{fig:2d_01A_Re} and \ref{fig:2d_01B_Re} show
the convergence history of {\tt NMGM} versus the  {\tt NMGM} iteration number $N_{\rm step}$ and CPU time
on three successively refined meshes for {\em channel I} and {\em channel II}, respectively.
It is seen that {\tt NMGM} works successfully on those meshes
and the convergence behaviors of {\tt NMGM} are independent
on the cell number $N$.
Compared to {\em channel I}, the computation of the steady solution for {\em channel II} seems
more difficult for {\tt NMGM} and takes more iteration steps.
In  those computations,  $\epsilon_p=2\times 10^2$.

%the solutions are non-oscillatory and the symmetry is preserved.

\begin{figure}[htbp]
  \centering
  \includegraphics[width=0.46\textwidth,height=5cm]{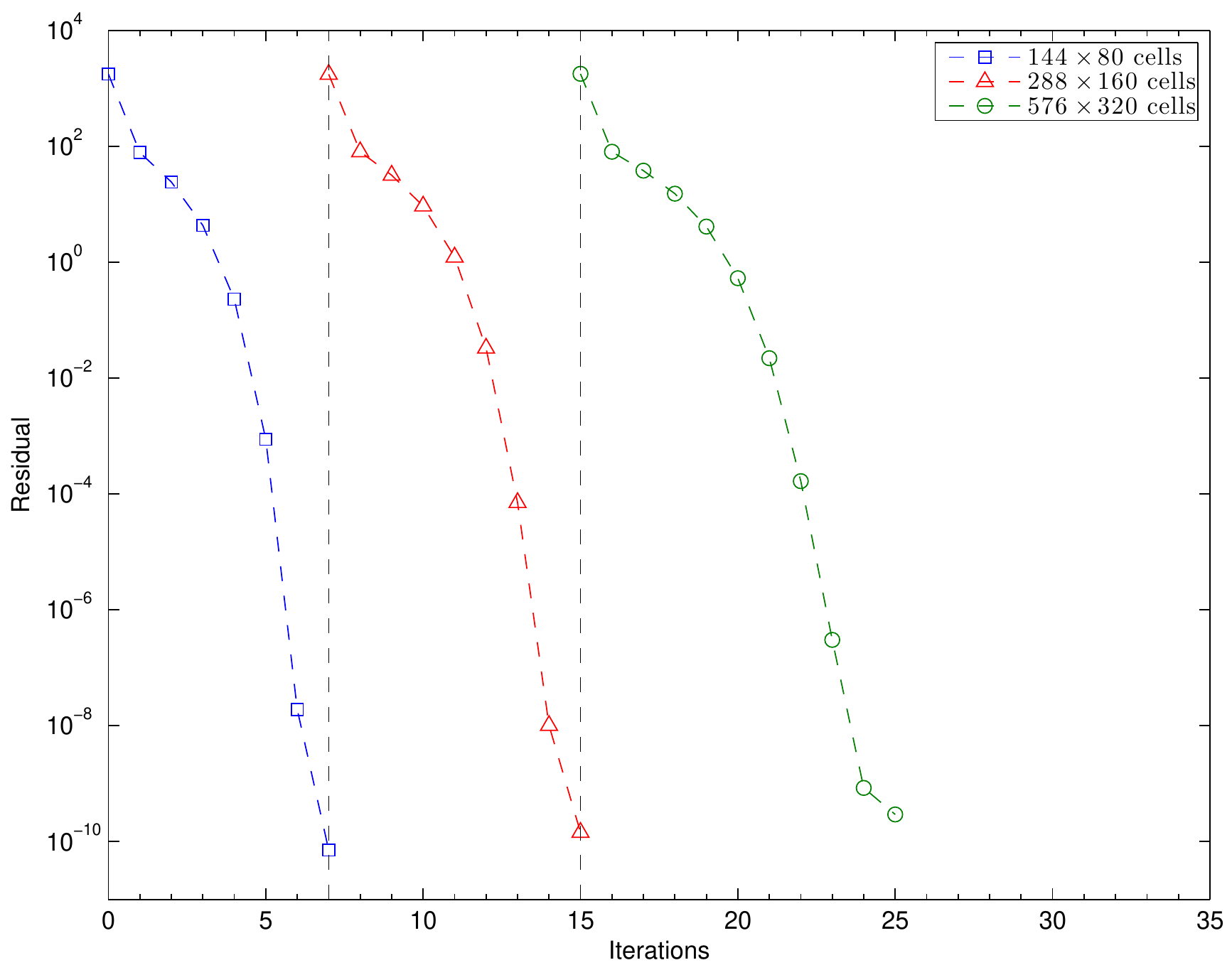}
  \includegraphics[width=0.46\textwidth,height=5cm]{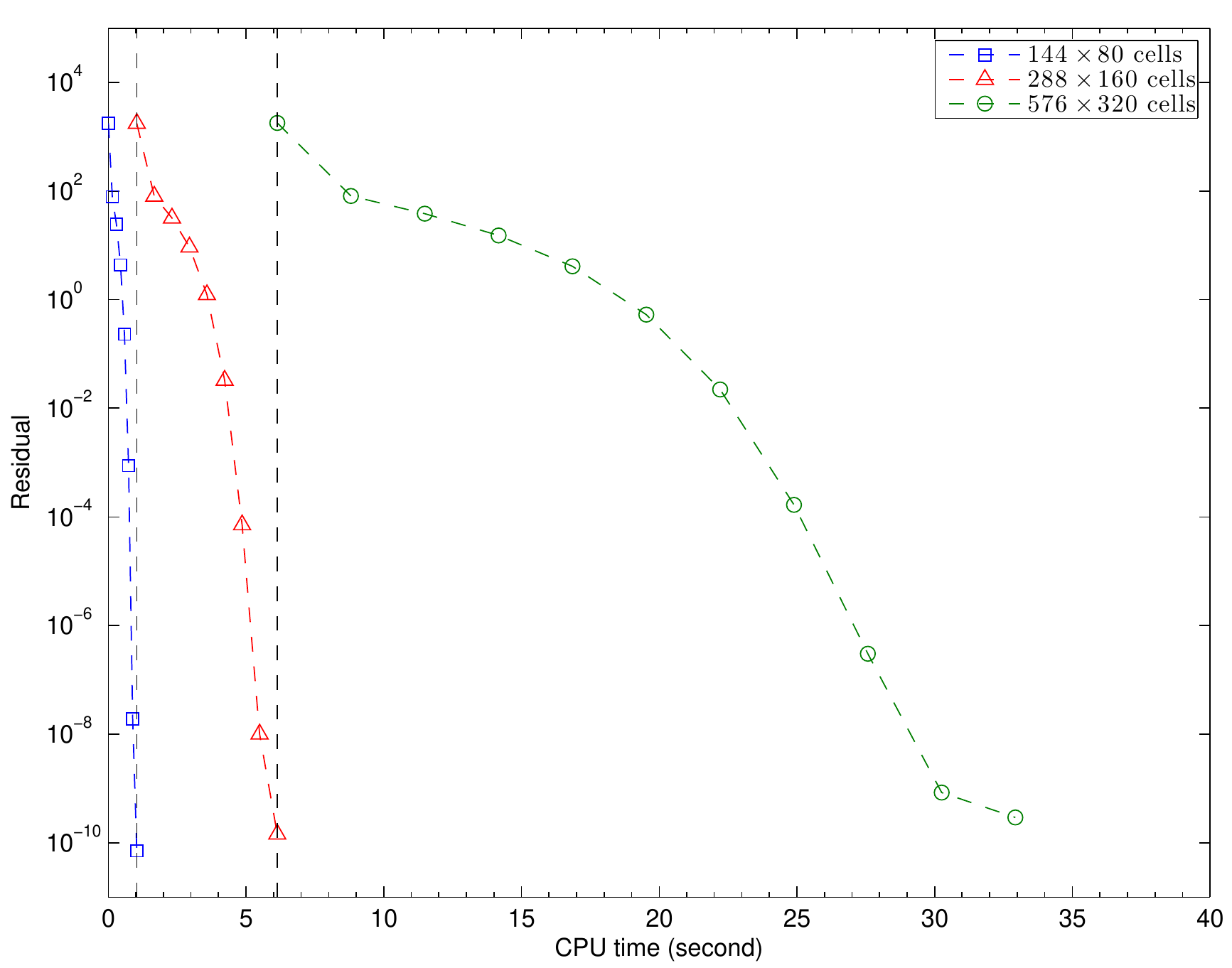}
  \caption{\small Example \ref{example2Dsym}:
Convergence history in terms of  the  {\tt NMGM} iteration number  $N_{\rm step}$  (left) and CPU time (right) on three meshes
 for {\em channel I}.
}\label{fig:2d_01A_Re}
\end{figure}

\begin{figure}[htbp]
  \centering
  \includegraphics[width=0.46\textwidth,height=5cm]{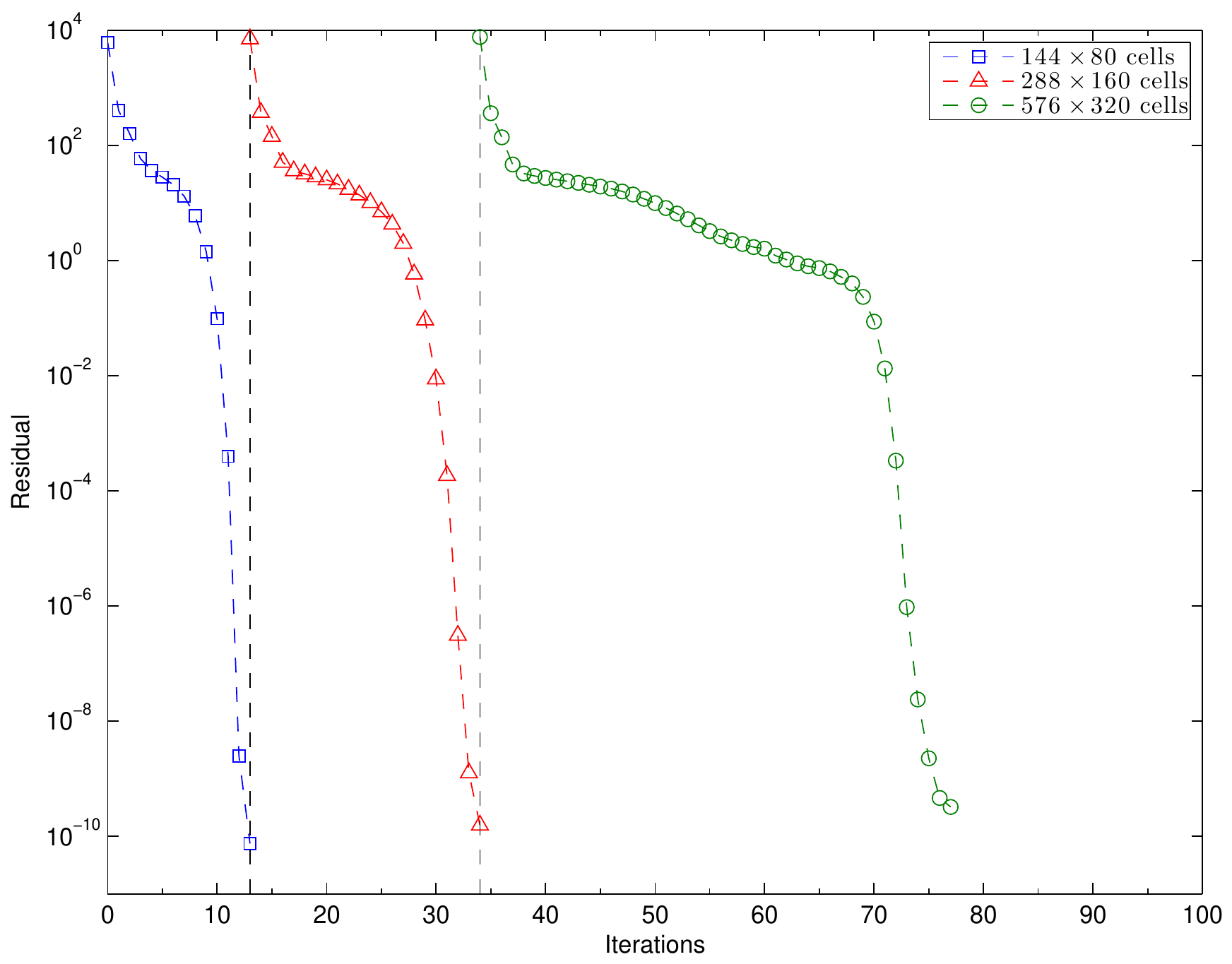}
  \includegraphics[width=0.46\textwidth,height=5cm]{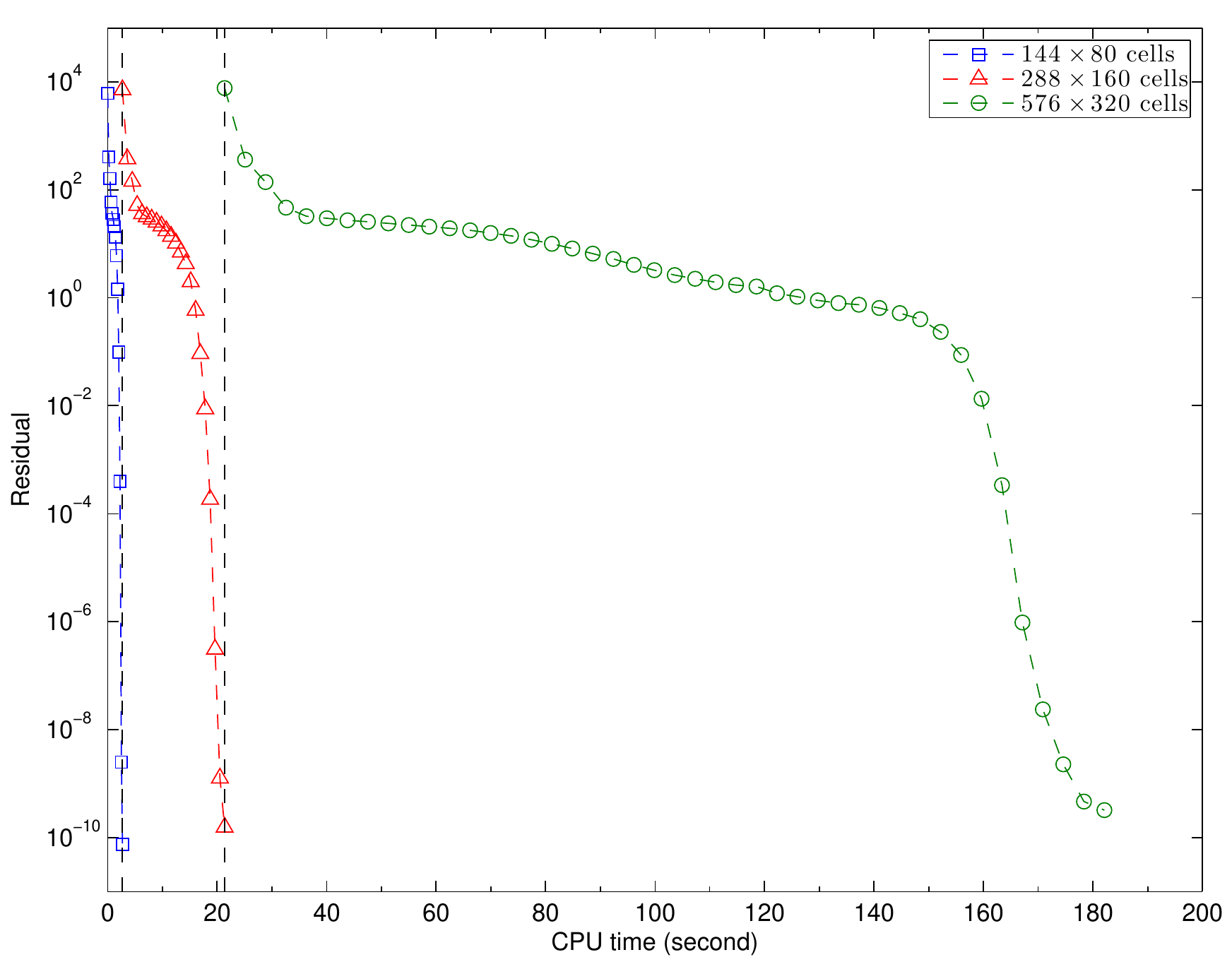}
  \caption{\small %Example \ref{example2Dsym}:
Same as Fig.\ref{fig:2d_01A_Re}, except for {\em channel II}.
}\label{fig:2d_01B_Re}
\end{figure}

\end{example}

%\begin{figure}[htbp]
%  \centering
%  \includegraphics[width=0.65\textwidth]{fig2D/2D02/mesh96x32}
%  \caption{\small Example \ref{example2Dsym2}:
% The detailed geometries and the  type of structured mesh (with $96 \times 32$ cells shown here)
%used in our computations.
%}
%\label{fig:2d_02_mesh}
%\end{figure}

\begin{figure}[htbp]
  \centering
  \subfigure[Structured mesh of $96 \times 32$]
        {
    \includegraphics[width=0.47\textwidth]{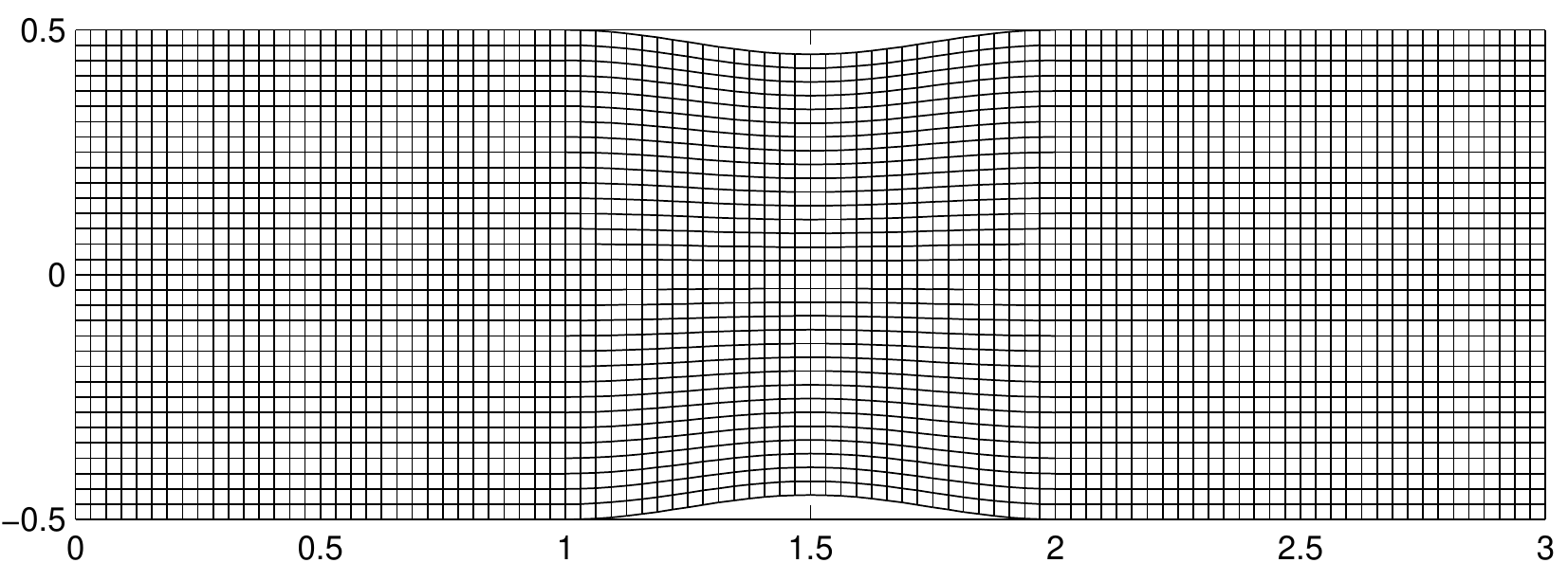}
    }
  \subfigure[$F_{\rm in} = 0.5$]
        {
  \includegraphics[width=0.47\textwidth]{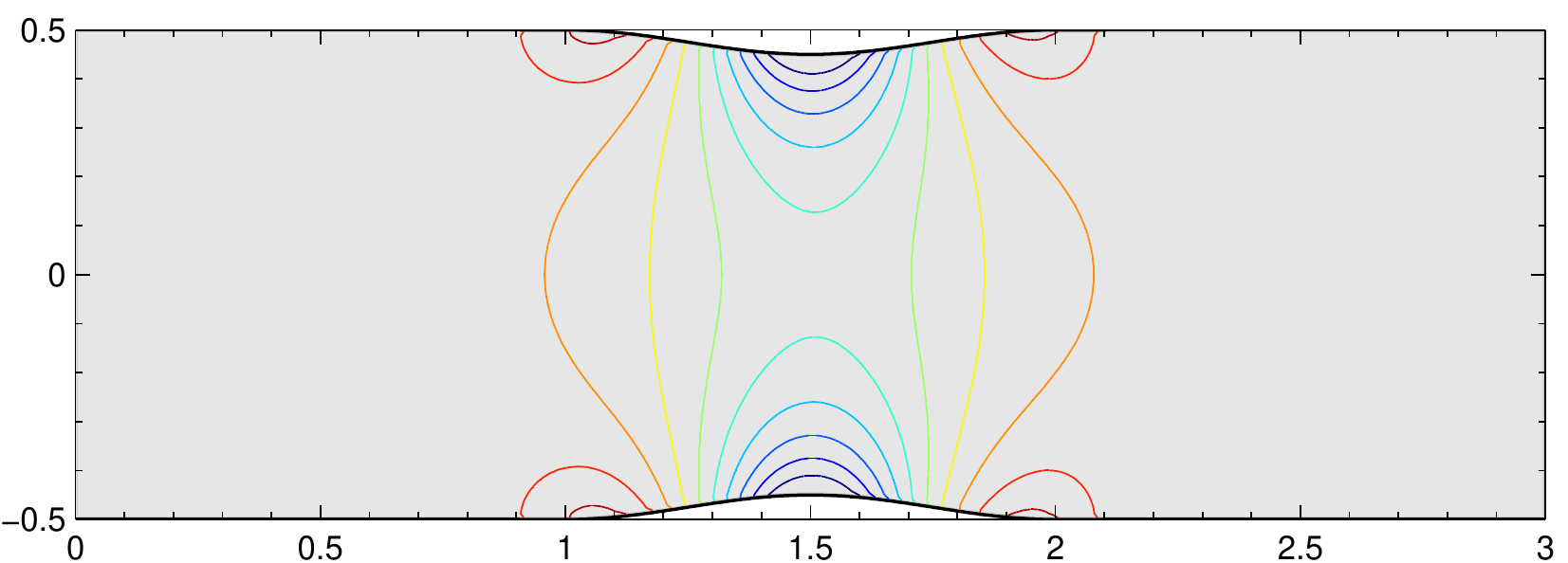}
  }
  \subfigure[$F_{\rm in} = 0.67$]
        {
  \includegraphics[width=0.47\textwidth]{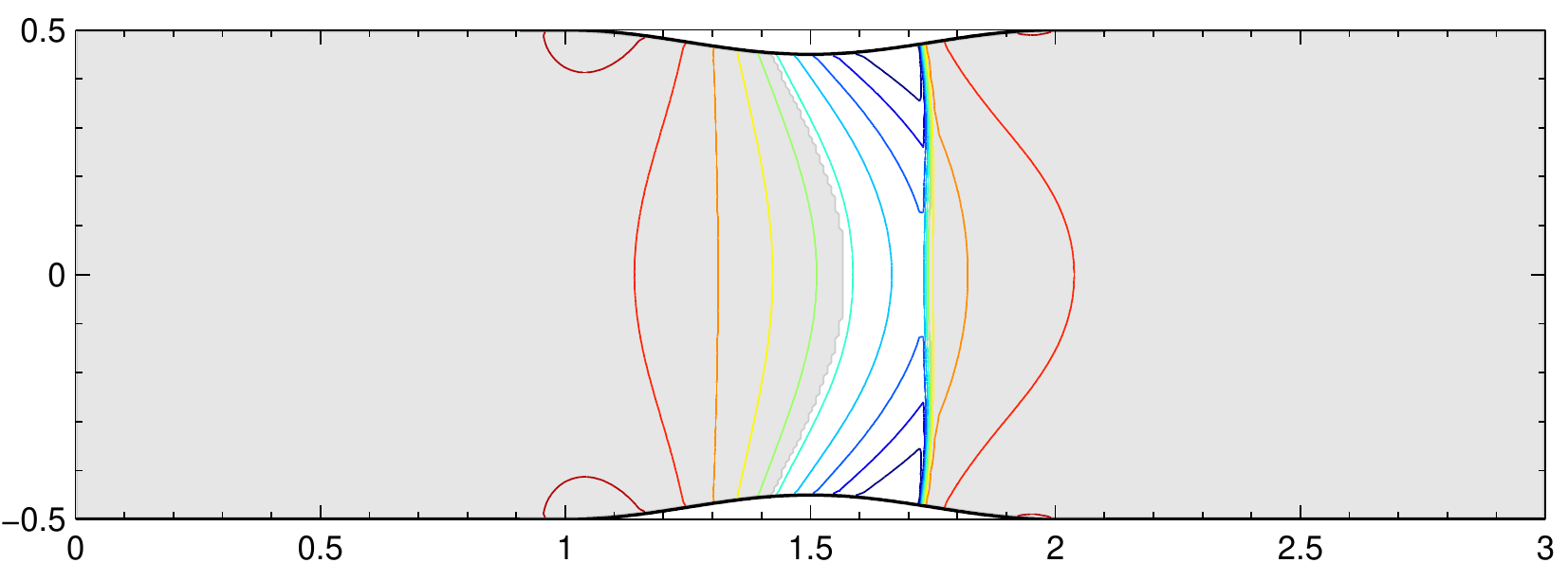}
  }
  \subfigure[$F_{\rm in} = 1.2$]
        {
        \includegraphics[width=0.47\textwidth]{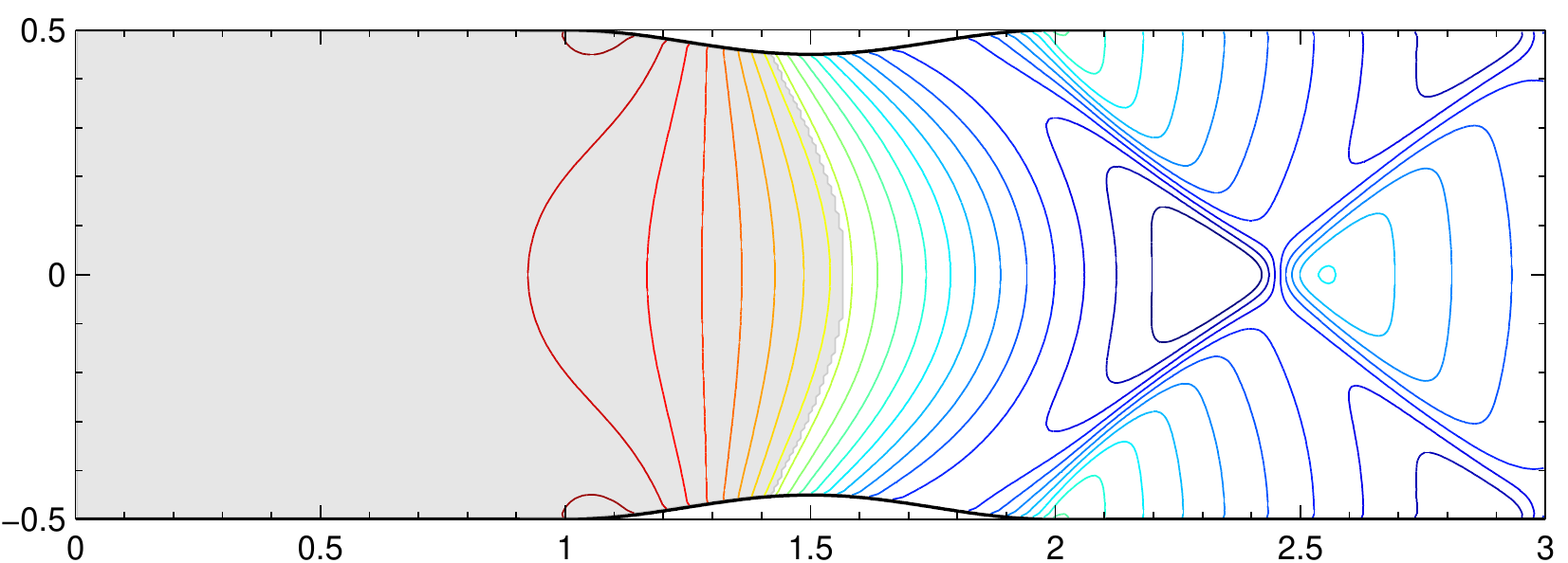}
        }
 \subfigure[$F_{\rm in} =1.7$]
       {
  \includegraphics[width=0.47\textwidth]{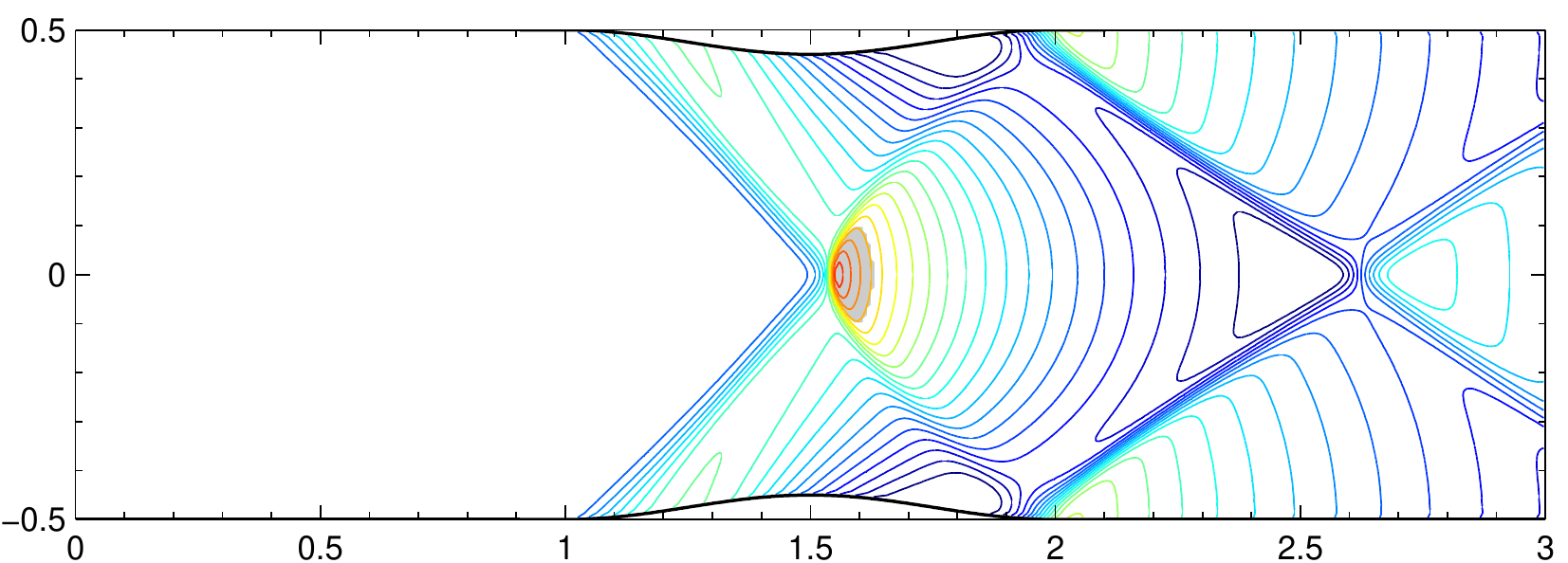}
  }
\subfigure[$F_{\rm in} = 2$]
      {
  \includegraphics[width=0.47\textwidth]{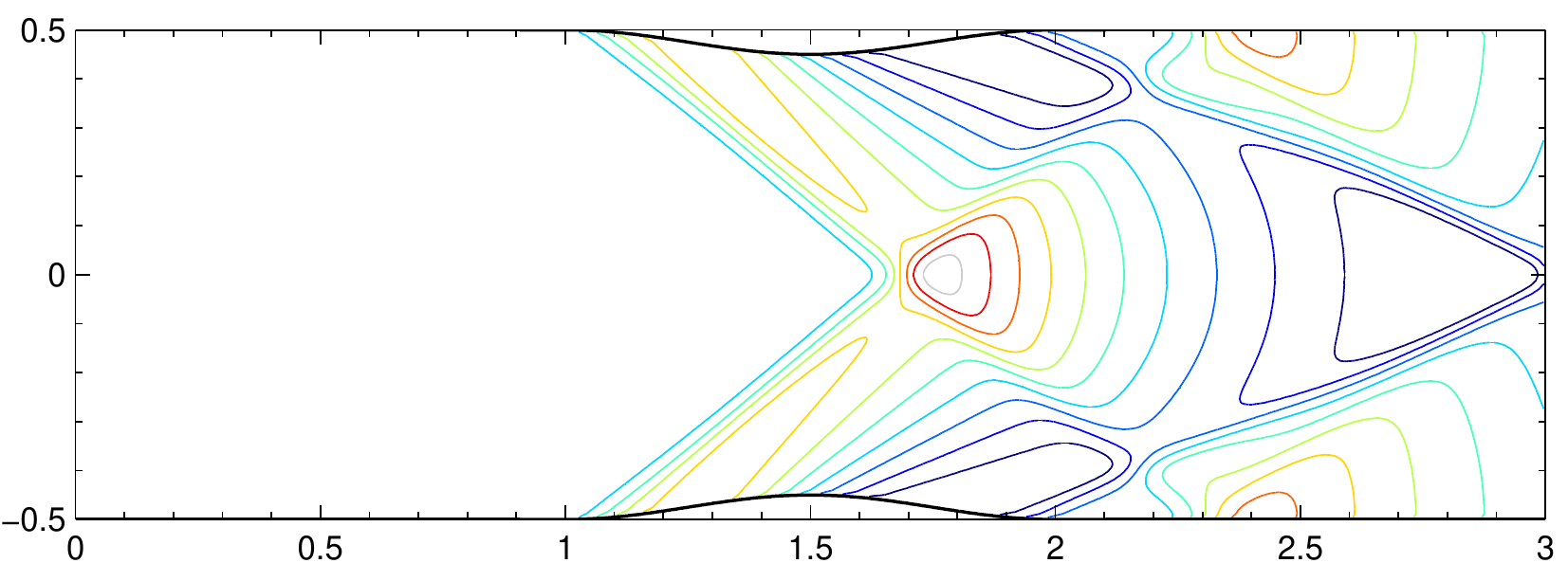}
  }
  \caption{\small Example \ref{example2Dsym2}:
  	The structured mesh of $96 \times 32$ cells and
the contours of steady water depth $h$ with the shaded subcritical regions obtained by {\tt NMGM} on the mesh of $384\times 128$ cells.
}\label{fig:2d_02_solu}
\end{figure}

\begin{example} \label{example2Dsym2}\rm
This example  is to solve another  channel constriction problem \cite{Hubbard2001},
in which a  channel of length $3$ units
is  with a symmetric constriction of length $1$ unit at its center $x=1.5$ and  the variable width
%breadth=¿í¶È
\begin{equation*}
W(x)=
\begin{cases}
1 - (1-W_{\rm min}) \cos^2(\pi (x-1.5 )) , & |x-1.5| \le 0.5,\\
1,&    \rm{otherwise},\end{cases}
\end{equation*}
where $W_{\rm min}$ is the minimum channel width,
%and $x$ is the distance into the channel (so the throat is positioned at the midpoint of the
%constriction),
see Fig. \ref{fig:2d_02_solu}(a)  for its geometry and corresponding  structured mesh of $96 \times 32$ cells.
%The steady solutions of this example have been studied in %\cite{Hubbard2001}
%minutely with different boundary conditions.
Boundary conditions are specified as follows:
\begin{align*}
& hu = h_0F_{\rm in}\sqrt{gh_0} ,\ hv=0  \mbox{  for subcritical inflow},\\
& h=h_0   \mbox{  for   subcritical outflow},\\
& hu = h_0F_{\rm in}\sqrt{gh_0} , \ hv=0,\ h=h_0
 \mbox{ for a supercritical inflow},\\
&  \mbox{supercritical outflow boundary condition at $x=3$},\\
& \mbox{slip boundary conditions on the top and bottom boundaries},
\end{align*}
%$hu = h_0F_{\rm in}\sqrt{gh_0} , hv=0, h=h_0$
%at a supercritical inflow boundary  and outflow boundary condition
%is taken at $x=3$; and slip boundary condition in the $y$--direction.
Our computations take $h_0= { 1,}$ $W_{\rm min}=0.9$, and $F_{\rm in} = 0.5, 0.67, 1.2, 1.7$, and 2,  respectively,
 investigate the effect of the Froude number $F_{\rm in}$  on the convergence of {\tt NMGM}, and demonstrate the robustness of {\tt NMGM}.
Figs. \ref{fig:2d_02_solu}(b-f) {display} the contours of steady water depth $h$ obtained by {\tt NMGM}
on the mesh of $384\times 128$ cells,
where the subcritical region (in which the Froude number $F_r(\vec x)<1$) has been shaded.
Those plots show that the steady solutions resolved by {\tt NMGM}
are non-oscillatory and the symmetry of the flow is well preserved.
Those results are very similar to those in \cite{Hubbard2001}, and
five different types of steady flow pattern can be observed as follows.
\begin{itemize}
  \item {\tt Case 1} ($F_{\rm in} = 0.5$): The  flow is smooth and purely subcritical.
  \item {\tt Case 2} ($F_{\rm in} = 0.67$): The flow is subcritical at inflow and outflow,
critical at the channel throat, and  with a steady discontinuity in the {divergent} region of the channel.
  \item {\tt Case 3} ($F_{\rm in} = 1.2$): The   flow shows the cross-wave pattern with the oblique downstream jumps,
and is subcritical at inflow, critical at the throat, and supercritical at outflow.
  \item {\tt Case 4} ($F_{\rm in} = 1.7$): The flow is supercritical at inflow and outflow, and with oblique jumps and a
subcritical pocket within the constriction.
  \item {\tt Case 5} ($F_{\rm in} = 2$): The   flow is purely supercritical everywhere throughout the channel.
\end{itemize}

The convergence history of {\tt NMGM} in terms of the  {\tt NMGM} iteration number  $N_{\rm step}$  on three  meshes
are presented in Figs. \ref{fig:2d_02_Re}(a)-(e).
%And in these computations, HLLC flux is used, the  number of levels of
%the multigrid is set to be $4$, and the V-cycle is adopted.
It is seen that {\tt NMGM} exhibits good robustness and
works well on those  meshes for the above steady flows.
 {\tt Case 2}  requires more iteration steps than the other cases.
Fig. \ref{fig:2d_02_Re}(f) gives
histogram of the {\tt NMGM} iteration number $N_{\rm step}$
of {\tt NMGM} in terms of the Froude number $F_{\rm in}$ on the mesh
of $384 \times 128 $ cells.
We further see that  {\tt NMGM}
requires more iteration steps when
the Froude number $F_{\rm in}$ is less than 0.3
or equal to 0.7.
%
%the low-speed steady
% subcritical flow and the steady transcritical flow
% which is subcritical at inflow and outflow
%are relatively difficult for {\tt NMGM}.

%In  computations,  $\epsilon_p=0.2$.

\begin{figure}[htbp]
  \centering
\subfigure[$F_{\rm in}=0.5$]
      {
  \includegraphics[width=0.46\textwidth,height=5cm]{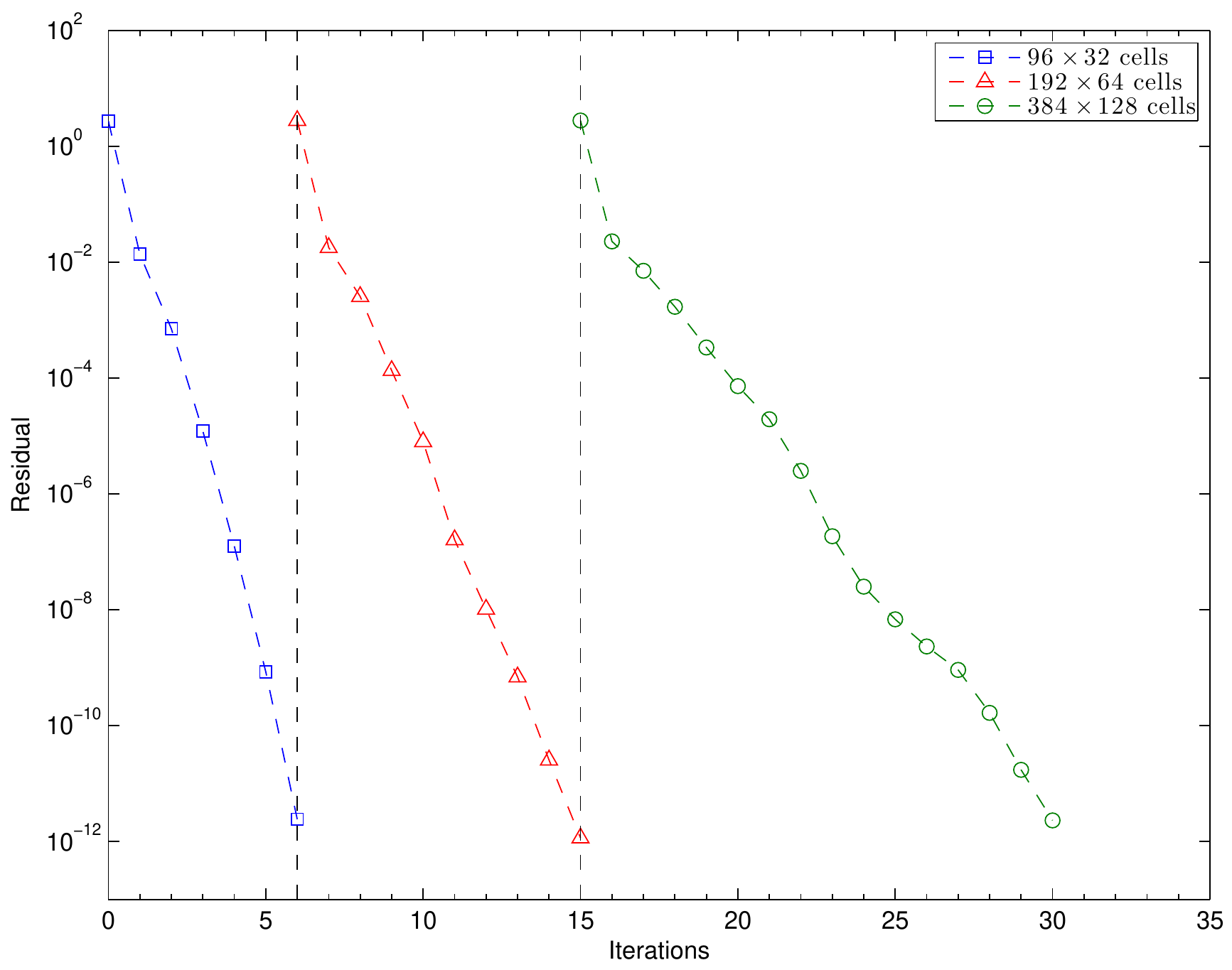}}
\subfigure[$F_{\rm in}=0.67$]
      {
  \includegraphics[width=0.46\textwidth,height=5cm]{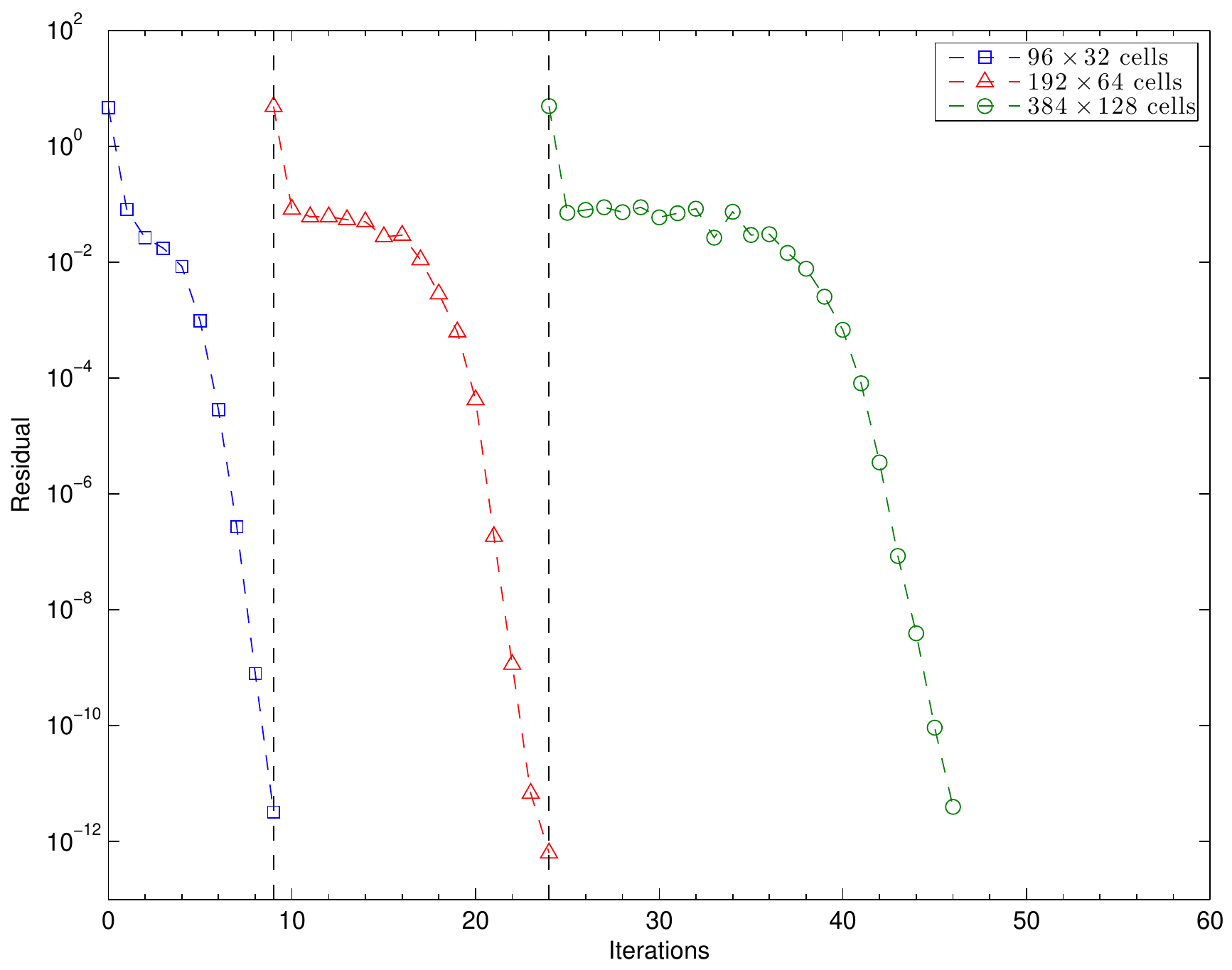}}
\subfigure[$F_{\rm in}=1.2$]
      {
  \includegraphics[width=0.46\textwidth,height=5cm]{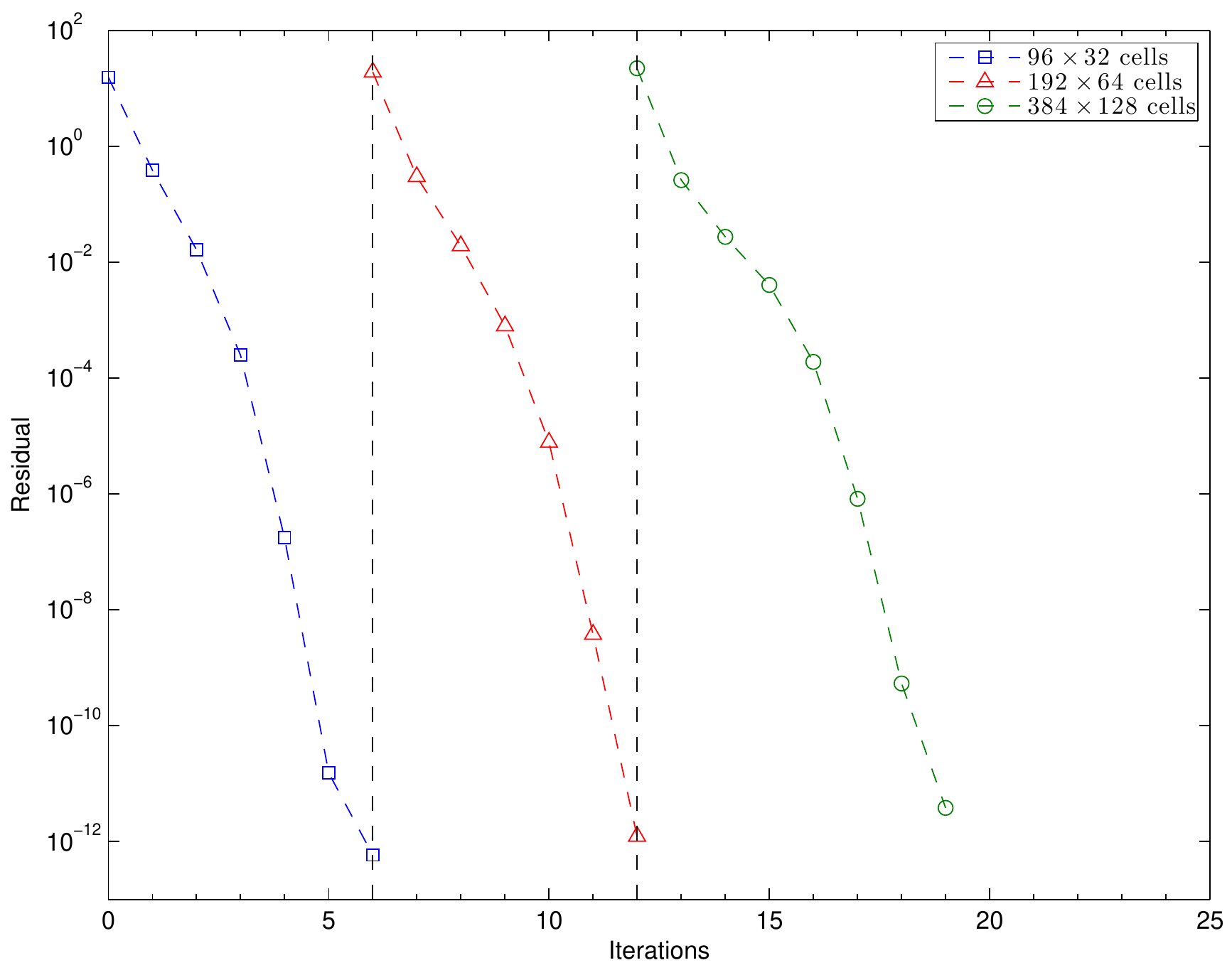}}
\subfigure[$F_{\rm in}=1.7$]
      {
  \includegraphics[width=0.46\textwidth,height=5cm]{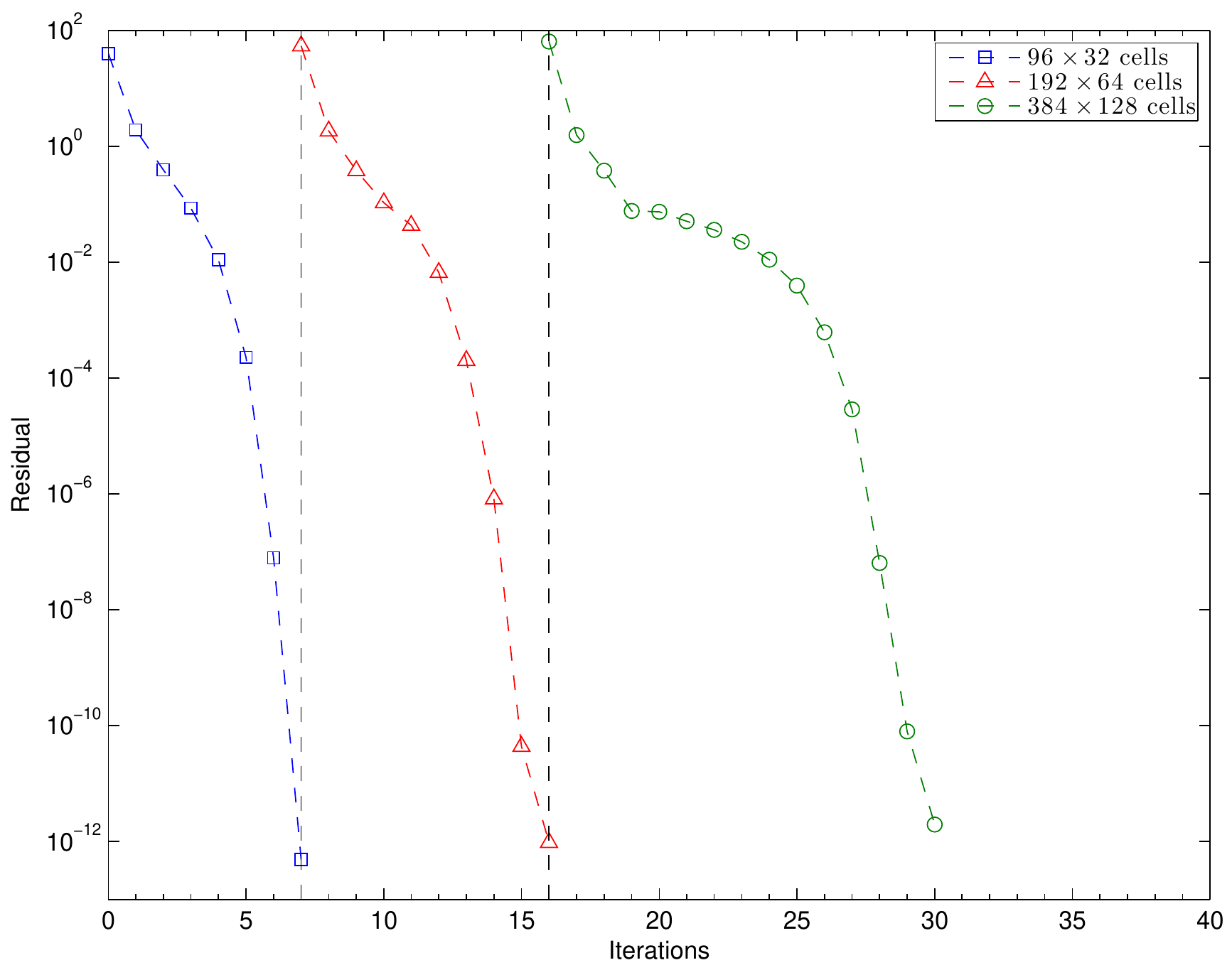}}
\subfigure[$F_{\rm in}=2$]
      {
\includegraphics[width=0.46\textwidth,height=5cm]{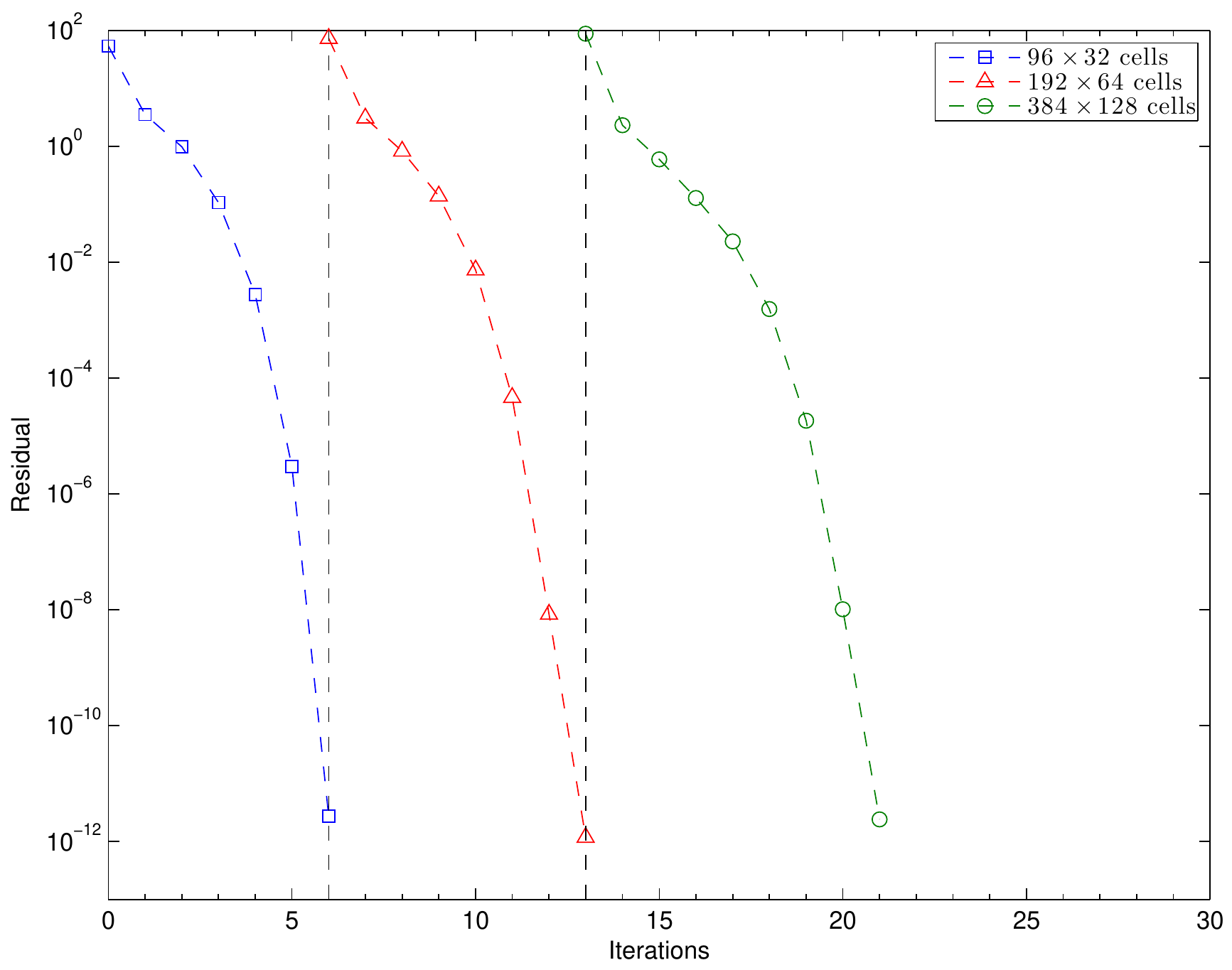}}
\subfigure[Histogram of
 $N_{\rm step}$ versus $F_{\rm in}$]
      {
\includegraphics[width=0.46\textwidth,height=5cm]{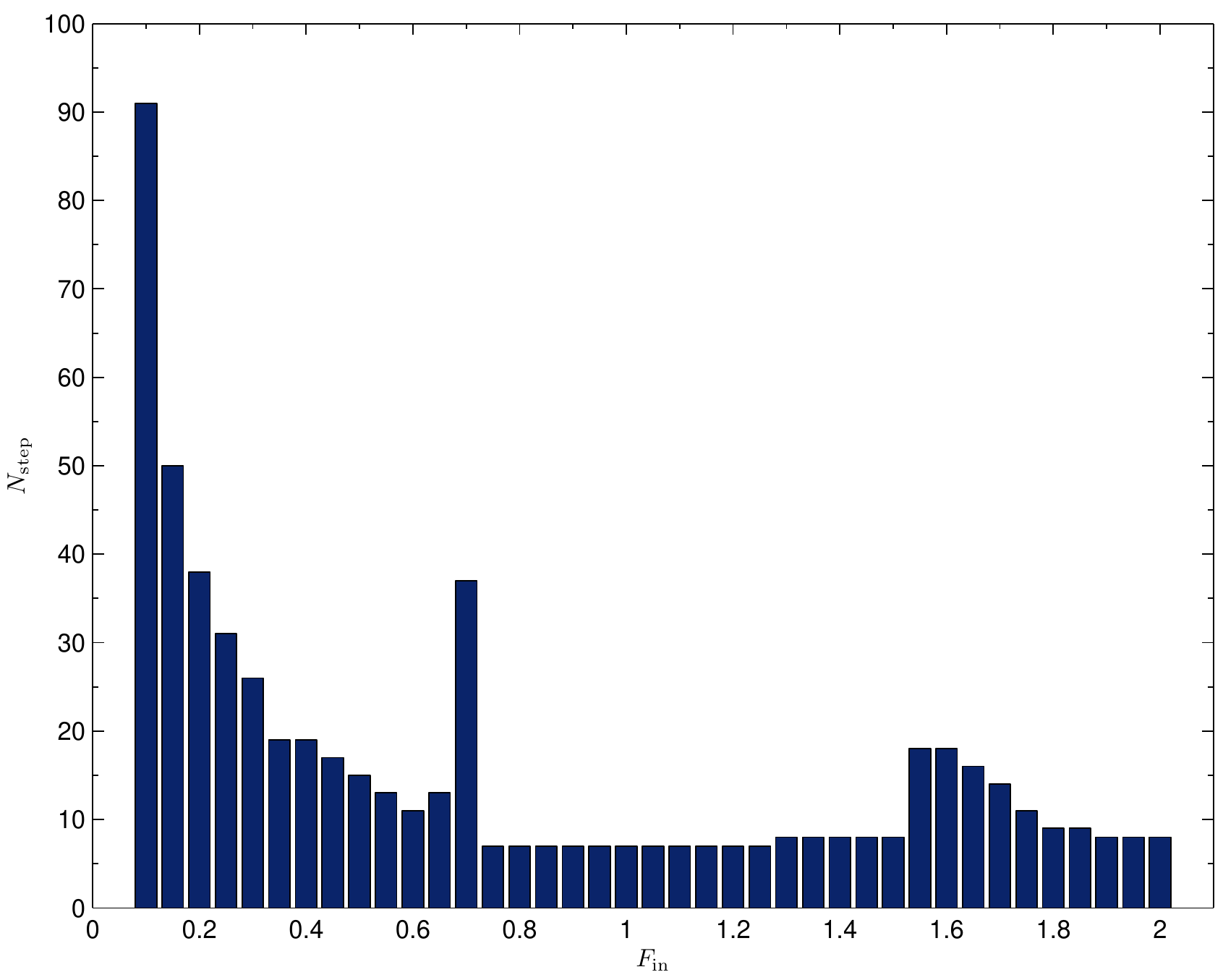}}
  \caption{\small Example \ref{example2Dsym2}: (a)-(e).
Convergence history in terms of the  {\tt NMGM} iteration number  $N_{\rm step}$  on three  meshes;
(f). Histogram of
the {\tt NMGM} iteration number $N_{\rm step}$  on the mesh
of $384 \times 128 $ cells versus  $F_{\rm in}$ uniformly varying from 0.1 to 2.
}\label{fig:2d_02_Re}
\end{figure}

\end{example}

\begin{figure}[htbp]
  \centering
  \includegraphics[width=0.65\textwidth]{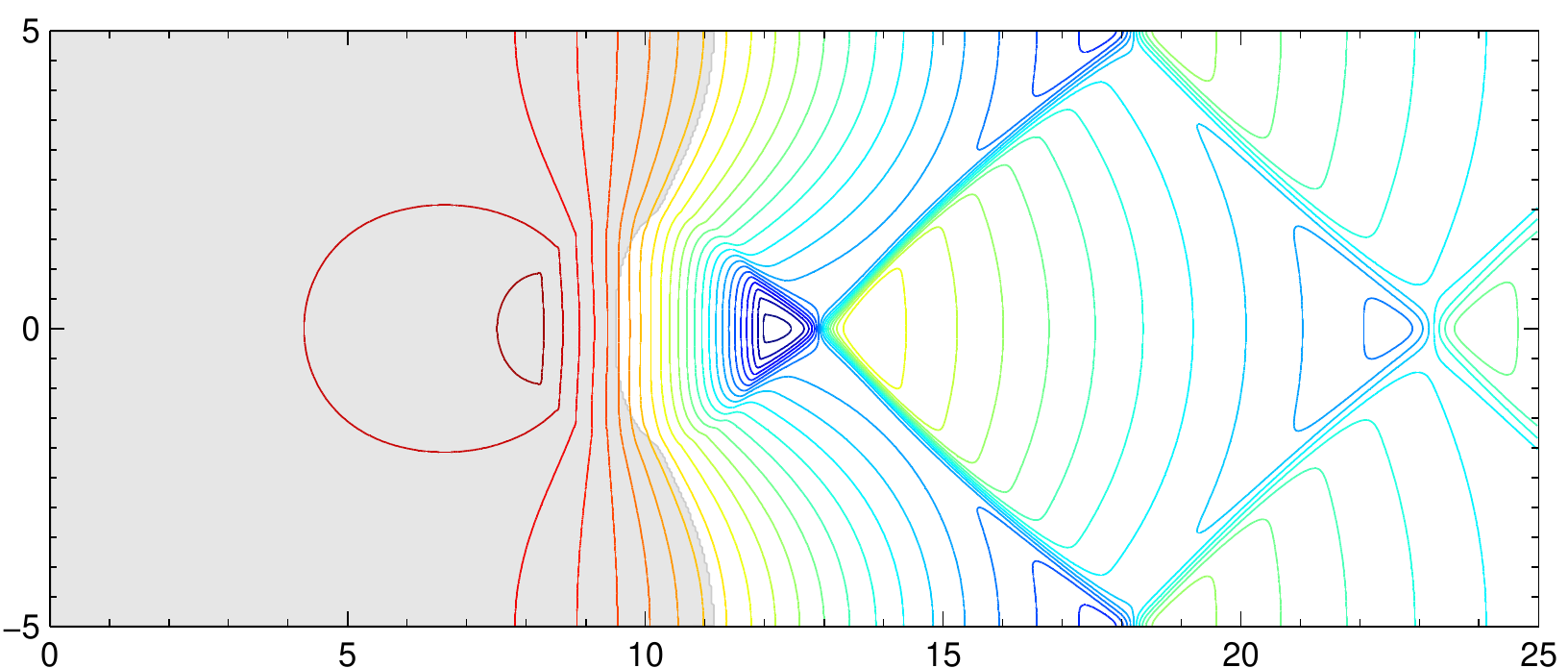}
  \caption{\small Example \ref{example2Dtran}:
The contours of the steady water surface $h+z$ with the shaded subcritical region
 obtained by {\tt NMGM} on the mesh of $640\times 320$ cells.  25 equally spaced
contour lines are used.
}\label{fig:2d_04_solu}
\end{figure}

\begin{figure}[htbp]
  \centering
  \includegraphics[width=0.46\textwidth,height=5cm]{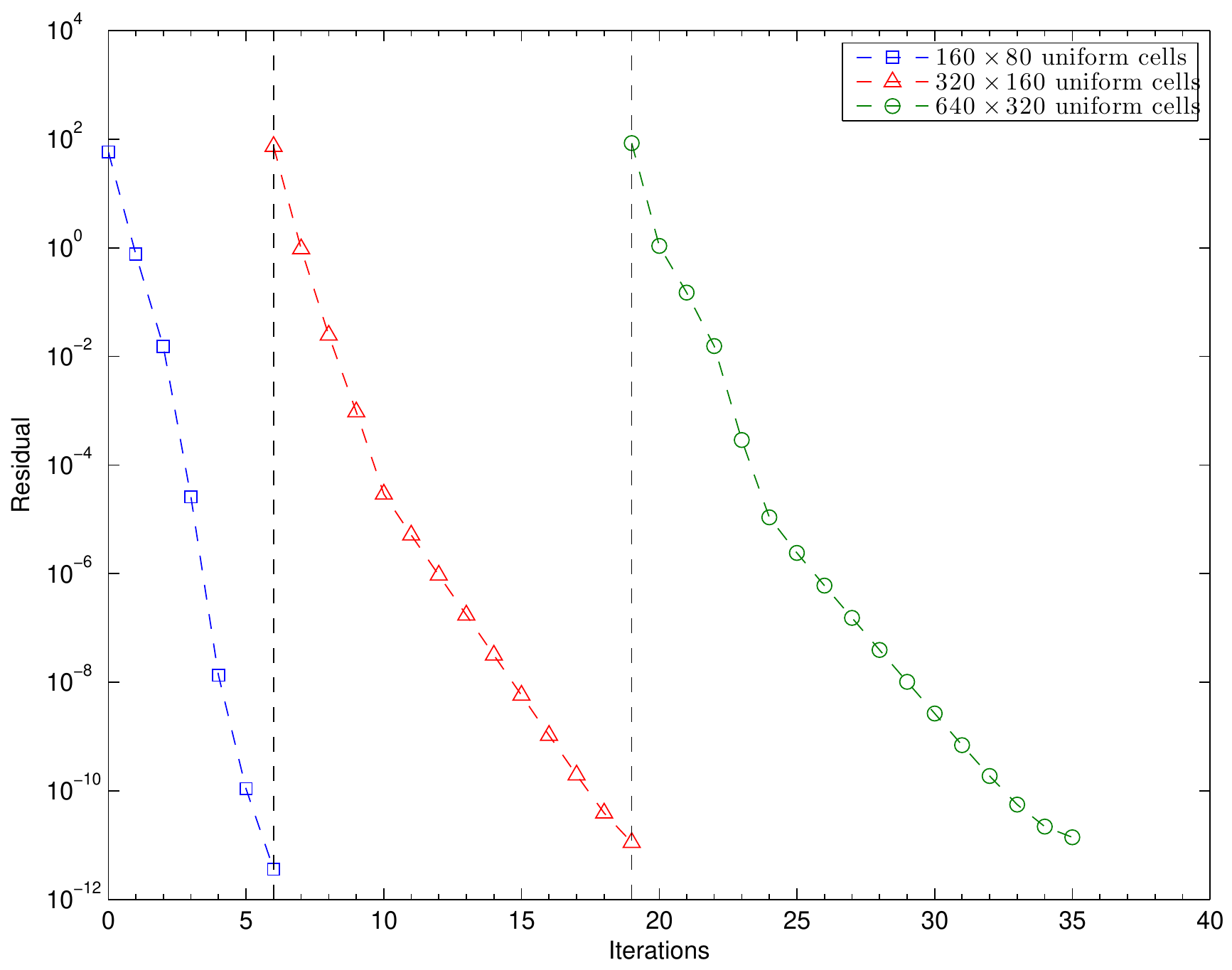}
  \includegraphics[width=0.46\textwidth,height=5cm]{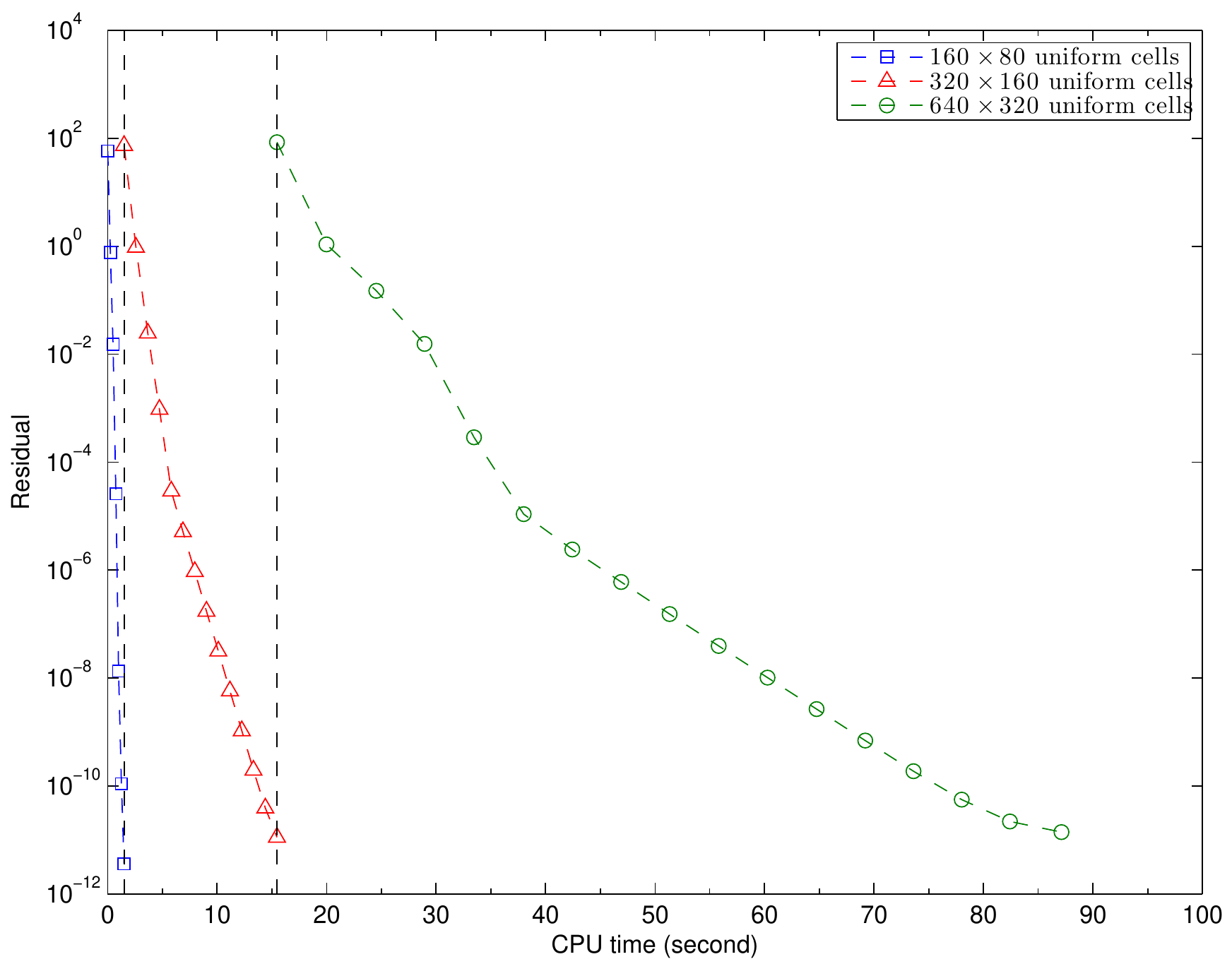}
  \caption{\small Example \ref{example2Dtran}:
Convergence history in terms of the {\tt NMGM} iteration number $N_{\rm step}$ (left) and CPU time (right) on three
  meshes.
}\label{fig:2d_04_Re}
\end{figure}

%%%%%%%%%%%%%%%%%%%%%%%%%%%%%%%%%% Example 2D 04 %%%%%%%%%%%%%%%%%%%%%%%%%%%%%%%%%%%%%%

\begin{example} \label{example2Dtran}\rm
The test  describes a transcritical flow over a 2D bump, which is a
2D extension of  Example \ref{example1Dtrans}.
The bottom topography is defined by
\begin{equation*}
z(x,y)=
\begin{cases}
0.2 - 0.05\left((x-10)^2+y^2 \right), &   (x-10)^2 +y^2 <4,\\
0,&    \rm{otherwise},\end{cases}
\end{equation*}
in the channel of $[0, 25]\times [-5,5]$.
The discharge $hu = 1.53, hv =0$ is imposed at $x = 0$, the water height $h = 0.52$ is imposed at
the downstream $x = 25$, and the reflective boundaries are specified at $y = \pm 5$ in the $y$--direction.

The {contours} of the steady water surface $h+z$ obtained by
{\tt NMGM}
are displayed in Fig. \ref{fig:2d_04_solu} on the mesh of $640 \times 320$ uniform cells, where the subcritical region is shaded.
In the steady flow, the cross-wave pattern with the oblique downstream jumps is
observed. Moreover,
the flow varies from subcritical at inflow to supercritical at outflow, and
is similar to the 1D Example \ref{example1Dtrans}(I).
The efficiency of {\tt NMGM} are demonstrated in Fig. \ref{fig:2d_04_Re}
by displaying its convergence history
 versus iterations and CPU time on three meshes of $ 160 \times 80 $, $320 \times 160$,
and $640 \times 320$ uniform cells, respectively.
In those computations, $\epsilon_p=2$.
%, HLLC flux is used, the  number of levels of
%the multigrid is set to be $4$, and the V-cycle is adopted.

\end{example}

%%%%%%%%%%%%%%%%%%%%%%%%%%%%%%%%%% Example 2D 06 %%%%%%%%%%%%%%%%%%%%%%%%%%%%%%%%%%%%%%

\begin{figure}[htbp]
  \centering
  \includegraphics[width=0.47\textwidth]{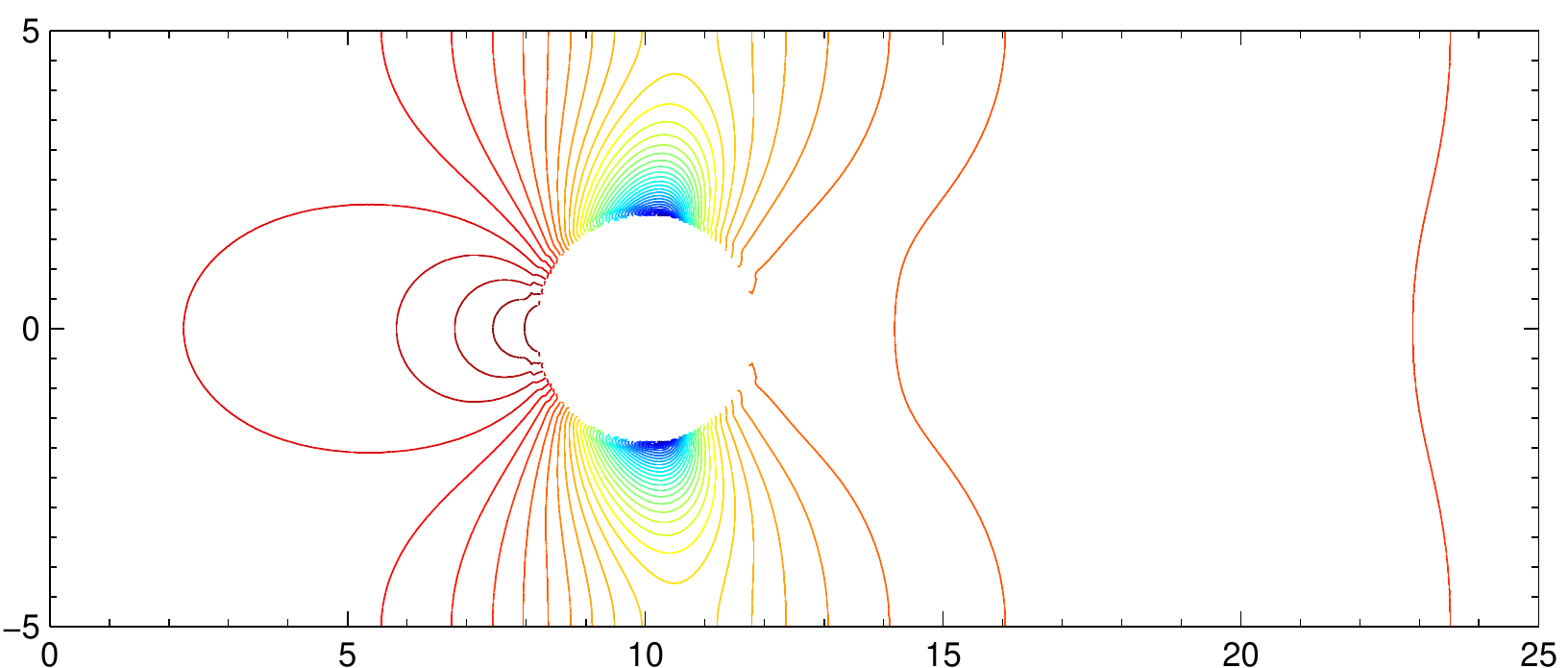}
  \includegraphics[width=0.47\textwidth]{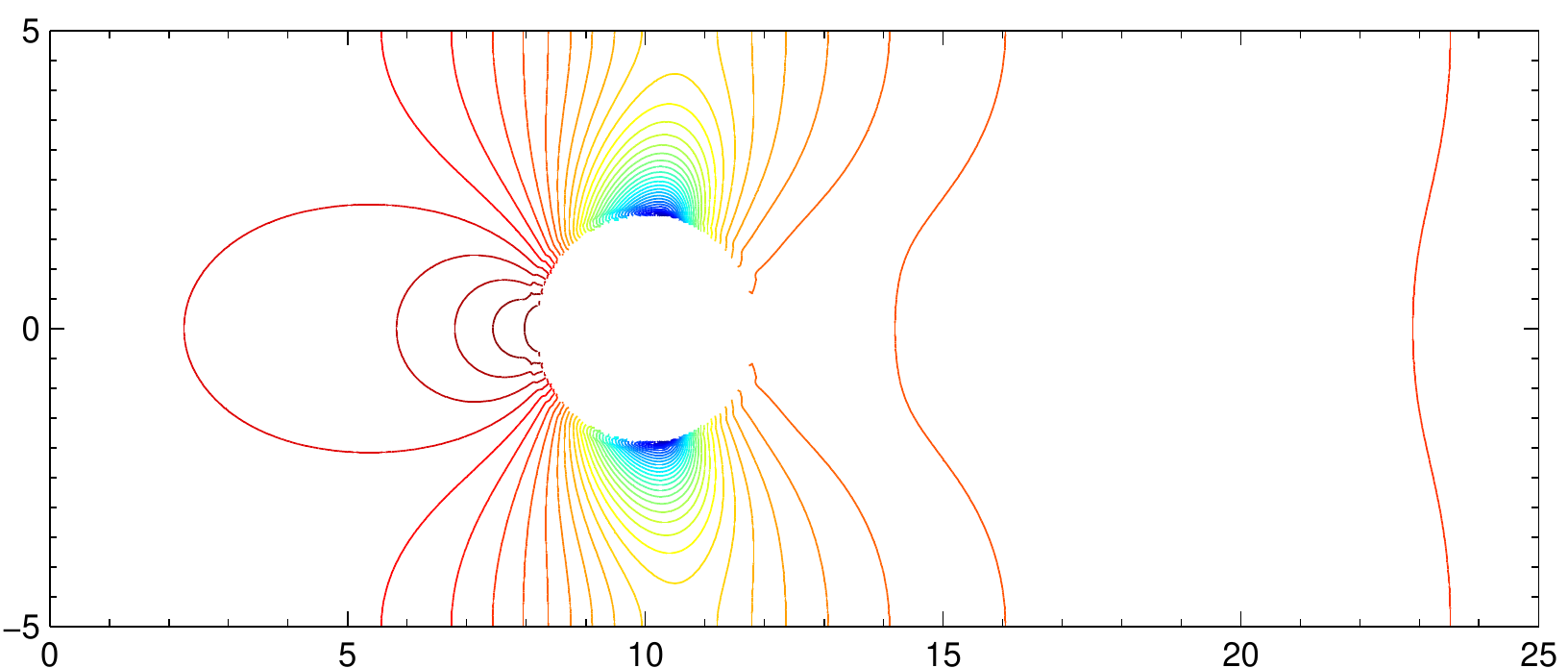}
  \caption{\small Example \ref{example2Ddry}:
The contours of the steady water surface $h+z$ (only the wet region is shown)
respectively obtained by using {\tt NMGM} (left)
and an explicit scheme to solve the time-dependent SWEs (right)
on the mesh of $512\times 256$ cells.  40 equally spaced
contour lines from $0.1067$ to $0.2201$ are used.
}\label{fig:2d_06_solu}
\end{figure}

\begin{figure}[htbp]
  \centering
  \includegraphics[width=0.46\textwidth,height=5cm]{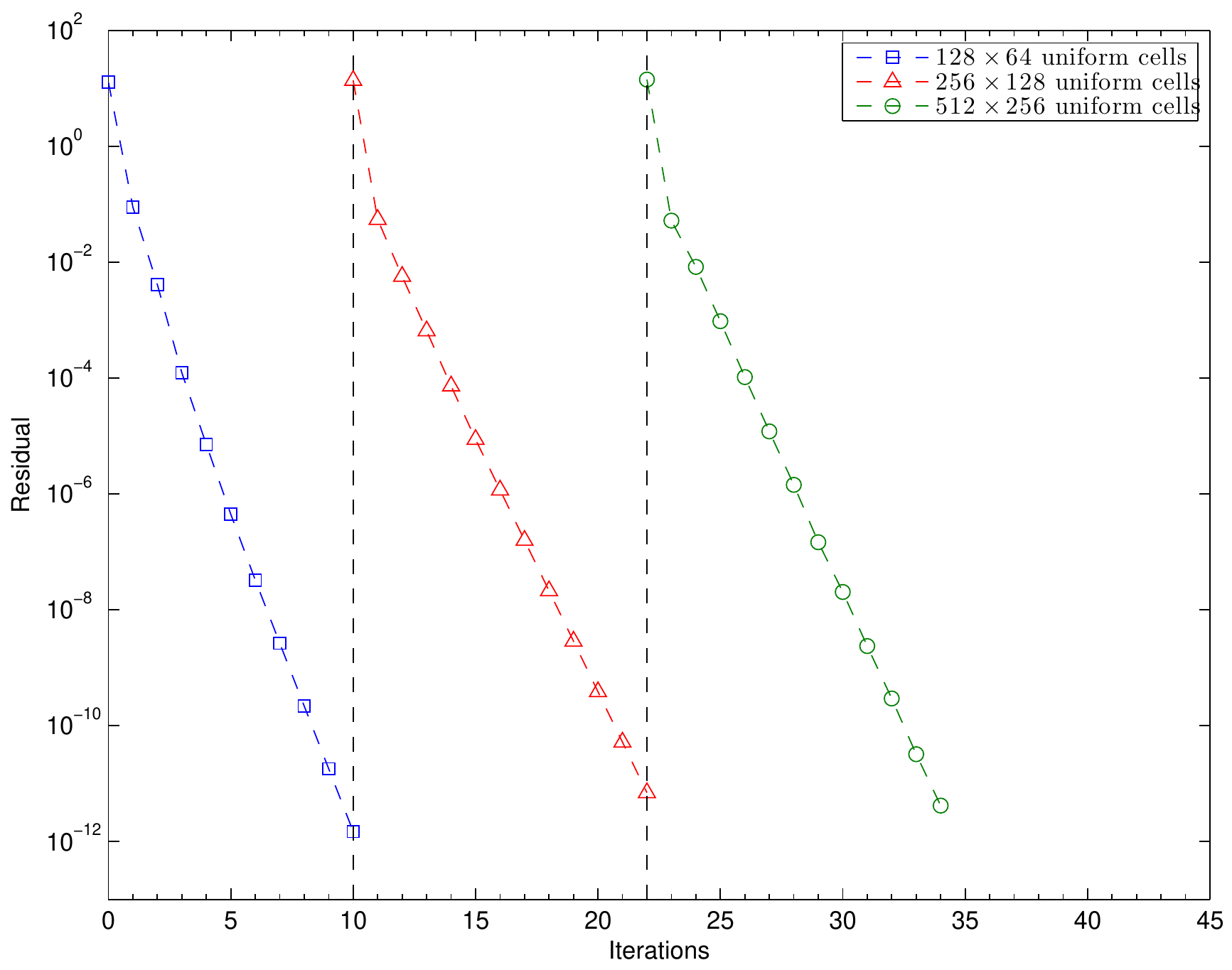}
  \includegraphics[width=0.46\textwidth,height=5cm]{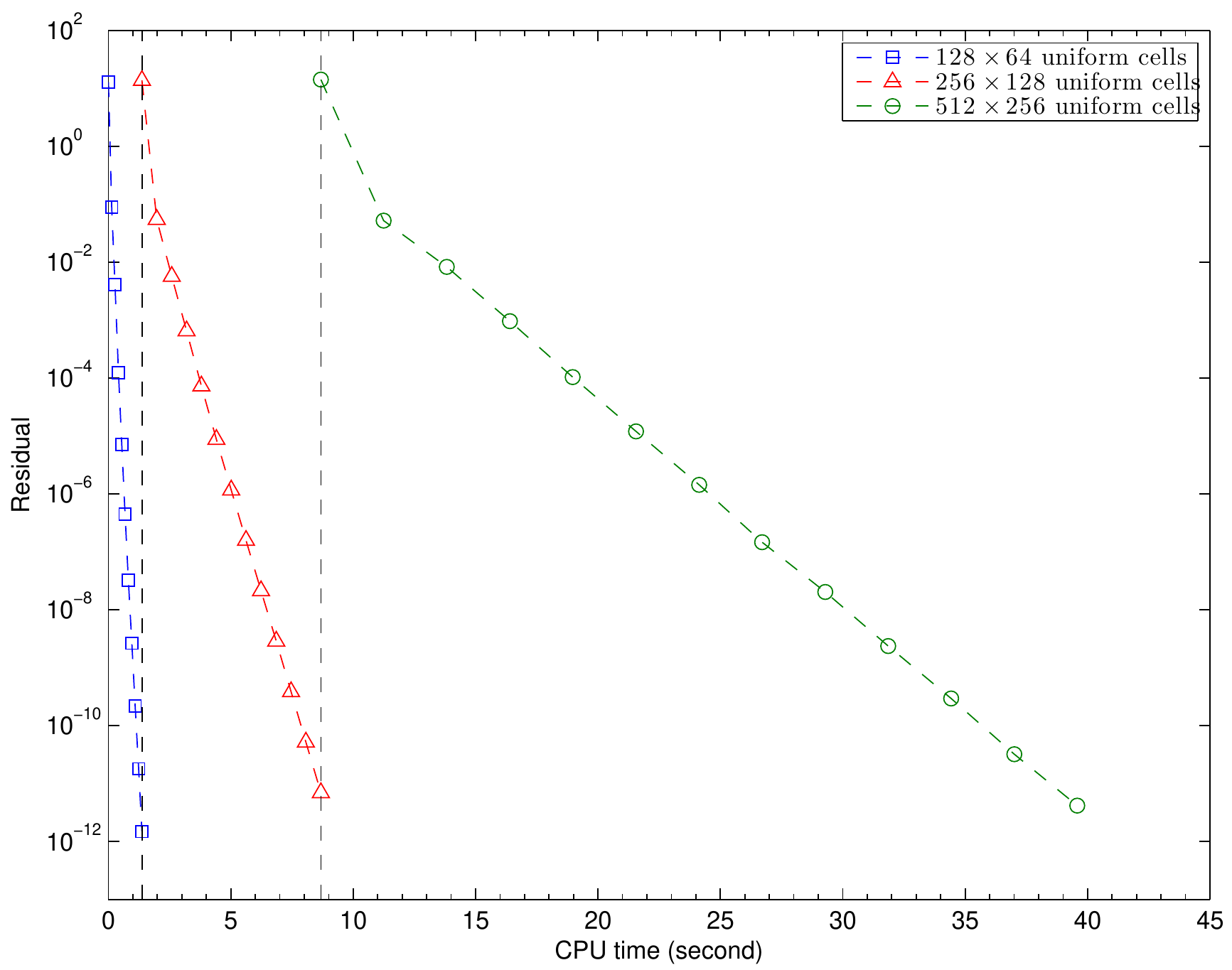}
  \caption{\small Example \ref{example2Ddry}:
Convergence history in terms of the {\tt NMGM} iteration number $N_{\rm step}$ (left) and CPU time (right) on three
  meshes.
}\label{fig:2d_06_Re}
\end{figure}

\begin{example} \label{example2Ddry}\rm
The last example simulates a steady channel flow around a hill (a relatively high bump) which
is defined by
\begin{equation*}
z(x,y)=
\begin{cases}
1.2 - 0.3\left((x-10)^2+y^2 \right), &   (x-10)^2 +y^2 <4,\\
0,&    \rm{otherwise},\end{cases}
\end{equation*}
in the channel of $[0, 25]\times [-5,5]$ with boundary conditions:
The discharge $hu = 0.1, hv =0$  at $x = 0$,
 the water height $h = 0.2$  at
the downstream $x = 25$, and the reflective boundaries conditions
at $y = \pm 5$ in the $y$--direction.
Initially, the channel is full of static flow with the water height $h(x,y)=\max \{0.2-z(x,y),0\}$.
Since  the steady solution  involves wet/dry transition,
 thus it becomes much more difficult to obtain such the steady solution by solving the time-dependent SWEs with numerical schemes.
{\tt NMGM} and the BLU-SGS iteration with the HLLC or Roe
flux  fail to work.
%In the computations of {\tt NMGM},  $\epsilon_p=0.2$.

Fig. \ref{fig:2d_06_solu} displays the contours of
the steady water surface $h+z$ obtained by
{\tt NMGM} with the LLF flux on the mesh of $512 \times 256$ uniform cells,
in comparison to the result given by using a first order accurate explicit finite volume scheme with
LLF flux to solve the time-dependent SWEs \eqref{eq:GNeqs}
 on the same mesh. The latter takes 412652 time steps
and CPU time of 13804 seconds to reduce the residual
to $1.5 \times 10^{-12}$ at physical time $t=1701.96$.
The results show that the steady solution given by {\tt NMGM}
agree well with that given by the explicit scheme,
and {\tt NMGM} with the LLF flux
works efficiently and robustly for the steady solution with wet/dry transitions.
Fig. \ref{fig:2d_06_Re} gives the convergence history of {\tt NMGM} with LLF flux
 in terms of the {\tt NMGM} iteration number $N_{\rm step}$ and CPU time on three   meshes of
$128 \times 64$, $256 \times 128$, and $512 \times 256$ uniform cells, respectively.
It is seen that the convergence behaviors of {\tt NMGM} are similar
on those  meshes, and $N_{\rm step}$  scarcely changes with the
mesh refinement.

%, the  number of levels of
%the multigrid is set to be $4$, and the W-cycle is adopted.

\end{example}

\section{Conclusions}\label{sec:conclud}

The paper developed a Newton multigrid method  for 1D and 2D steady shallow water equations ({SWEs}) with topography and dry areas.
The steady-state SWEs were first
approximated by using the well-balanced finite volume discretization
based on the hydrostatic reconstruction technique.
The resulting  nonlinear system was linearized  by using Newton's method, and the geometric MG method  with the  block symmetric Gauss-Seidel (SGS) smoother was used to solve the linear system.
The proposed {Newton MG}  method made use of the local residual to regularize the Jacobian matrix of the Newton iteration,  and could  handle the steady-state problem with wet/dry transitions.
Several numerical experiments were conducted to demonstrate the
efficiency, robustness, and  well-balanced property  of the proposed  method. Moreover,
the relation between the convergence behavior of  the proposed method
with the HLL, LLF, or Roe flux and the distribution
of the eigenvalues of the iteration matrix was detailedly discussed, in comparison to {the block LU-SGS method.}
Numerical results showed that
%{\tt NMGM} could preserve the well-balanced property
%and handle the steady state problem with wet/dry transition,
 convergence behaviors of {\tt NMGM} depended on the numerical flux
 and  the Froude number of the flow.
% different convergence behaviors were obtained for
%subcritical flow, transcritical flow and supercritical flow.

\section*{Acknowledgements}

The work was partially supported by
the National Natural Science Foundation
of China (Nos. 91330205  \& 11421101).

\end{document}